\documentclass[a4paper,11pt,reqno]{amsart}

\usepackage[margin=1in,includehead,includefoot,asymmetric,marginparwidth=1cm]{geometry}
\usepackage[utf8]{inputenc}
\usepackage[T1]{fontenc}
\usepackage{lmodern,microtype}
\usepackage{hyperref}
\usepackage{url}
\usepackage{booktabs}
\usepackage{mathtools,amssymb,amsthm}
\usepackage[foot]{amsaddr}
\usepackage{graphicx}
\usepackage[dvipsnames]{xcolor}
\usepackage[margin=5mm]{caption}
\usepackage{subcaption}
\usepackage{enumitem}
\usepackage{cleveref}
\usepackage[b]{esvect} 
\usepackage[textsize=footnotesize,textwidth=22mm,disable]{todonotes}
\usepackage{stackengine,scalerel}
\usepackage{subdepth} 
\usepackage{mycommands}

\graphicspath{{./graphics/}}

\newtheorem{theorem}{Theorem}
\newtheorem{lemma}[theorem]{Lemma}
\newtheorem{conjecture}[theorem]{Conjecture}
\newtheorem{proposition}[theorem]{Proposition}

\crefname{figure}{Figure}{Figures}
\crefname{theorem}{Theorem}{Theorems}
\crefname{lemma}{Lemma}{Lemmas}
\crefname{corollary}{Corollary}{Corollaries}
\crefname{section}{Section}{Sections}
\crefname{appendix}{Appendix}{Appendices}
\crefname{remark}{Remark}{Remarks}
\crefname{claim}{Claim}{Claims}
\crefname{conjecture}{Conjecture}{Conjectures}
\crefname{observation}{Observation}{Observations}

\definecolor{amethyst}{rgb}{0.7,0.5,0.9}
\newcommand{\TM}[1]{\todo[color=red]{TM: #1}}%


\hypersetup{
  pdftitle={Listing faces of polytopes},
  pdfauthor={Nastaran Behrooznia, Sofia Brenner, Arturo Merino, Torsten M\"utze, Christian Rieck, Francesco Verciani}
}

\title{Listing faces of polytopes}

\author{Nastaran Behrooznia}
\address[Nastaran Behrooznia]{Department of Computer Science, University of Warwick, United Kingdom}
\email{nastaran.behrooznia@warwick.ac.uk}

\author{Sofia Brenner}
\address[Sofia Brenner]{Institut f\"ur Mathematik, Universit\"at Kassel, Germany}
\email{sofia.brenner@mathematik.uni-kassel.de}

\author{Arturo Merino}
\address[Arturo Merino]{Universidad de O'Higgins, Rancagua, Chile}
\email{arturo.merino@uoh.cl}

\author{Torsten M\"utze}
\address[Torsten M\"utze]{Institut f\"ur Mathematik, Universit\"at Kassel, Germany}
\email{tmuetze@mathematik.uni-kassel.de}

\author{Christian Rieck}
\address[Christian Rieck]{Institut f\"ur Mathematik, Universit\"at Kassel, Germany}
\email{christian.rieck@mathematik.uni-kassel.de}

\author{Francesco Verciani}
\address[Francesco Verciani]{Institut f\"ur Mathematik, Universit\"at Kassel, Germany}
\email{francesco.verciani@mathematik.uni-kassel.de}

\thanks{This project was supported by German Science Foundation grant~522790373.
Arturo Merino was supported by ANID FONDECYT Iniciaci\'on No.~11251528.}

\thanks{An extended abstract of this paper appeared in SODA~2026~\cite{MR5027189}.}

\begin{document}

\begin{abstract}
This paper investigates the problem of listing faces of polytopes that represent combinatorial objects, such as hypercubes, permutahedra, associahedra, and their generalizations.
Firstly, we consider the face lattice, which is the inclusion order of all faces of a polytope, and we seek a Hamiltonian cycle in its cover graph, i.e., for any two consecutive faces, one must be a subface of the other, and their dimensions differ by $1$.
We construct such Hamiltonian cycles for hypercubes, permutahedra, $B$-permutahedra, associahedra, cyclic polytopes, 3-dimensional polytopes, graph associahedra of chordal graphs, and quotientopes.
Secondly, we consider facet-Hamiltonian cycles, which are cycles on the (1-)skeleton of a polytope that enter and leave every facet exactly once.
This notion was recently introduced by Akitaya, Cardinal, Felsner, Kleist, and Lauff~[SODA~2025], where the authors conjectured that $B$-permutahedra admit a facet-Hamiltonian cycle for all dimensions.
We construct such facet-Hamiltonian cycles in this paper, thus establishing their conjecture as a theorem.
A key tool we use are so-called rhombic strips, which are planar spanning subgraphs of the cover graph of the face lattice in which every face is a 4-cycle.
Specifically, we construct a rhombic strip in the face lattice of the hypercube of any dimension, and characterize the existence of rhombic strips in the face lattice of 3-dimensional polytopes.
Our constructions yield time- and space-efficient algorithms for computing the aforementioned cycles and thus for listing the corresponding combinatorial objects, including ordered set partitions and dissections of a convex polygon.
\end{abstract}

\maketitle

\section{Introduction}

\subsection{Combinatorial polytopes}

\begin{figure}[h!]
\makebox[0cm]{ 
\includegraphics{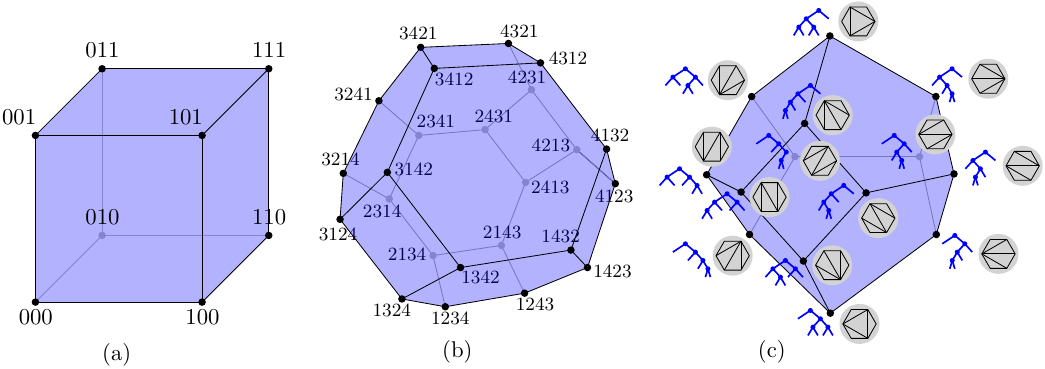}
}
\caption{Examples of 3-dimensional combinatorial polytopes: (a) the hypercube; (b) the permutahedron; (c) the associahedron.}
\label{fig:3poly}
\end{figure}

Convex polytopes are objects of fundamental interest that tie together geometric, combinatorial, algebraic and algorithmic concepts and problems.
In this paper, we are particularly interested in polytopes whose vertices represent combinatorial objects, such as the set of all bitstrings of length~$n$, the set of all permutations of~$[n]\ass\{1,\ldots,n\}$, or the set of all binary trees with $n$ vertices.
The corresponding polytopes are the well-known \defi{hypercube}, whose edges connect pairs of bitstrings that differ in a single bit, the \defi{permutahedron}, whose edges connect pairs of permutations that differ in an adjacent transposition, and the \defi{associahedron}, respectively, whose edges connect pairs of binary trees that differ in a tree rotation; see Figure~\ref{fig:3poly}.
An equivalent model of the associahedron is to consider the dual graphs of the binary trees, namely triangulations of a convex polygon, and each tree rotation translates to a flip operation that changes exactly one diagonal in the triangulation.

Recently, there has been an extensive line of work to define and analyze vast classes of polytopes that generalize the aforementioned three special polytopes.
One of these generalizations are \defi{graph associahedra} \cite{MR2239078,MR2487491,MR2479448}, which are parameterized by an underlying graph~$H$ and have as vertices all elimination trees of~$H$, with edges connecting pairs of elimination trees that differ in a tree rotation (if $H$ is a perfect matching, a complete graph, or a path, respectively, we obtain the hypercube, permutahedron, and associahedron as special cases).
Another important generalization are \defi{quotientopes}~\cite{MR3964495,MR4584712}, which arise from lattice congruences of the weak order on permutations~\cite{MR2142177}.

\subsection{Shortest paths on the skeleton}

The \defi{(1-)skeleton} of a polytope~$P$ is the graph~$G(P)$ formed by the vertices and edges of~$P$.
If $P$ is $d$-dimensional, the graph~$G(P)$ is known to be $d$-connected by Balinski's theorem.
The skeleton captures much information relevant for various applications.

A fundamental algorithmic problem is to compute a shortest path on the skeleton of a polytope, for a given pair of input vertices.
For the associahedron, we are given two binary trees and ask for a shortest sequence of rotations to transform one tree into the other (this is relevant in the context of balancing binary search trees).
This long-standing problem was shown to be NP-hard in a recent breakthrough by Dorfer~\cite{dorfer_preprint}.
Nonetheless, a 2-approximation algorithm~\cite{MR2740180} and various fixed-parameter algorithms are known~\cite{MR2541971,MR4880181}.
Recently, Cunha, Sau, Souza and Valencia-Pabon~\cite{DBLP:conf/icalp/CunhaSSV25} proved that the shortest path problem is fixed-parameter tractable even on graph associahedra, a problem that generalizes the rotation distance between binary trees.
Complementing this result, Cardinal, Pournin and Valencia-Pabon~\cite{MR4668301} found a polynomial-time algorithm for computing rotation distance in graph associahedra of complete split graph, and this is the only known non-trivial family of graphs for which such a positive result is known.
Furthermore, Cardinal and Steiner~\cite{MR4962117} showed that the shortest path problem is NP-hard on base polytopes of polymatroids, even if the polytope is defined by few inequalities.
They also proved inapproximability when the polymatroid is a hypergraphic polytope, whose vertices are in bijection with acyclic orientations of a given hypergraph; see also~\cite{cunha-et-al_preprint}.

Ito, Kakimura, Kamiyama, Kobayashi, and Okamoto~\cite{MR4413061} considered the perfect matching polytope, i.e., the polytope defined as the convex hull of the characteristic vectors of all perfect matchings of a graph.
They proved that shortest paths are NP-hard to compute, even if the underlying graph is planar or bipartite, but the problem can be solved efficiently if the graph is outerplanar.
Cardinal and Steiner~\cite{MR4872027} showed that shortest paths on the perfect matching polytope are also hard to approximate.

The \defi{diameter} of~$G(P)$ is the shortest distance, maximized over all pairs of vertices of~$P$.
This parameter has been heavily studied for various classes of combinatorial polytopes~\cite{MR928904,MR3197650,MR3649601,MR3874284,MR4742547,DBLP:conf/swat/Berendsohn22}, as it provides a lower bound for the running time of the simplex algorithm when optimizing a linear objective function over~$P$.

\begin{figure}
\makebox[0cm]{ 
\includegraphics[page=2]{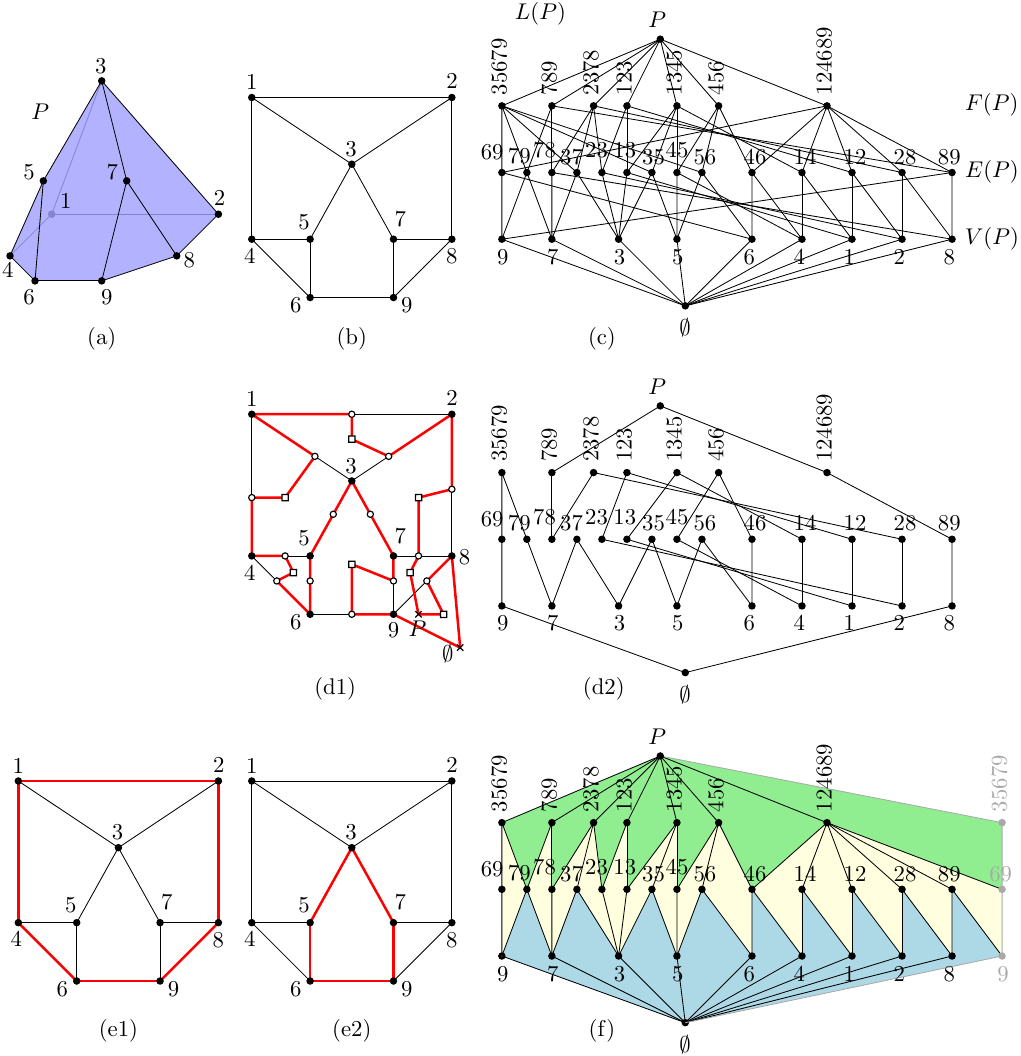}
}
\caption{(a) A 3-polytope~$P$; (b) a Schlegel diagram of~$P$, a 3-connected plane graph; (c) the face lattice of~$P$; (d1)+(d2) a Hamiltonian cycle in $G(L(P))$; (e1)+(e2) two distinct facet-Hamiltonian cycles of~$P$ (of different lengths); (f) a rhombic strip of~$G(L(P))$.
The rhombi are colored according to ranks for clarity.
In this and all of the following pictures of rhombic strips in our paper, we display a grayed-out copy of the leftmost vertices at the right-hand side, in order to depict the `wrap-around' edges on the sphere.}
\label{fig:example}
\end{figure}

\subsection{Hamiltonicity of the skeleton}

Another fundamental property with algorithmic applications is to search for a longest path instead of a shortest path on the skeleton.
Formally, a \defi{Hamiltonian path/cycle} is a path/cycle on the skeleton that visits every vertex exactly once.
The hypercube, permutahedron and associahedron are known to have Hamiltonian cycles by classical algorithms.
Specifically, the well-known \defi{binary reflected Gray code}, first described in patents by the Bell Labs researchers George R.~Stibitz and Frank Gray from 1943 and 1953, respectively, computes a Hamiltonian cycle in the hypercube.
The \defi{Steinhaus-Johnson-Trotter} algorithm \cite{MR0157881,MR0159764,DBLP:journals/cacm/Trotter62} computes a Hamiltonian cycle in the permutahedron.
Similar constructions of Hamiltonian cycles are also known for the associahedron~\cite{MR1239499,MR1723053}.
More generally, all graph associahedra admit a Hamiltonian cycle if the underlying graph~$H$ has at least two edges~\cite{MR3383157}.
If $H$ is chordal, then a Hamiltonian path can be computed efficiently by a simple greedy algorithm~\cite{DBLP:journals/talg/CardinalMM25}.
Furthermore, a variant of that algorithm can be used to compute a Hamiltonian path in all quotientopes~\cite{MR4344032}.

The aforementioned Hamiltonian paths and cycles on combinatorial polytopes are instances of so-called \defi{combinatorial Gray codes}~\cite{MR1491049,MR4649606}.
This term refers to a listing of combinatorial objects such that any two consecutive objects differ in a `small change'.
In the case of hypercubes, permutahedra, and associahedra, respectively, this is a single bit being flipped, an adjacent transposition or tree rotation being applied.
The purpose of such `small change' listings is to enable fast generation algorithms for the combinatorial objects, ideally in time~$\cO(1)$ per generated object.
Such algorithms are often referred to as \defi{loopless}.

\subsection{Hamiltonicity of the face lattice}

In this work, we aim to list not only vertices and edges of a polytope~$P$, but also the faces of all other dimensions of~$P$.
For this we consider the \defi{face lattice}~$L(P)$, i.e., the inclusion order of all faces of~$P$, which captures the complete combinatorial structure of~$P$; see Figure~\ref{fig:example}~(a)--(c).
The cover relations in the face lattice are pairs of faces whose dimension differs by~1 where one is a subface of the other.
The face lattice includes two special faces, called \defi{trivial} faces, namely the empty set~$\emptyset$ that constitutes the unique minimum of the face lattice~$L(P)$, and the entire polytope~$P$ that constitutes the unique maximum of~$L(P)$.
The dimension of the trivial face~$\emptyset$ is defined to be~$-1$.

We seek a Hamiltonian cycle in the cover graph of the face lattice, denoted $G(L(P))$, i.e., a cyclic listing of all faces of~$P$, in which every face appears exactly once, such that the dimension of any two cyclically consecutive faces differs by 1, and one is a subface of the other; see Figure~\ref{fig:example}~(d1)+(d2).
In other words, this is a Gray code listing of all the faces of~$P$ with respect to the inclusion order.
While obviously there is a large variety of computational problems related to face lattices of polytopes (see, e.g., \cite{MR1289077,MR1482230,MR1927137}), the question for Hamiltonicity of its cover graph appears to be novel to the best of our knowledge.

\subsection{A brave conjecture and some evidence}
\label{sec:brave}

We raise the following brave conjecture.

\begin{conjecture}
\label{conj:LP-HC}
For any polytope~$P$ of dimension~$d\geq 1$, the graph~$G(L(P))$ has a Hamiltonian cycle.
\end{conjecture}

\begin{figure}
\begin{tabular}{ccccccc}
(a1) & (a2) & (a3) & (a4) & (b2) & (b3) & (b4) \\
\raisebox{-\height}{\includegraphics{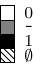}} &
\raisebox{-\height}{\includegraphics{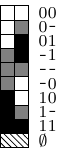}} &
\raisebox{-\height}{\includegraphics{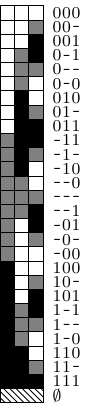}} &
\raisebox{-\height}{\includegraphics{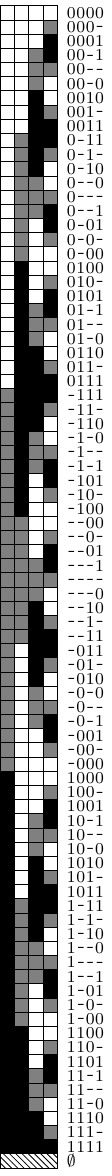}} &
\raisebox{-\height}{\includegraphics{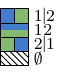}} &
\raisebox{-\height}{\includegraphics{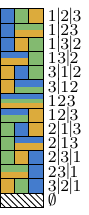}} &
\raisebox{-\height}{\includegraphics{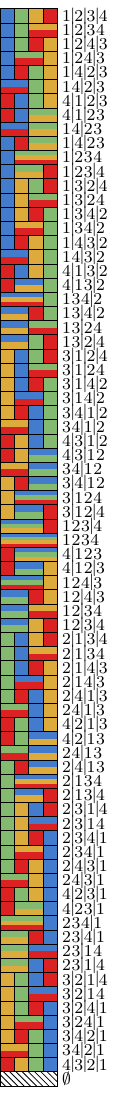}} \\
\end{tabular}
\caption{Hamiltonian cycles in the cover graph of the face lattice of hypercubes~$Q_n$ for $n=1,2,3,4$ (a1)--(a4) and permutahedra~$\Pi_n$ for $n=2,3,4$ (b2)--(b4).
The encoding for hypercubes is 0=white, 1=black, \hyph=gray and for permutahedra it is 1=blue, 2=green, 3=yellow, 4=red, where for values in the same block of an ordered set partition the corresponding rectangle is striped horizontally with the colors of values in that block.
}
\label{fig:hc-cube-perm}
\end{figure}

\begin{figure}
\includegraphics{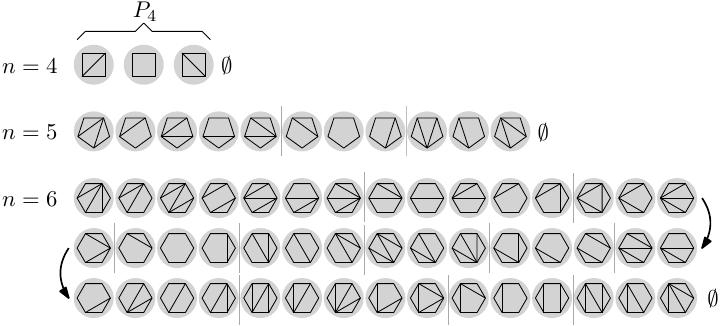}
\caption{Hamiltonian cycles in the cover graph of the face lattice of the associahedron~$A_n$ for $n=4,5,6$.
The vertical bars delimit groups of dissections obtained from the same parent dissection in the previous sequence.}
\label{fig:hc-asso}
\end{figure}

For polytopes of dimension~1 and~2, Conjecture~\ref{conj:LP-HC} is trivially true.
For a simplex~$P$, the face lattice~$L(P)$ is the hypercube, and therefore a Hamiltonian cycle in~$G(L(P))$ is given by the aforementioned binary reflected Gray code.
In this paper, we verify Conjecture~\ref{conj:LP-HC} for a number of further interesting polytopes:
\begin{itemize}[leftmargin=4mm]
\item hypercubes (\ourresult{Theorem~\ref{thm:LQ-HC}}; see Figure~\ref{fig:hc-cube-perm}~(a1)--(a4));
\item permutahedra and $B$-permutahedra (\ourresult{Theorems~\ref{thm:LPi-HC} and~\ref{thm:LPib-HC}}; see Figure~\ref{fig:hc-cube-perm}~(b2)--(b4));
\item associahedra (\ourresult{Theorem~\ref{thm:LAsso-HC}}; see Figure~\ref{fig:hc-asso});
\item cyclic polytopes (\ourresult{Theorem~\ref{thm:cyclic}});
\item 3-dimensional polytopes (\ourresult{Theorem~\ref{thm:3d-HC}});
\item graph associahedra of chordal graphs (\ourresult{Theorem~\ref{thm:LGAsso-HC}});
\item quotientopes (\ourresult{Theorem~\ref{thm:LSn-HC}}).
\end{itemize}
In this list and in the rest of the introduction, the new theorems established in this work are highlighted in magenta.
For a polytope~$P$, we write $P^*$ for the \defi{polar} polytope of~$P$.
The face lattice~$L(P^*)$ is obtained by turning the face lattice~$L(P)$ upside down, an operation that clearly preserves Hamiltonian cycles.
Consequently, the polars of all the polytopes mentioned before also satisfy Conjecture~\ref{conj:LP-HC}, including for example all cross-polytopes, which are the polars of hypercubes.
As mentioned before, while hypercubes, permutahedra and associahedra are special cases of graph associahedra and quotientopes, it is still worth to discuss them separately, as this gives new explicit listings and algorithms.

\begin{figure}[h!]
\includegraphics{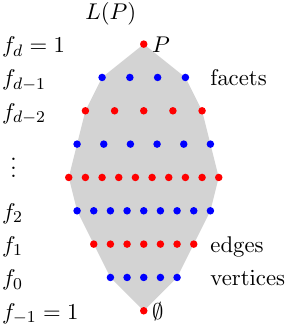}
\caption{The two partition classes of the bipartite graph~$G(L(P))$.}
\label{fig:fvector}
\end{figure}

Strikingly, we did not find any counterexamples to Conjecture~\ref{conj:LP-HC}.
Additional hope for the conjecture is fueled by the following observation:
The cover graph of the face lattice~$G(L(P))$ is bipartite, with the partition classes given by the parity of the dimension of the faces; see Figure~\ref{fig:fvector}.

We write $d\ass\dim(P)$ for the dimension of~$P$, and $f_i$, $i=-1,0,1,\ldots,d$ for the number of $i$-dimensional faces of~$P$, i.e., $f=(f_{-1},f_0,f_1,\ldots,f_{d-1},f_d)=(1,f_0,f_1,\ldots,f_{d-1},1)$ is the \defi{$f$-vector} of~$P$.
A necessary condition for~$G(L(P))$ to have a Hamiltonian cycle is that its partition classes have the same size, i.e.,
\[
\sum_{\substack{i=-1,\ldots,d \\ i\text{ odd}}} f_i = \sum_{\substack{i=-1,\ldots,d \\ i\text{ even}}} f_i
\]
which is equivalent to Euler-Poincar\'e's famous formula
\begin{equation}
\label{eq:euler}
\sum_{i=-1}^d (-1)^i f_i=0,
\end{equation}
valid \emph{for any polytope~$P$}.
In fact, our construction of a Hamiltonian cycle in~$G(L(P))$ for 3-dimensional polytopes yields another proof of~\eqref{eq:euler} for the case~$d=3$, i.e., for the formula
\[-1+f_0-f_1+f_2-1=0,\]
which can be rewritten in the more familiar form
\[|V(P)|-|E(P)|+|F(P)|=2,\]
where $V(P)$, $E(P)$, and~$F(P)$ are the sets of vertices, edges, and faces of~$P$, respectively.
Specifically, a bijection between the two partition classes of~$G(L(P))$ is obtained by taking every second edge from our Hamiltonian cycle.

Similarly, proving Conjecture~\ref{conj:LP-HC} for larger values of~$d$ (maybe $d=4$ as a next step), would provide another proof of~\eqref{eq:euler} in the more general $d$-dimensional case.

\subsection{Facet-Hamiltonian cycles}

Akitaya, Cardinal, Felsner, Kleist and Lauff~\cite{MR4863586} investigated a new type of cycles in polytopes.
Given a polytope~$P$ of dimension~$d$, each $(d-1)$-dimensional face is called a \defi{facet}.
A \defi{facet-Hamiltonian cycle} in~$P$ is a cycle~$C$ in the skeleton of~$P$ that enters and leaves every facet of~$P$ exactly once; see Figure~\ref{fig:example}~(e1)+(e2).
Formally, for every facet~$F$ of~$P$, the intersection~$C\cap F$ is connected and nonempty.
The notion of a \defi{facet-Hamiltonian path} is defined analogously.
The authors show that permutahedra and $A$-, $B$-, $C$-, $D$-associahedra admit facet-Hamiltonian cycles, and the same is true for graph associahedra of wheels, fans and complete split graphs.
For graph associahedra of complete bipartite graphs and caterpillars, they construct facet-Hamiltonian paths.

\begin{figure}[htb]
\includegraphics{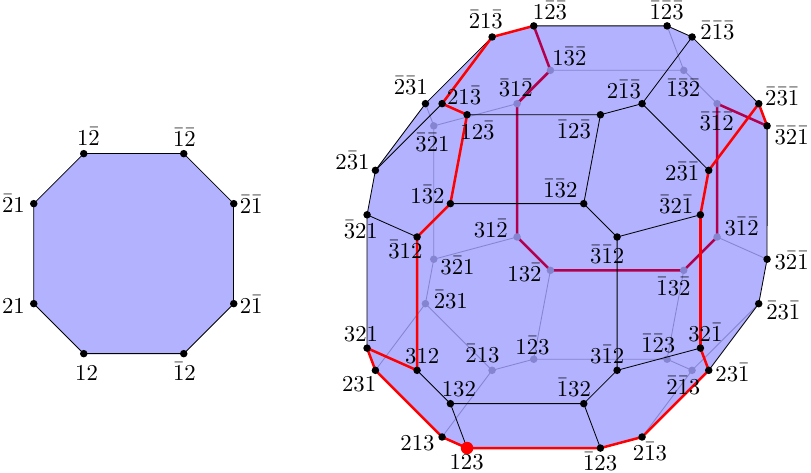}
\caption{The 2- and 3-dimensional $B$-permutahedra.
A facet-Hamiltonian cycle in the 3-dimensional $B$-permutahedron is highlighted.}
\label{fig:bperm}
\end{figure}

One of the conjectures from~\cite{MR4863586} concerns the \defi{$B$-permutahedron}, which has as vertices all signed permutations of~$[n]$, i.e., permutations of~$[n]$ in which every entry has a positive/negative sign; see Figure~\ref{fig:bperm}.
In the figure, entries of the permutation with a negative sign are overlined.
The edges of the $B$-permutahedron connect pairs of signed permutations that either differ in an adjacent transposition, preserving all signs, or in a complementation of the sign of the first entry.

\begin{conjecture}[{\cite[Conj.~1]{MR4863586}}]
\label{conj:bperm}
The $B$-permutahedron of any dimension $d\geq 2$ has a facet-Hamiltonian cycle.
\end{conjecture}

In this work, we resolve Conjecture~\ref{conj:bperm} affirmatively and establish it as a theorem (\ourresult{Theorem~\ref{thm:BPi-fHC}}).

\subsection{Rhombic strips}
\label{sec:strips}

One of the key concepts introduced in~\cite{MR4863586} is that of a \defi{rhombic strip}, which is a plane spanning subgraph of the cover graph of the face lattice of some polytope such that every face in the embedding is a \defi{diamond}, i.e., a 4-cycle that spans three consecutive levels of the face lattice; see Figure~\ref{fig:example}~(f).
Furthermore, the rhombic strip wraps around at the left and right boundary like on a sphere with $\emptyset$ and the entire polytope~$P$ at the south and north pole, respectively.
Note that if $G(L(P))$ admits a rhombic strip, then for any fixed level/dimension~$k\in\{1,\ldots,d-1\}$, $d\ass\dim(P)$, we may consider all diamonds whose two middle vertices lie on level~$k$.
In Figure~\ref{fig:example}~(f), these cyclic orderings of diamonds for $k=1,2,3$ are drawn blue, yellow and green, respectively.
The plane embedding of the rhombic strip forces a cyclic ordering of these diamonds, which forces a cyclic ordering of all $k$-dimensional faces of~$P$, with the property that any two consecutive faces have a common superface of dimension~$k+1$ (this is the top vertex of the corresponding diamond) and a common subface of dimension~$k-1$ (this is the bottom vertex of the corresponding diamond), i.e., we obtain a Gray code ordering of the faces of dimension~$k$.
As the rhombic strip contains the faces of all dimensions, we obtain such orderings simultaneously for all~$k=1,\ldots,d-1$.
In particular, for $k=1$ we obtain a cyclic ordering of the vertices of~$P$, where any two consecutive vertices are joined by an edge, which is precisely a Hamiltonian cycle in the skeleton~$G(P)$, i.e., a rhombic strip in~$G(L(P))$ implies a Hamiltonian cycle in~$G(P)$.
Furthermore, by the earlier remarks about the polar polytope~$P^*$, a rhombic strip in~$G(L(P))$ implies a Hamiltonian cycle in~$G(P^*)$.
Combining these observations, we see that a necessary condition on~$P$ for $G(L(P))$ to admit a rhombic strip is that~$G(P)$ and~$G(P^*)$ both admit a Hamiltonian cycle.
Note that for 3-dimensional polytopes~$P$, the skeleton~$G(P^*)$ is simply the dual graph of~$G(P)$.

Akitaya, Cardinal, Felsner, Kleist, and Lauff~\cite{MR4863586} proved that in order to resolve Conjecture~\ref{conj:bperm}, it is sufficient to find a rhombic strip in the cover graph of the face lattice of the hypercube.
Our proof of Conjecture~\ref{conj:bperm} follows this route, by constructing a rhombic strip in the face lattice of the hypercube of any dimension (\ourresult{Theorem~\ref{thm:LQ-RS}}).
Furthermore, we characterize 3-dimensional polytopes~$P$ for which $G(L(P))$ admits a rhombic strip (\ourresult{Theorem~\ref{thm:3d-RS}}), and this characterization is a strengthening of the aforementioned necessary condition that~$G(P)$ and its dual graph must admit a Hamiltonian cycle.
Using this characterization, we provide infinitely many examples of 3-dimensional polytopes~$P$ for which the necessary condition holds but~$G(L(P))$ does not admit a rhombic strip (\ourresult{Theorem~\ref{thm:3d-RS-counter}} and \ourresult{Proposition~\ref{prop:fano}}, illustrated in Figures~\ref{fig:trunc-tetra1} and~\ref{fig:fano}, respectively).

\subsection{Efficient Gray code algorithms}
\label{sec:algo}

Several of our aforementioned results yield new and efficient Gray code algorithms for listing various combinatorial objects.
Specifically, the constructions for establishing Hamiltonian cycles in the cover graphs of face lattices of hypercubes, permutahedra, $B$-permutahedra, associahedra and cyclic polytopes (Theorems~\ref{thm:LQ-HC}, \ref{thm:LPi-HC}, \ref{thm:LPib-HC}, \ref{thm:LAsso-HC} and~\ref{thm:cyclic} respectively) can be turned straightforwardly into loopless algorithms to compute those cycles, i.e., they work in optimal time~$\cO(1)$ per visited vertex.
Furthermore, the algorithms only store the currently visited vertex in memory, plus a few additional data structures.
In particular, from the permutahedron we obtain a loopless Gray code algorithm for listing all ordered set partitions of~$[n]$, where consecutive partitions differ in either merging or splitting two adjacent sets/blocks in the ordered partition; see Figure~\ref{fig:hc-cube-perm}~(b2)--(b4).
Furthermore, from the associahedron we obtain a loopless Gray code algorithm for listing all dissections of a convex $n$-gon, where consecutive dissections differ in either adding or removing a single diagonal; see Figure~\ref{fig:hc-asso}.
Implementations of these Gray codes in C++ are available for download and experimentation on the Combinatorial Object Server website~\cite{cos_part}, \cite{cos_kary} and~\cite{cos_cyclic}, respectively.

The construction of a facet-Hamiltonian cycle in the $B$-permutahedron described in the proof of Conjecture~\ref{conj:bperm} (Theorem~\ref{thm:BPi-fHC}) also translates to a polynomial-time and -space algorithm for computing the cycle.
An implementation in C++ can be found at~\cite{cos_sperm}.

\subsection{Key ideas and methods}
\label{sec:ideas}

The constructions of Hamiltonian cycles in the face lattices of most of the polytopes listed after Conjecture~\ref{conj:LP-HC} (Theorems~\ref{thm:LQ-HC}, \ref{thm:LPi-HC}, \ref{thm:LPib-HC}, \ref{thm:LAsso-HC}, \ref{thm:LGAsso-HC} and~\ref{thm:LSn-HC}) can be seen as applications of the zigzag framework for combinatorial generation pioneered in~\cite{MR4391718} and extensively used in~\cite{MR4344032,MR4598046,DBLP:journals/talg/CardinalMM25,MR4614413,MR4775168}, thus further extending the reach of this framework.
The idea is to construct a Hamiltonian cycle by induction on the dimension, where in the induction step~$d\rightarrow d+1$, each face in the $d$-dimensional polytope gives rise to a sequence of `children', namely a sequence of faces in the $(d+1)$-dimensional polytope which form a path in the cover graph of the face lattice.
Consequently, those children faces can be visited consecutively, alternatingly in forward (=zig) or backward (=zag) order, and the transitions between the children sequences are inherited inductively from the previous Hamiltonian cycle.

On the other hand, our result for cyclic polytopes (Theorem~\ref{thm:cyclic}) is obtained by restricting the binary reflected Gray code to the subset of bitstrings that describe faces of the polytope.
Furthermore, our results for 3-dimensional polytopes (Theorems~\ref{thm:3d-HC}, \ref{thm:3d-RS}, \ref{thm:3d-RS-counter} and Proposition~\ref{prop:fano}), are proved by direct combinatorial arguments on 3-connected planar graphs, using Steinitz' theorem.

As outlined before, Conjecture~\ref{conj:bperm} is proved following the approach suggested in~\cite{MR4863586}, by constructing a rhombic strip in the face lattice of the hypercube of any dimension.
This is based on the observation that the \mbox{$B$-permutahedron} is obtained by truncating every (non-trivial) face of the hypercube; see~Figure~\ref{fig:truncate}.
Therefore, a vertex of the $B$-permutahedron corresponds to a maximal chain in the face lattice of the hypercube, and a facet of the $B$-permutahedron corresponds to the union of all chains through a certain face of the hypercube, i.e., the union of the upset and downset of this face, forming an `hourglass' centered at this face.
Consequently, sweeping a maximal chain from left to right through the rhombic strip enters and leaves each `hourglass', i.e., each facet of the $B$-permutahedron, exactly once.

\subsection{Outline of this paper}

In Section~\ref{sec:prelim}, we provide some terminology and notation that will be used throughout this paper.
We prove our results for the different types of polytopes starting with the elementary ones, namely simplices, hypercubes, permutahedra, associahedra, and cyclic polytopes in Sections~\ref{sec:simplex}, \ref{sec:cube}, \ref{sec:perm}, \ref{sec:asso}, \ref{sec:cyclic} respectively, followed by the 3-dimensional polytopes in Section~\ref{sec:3dim}, before proceeding to more advanced ones, namely graph associahedra and quotientopes in Sections~\ref{sec:Gasso} and~\ref{sec:quotient}, respectively.
An interlude is Section~\ref{sec:truncate}, where we present the proof of Conjecture~\ref{conj:bperm}, i.e., we establish the existence of facet-Hamiltonian cycles in $B$-permutahedra.

In Section~\ref{sec:open}, we conclude with some open questions and ideas on how to tackle Conjecture~\ref{conj:LP-HC} in general.

\section{Preliminaries}
\label{sec:prelim}

Let $(P,<)$ be a poset.
We say that two distinct elements $x,y\in P$ with $x<y$ are in a \defi{cover relation}, denoted $x\lessdot y$, if there is no~$z\in P$ with $x<z<y$.
The \defi{cover graph} of~$P$, denoted $G(P)$, has as vertices all elements of~$P$, and an edge~$(x,y)$ for every cover relation~$x\lessdot y$.
The \defi{downset} of some~$x\in P$ is the set of all~$y\in P$ for which $y\leq x$.
Similarly, the \defi{upset} of some~$x\in P$ is the set of all~$y\in P$ for which $x\leq y$.
An \defi{interval}~$[x,y]$ is the set of all~$z\in P$ with $x\leq z\leq y$, i.e., it is the intersection of the upset of~$x$ and the downset of~$y$.
A \defi{chain} in~$P$ is a sequence of elements $(x_1,\ldots,x_k)$ from~$P$ such that $x_1\lessdot x_2\lessdot\cdots\lessdot x_k$.
It corresponds to a monotonically increasing path in the cover graph.
For any $x,y\in P$, the \defi{join} $x\vee y$ of~$x$ and~$y$ denotes the unique smallest element $z\in P$, such that $x\leq z$ and $y\leq z$, and the \defi{meet} $x\wedge y$ of~$x$ and~$y$ denotes the unique largest element $z\in P$ such that $z\le x$ and $z\le y$ (if such elements exist).
If for all $x,y\in P$, both the join~$x\vee y$ and meet~$x\wedge y$ exist and are unique, then $(P,<)$ is called a \defi{lattice}.

A~poset $(P,<)$ is \defi{graded} if there is a function~$\rho:P\rightarrow \mathbb{Z}$ such that $\rho(y)=\rho(x)+1$ if $(x,y)$ is a cover relation in~$P$.
The function~$\rho$ is called \defi{rank function}, and $\rho(x)$ is called the \defi{rank} of~$x\in P$.
A poset is \defi{$M_3$-free} if there are no five distinct elements~$a,b_1,b_2,b_3,c$ with $a\lessdot b_1\lessdot c$, $a\lessdot b_2\lessdot c$ and $a\lessdot b_3\lessdot c$.

A \defi{polytope~$P$} is the convex hull of a finite set of points in~$\mathbb{R}^d$.
Equivalently, it is a bounded intersection of half-spaces.
Its \defi{dimension} is the dimension of the smallest affine space containing it.
If $P$ is $d$-dimensional, we sometimes refer to it as a \defi{$d$-polytope}.
We write~$G(P)$ for the graph of the (1-)skeleton of~$P$.
A \defi{face} of~$P$ is the intersection of~$P$ with a hyperplane such that all of~$P$ lies on the same side or on the hyperplane.
Each face is itself a polytope and thus has a dimension.
The combinatorial structure of~$P$ is captured by its~\defi{face lattice}, denoted~$L(P)$, which is the inclusion order of all faces of~$P$.
Figure~\ref{fig:example}~(a) shows a 3-dimensional polytope, and part~(c) of the figure shows its face lattice.
The face lattice has the empty set~$\emptyset$ as its unique minimum, and the full polytope~$P$ as its unique maximum.
We refer to these two special faces of~$P$ as \defi{trivial} faces.
The face lattice is graded, where we can take as the rank function the dimensions of the faces.
The unique minimum~$\emptyset$ of~$L(P)$ is assigned the rank (dimension)~$-1$.
Furthermore, the 0-dimensional faces are called \defi{vertices} and the 1-dimensional faces are called \defi{edges}.
If $P$ is $d$-dimensional, then the faces of dimension~$d-1$ are called \defi{facets}.
We write~$F_k(P)$ for the set of $k$-dimensional faces of~$P$.
Given two distinct faces $F,G\in L(P)$, we write $F\subset G$ if $F$ is a subface of~$G$, i.e., a subset of~$G$ that is a face of~$P$.
Furthermore, we write $F\cov G$ if $F$ and~$G$ form a cover relation in~$L(P)$, which means that $F\subset G$ and the dimension/rank of~$G$ is exactly one higher than that of~$F$.
A \defi{flag} of~$P$ is a chain in~$L(P)$ that starts at~$\emptyset$ and ends at~$P$, i.e., this chain contains exactly one face from each rank.

For a 3-dimensional polytope~$P$, we write $V(P)\ass F_0(P)$, $E(P)\ass F_1(P)$, $F(P)\ass F_2(P)$ for the sets of vertices, edges, and 2-dimensional faces of~$P$, respectively.
Similarly, for a plane graph~$H$, we write $V(H)$, $E(H)$, and~$F(H)$ for the set of vertices, edges, and faces of~$H$, respectively, and we refer to the elements of these sets as \defi{cells} (so in the setting of plane graphs the word face is strictly reserved for the 2-dimensional cells).
The projection of a 3-dimensional polytope onto the plane through a point just outside one of its faces yields a plane graph, called its \defi{Schlegel diagram}.
By Steinitz' theorem, the Schlegel diagrams of 3-dimensional polytopes are precisely the 3-connected plane graphs.
For a plane graph~$H$, we write~$L(H)$ for the inclusion order of its cells~$\{\emptyset\}\cup V(H)\cup E(H)\cup F(H)\cup \{H\}$, which includes the two \defi{trivial} cells~$\emptyset$ and~$H$.
While $L(H)$ is a graded poset, it is in general not a lattice.
For example, if~$H$ is a cycle, then any two edges~$e,e'$ on the cycle have the inner face~$f$ and the outer face~$\ol{f}$ covering both of them, i.e., $e,e'\cov f$ and~$e,e'\cov\ol{f}$.
For a polytope~$P$, the face lattice~$L(P)$ is known to be $M_3$-free~\cite[Thm.~2.7~(iii)]{MR1311028}.
Similarly, if~$H$ is 2-connected, then~$L(H)$ is $M_3$-free.
The unbounded face of a plane graph~$H$ is referred to as the \defi{outer face}, and the remaining faces are called \defi{inner} faces.

A \defi{Hamiltonian cycle} in a graph is a cycle that visits every vertex exactly once.
A Hamiltonian cycle of~$G(L(P))$ corresponds to a cyclic listing of all faces of~$P$ such that any two consecutive faces form a cover relation in~$L(P)$; see Figure~\ref{fig:example}~(d1)+(d2).

A \defi{facet-Hamiltonian cycle} of a polytope~$P$ is a cycle~$C$ in the skeleton graph~$G(P)$ with the property that every facet of~$P$ has a nonempty and connected intersection with~$C$; see Figure~\ref{fig:example}~(e1)+(e2).
In terms of the face lattice~$L(P)$, such a cycle enters and leaves the downset of every facet exactly once.

Given a graded poset~$P$ with a unique minimum and maximum, a \defi{rhombic strip} is a spanning subgraph of the cover graph~$G(P)$ embedded on the sphere subject to the following conditions; see Figure~\ref{fig:example}~(f):
\begin{itemize}[leftmargin=4mm]
\item the height of vertices in the embedding is given by their rank, so that the unique minimum and maximum are embedded at the south and north pole, respectively;
\item edges are drawn as straight lines and there are no crossings between them;
\item every face is a \defi{diamond}, i.e., a 4-cycle $(a,b,c,d)$, where $a$ and~$c$ have the same rank~$k$, and $b$ and~$d$ have ranks $k-1$ and~$k+1$, respectively.
\end{itemize}

The following lemma formalizes the observations made in Section~\ref{sec:strips}.

\begin{lemma}
\label{lem:strip-order}
Let $d\geq 2$ and let $P$ be a graded poset with ranks~$-1,\ldots,d$ that has a unique minimum and maximum and that admits a rhombic strip~$R$.
Then we have the following:
\begin{enumerate}[label=(\roman*),leftmargin=8mm,topsep=1mm]
\item
For every rank~$k=0,\ldots,d-1$ and any two elements~$x,x'$ that are consecutive in the cyclic ordering of rank~$k$ elements given by~$R$, there is an element~$y$ (of rank~$k+1$) such that both $x\lessdot y$ and~$x'\lessdot y$ are edges in~$R$ and an element~$z$ (of rank~$k-1$) such that both $z\lessdot x$ and~$z\lessdot x'$ are edges in~$R$.
\item
For any element~$x$ of rank~$k\in\{-1,\ldots,d-1\}$ of~$P$, there is at least one edge $x\lessdot y$ in~$R$, and for any element~$x$ of rank~$k\in\{0,\ldots,d\}$ of~$P$, there is at least one edge $z\lessdot x$ in~$R$.
\item
If $P$ is $M_3$-free, then for any element~$x$ of rank~$k\in\{0,\ldots,d-1\}$ of~$P$, there are at least two edges $x\lessdot y$ and $x\lessdot y'$ or a least two edges $z\lessdot x$ and $z'\lessdot x$ in~$R$.
\end{enumerate}
As a consequence of~(i), we have the following:
\begin{enumerate}[start=4,label=(\roman*),leftmargin=8mm,topsep=1mm]
\item If $P$ is the face lattice of a polytope~$Q$, i.e., $P=L(Q)$, then $G(Q)$ and~$G(Q^*)$ both have a Hamiltonian cycle.
\item If $P$ is the inclusion order of the cells of a 2-connected plane graph~$H$, i.e., $P=L(H)$, then $H$ and its dual graph both have a Hamiltonian cycle.
\end{enumerate}
\end{lemma}

\begin{proof}
Part~(i) follows by considering the plane embedding of all diamonds that have two vertices on rank~$k$, as described in Section~\ref{sec:strips}.
Part~(ii) is a direct consequence of~(i).
Part~(iii) follows directly from~(ii), using that if there was exactly one edge~$x\lessdot y$ and exactly one edge~$z\lessdot x$, then there would be two diamonds~$(x',z,x,y)$ and~$(x,z,x'',y)$ that create a copy of~$M_3$.

Parts~(iv) and~(v) follow by applying~(i) to the ranks $k=1$ and~$k=d-1$.
\end{proof}

Given strings~$x$ and~$y$, we write $xy$ for their concatenation.
This operation extends to sets of strings~$Y$ in the natural way, specifically $xY\ass\{xy\mid y\in Y\}$.
Similarly, for a sequence of strings $Y=(y_1,\ldots,y_n)$, we define $xY\ass(xy_1,\ldots,xy_n)$.
For a string~$x$ and an integer~$n$, we write $x^n$ for the $n$-fold concatenation of~$x$ with itself.
Given a sequence~$x=(x_1,\ldots,x_n)$, we write $\rev(x)=(x_n,x_{n-1},\ldots,x_1)$ for the reversed sequence.

Given a positive integer $n$, we define $[n]\ass\{1,\ldots,n\}$, and given two positive integers $a\le b$, we define $[a,b]\ass\{n\in\mathbb{N}\mid a\le n\le b\}$ and $\openint{a,b}\ass [a,b]\setminus\{a,b\}$.

\section{Simplices}
\label{sec:simplex}

\subsection{The polytope and its face lattice}

\begin{figure}[h!]
\includegraphics{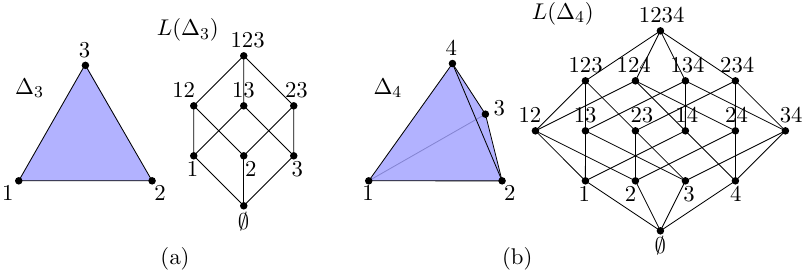}
\caption{The 3- and 4-simplex and their face lattices, the Boolean lattices of dimension~3 and~4, respectively.}
\label{fig:simplex}
\end{figure}

The \defi{Boolean lattice~$Q_n$} is the inclusion order of all $2^n$ subsets of~$[n]$.
We identify subsets of~$[n]$ by their characteristic vectors of length~$n$, i.e., we think of $G(Q_n)$ as the graph that has~$\{0,1\}^n$ as its vertex set, and an edge between any two bitstrings that differ in a single bit.
The \defi{\mbox{$n$-simplex~$\Delta_n$}} is the convex hull of the $n$ unit vectors in~$\mathbb{R}^n$; see Figure~\ref{fig:simplex}.
It has $n$ vertices, and any $k$-subset of vertices forms a $(k-1)$-face, for $k=1\ldots,n$, which is itself a simplex.
Specifically, the 1-simplex is a point, the 2-simplex is a line segment, the 3-simplex is a triangle, the 4-simplex is a tetrahedron, etc.
Furthermore, the face lattice of the $n$-simplex~$\Delta_n$ is the Boolean lattice~$Q_n$.

\subsection{Hamiltonian cycles and rhombic strips in the face lattice}

The face lattice of the simplex is known to admit a Hamiltonian cycle and a rhombic strip.

\begin{theorem}[Folklore]
\label{thm:Q-HC}
For any $n\geq 2$, the graph~$G(L(\Delta_n))=G(Q_n)$ has a Hamiltonian cycle.
\end{theorem}

\begin{theorem}[Folklore]
\label{thm:Q-RS}
For any $n\geq 2$, the graph~$G(L(\Delta_n))=G(Q_n)$ has a rhombic strip.
\end{theorem}

Even though the notion of rhombic strips in general posets was introduced only recently in~\cite{MR4863586}, the concept was implicitly used long before that in special cases.
In particular, we refer to Theorem~\ref{thm:Q-RS} as a folklore result because it follows directly from classical and well-known constructions of Venn diagrams, a connection that we will briefly explain at the end of this section.

We revisit two proofs for Theorem~\ref{thm:Q-RS}, and later combine the techniques presented here for constructing rhombic strips for the face lattice of the hypercube.
The following proofs are based on two classical constructions of the binary reflected Gray code, and in fact yield the same rhombic strip, under an isomorphism that renames elements of the ground set.

\begin{proof}[Proof~1 of Theorem~\ref{thm:Q-RS}]
A rhombic strip is a spanning subgraph of~$G(Q_n)$ with certain additional properties, and hence it has $\{0,1\}^n$ as its vertex set.
We specify such a subgraph by the set of edges of~$G(Q_n)$ that it contains.
Specifically, we construct a rhombic strip~$R_n\seq G(Q_n)$ by induction on~$n$.
The notations used in this construction are illustrated in Figures~\ref{fig:Q-RS1a} and~\ref{fig:Q-RS1b}.
The construction also maintains two chains~$C_n=(x_0,\ldots,x_n)\seq R_n$ (red in the figures) and $D_n=(y_0,\ldots,y_n)\seq R_n$ (blue in the figures) with $x_0=y_0=0^n$ and~$x_n=y_n=1^n$, where $C_n$ and $D_n$ contains the leftmost and rightmost vertices across all ranks, respectively, and both chains are connected by edges $Z_n\ass\{(x_i,y_{i+1})\mid i=1,\ldots,n-2\}\seq R_n$ (in the gray-shaded area in the figures).
We refer to the edges in~$Z_n$ as \defi{zipper edges}, as they connect the left and right side of the plane drawing to make a rhombic strip on the sphere, like a zipper that connects two pieces of fabric.

\begin{figure}[b!]
\centering
\includegraphics[page=1]{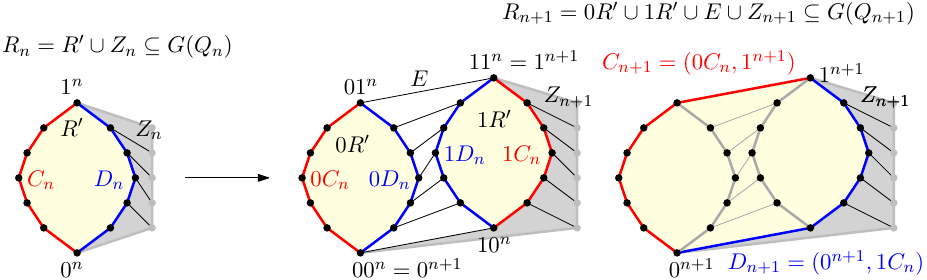}
\caption{Schematic illustration of proof~1 of Theorem~\ref{thm:Q-RS}.}
\label{fig:Q-RS1a}
\end{figure}

\begin{figure}[t!]
\centering
\includegraphics[page=2]{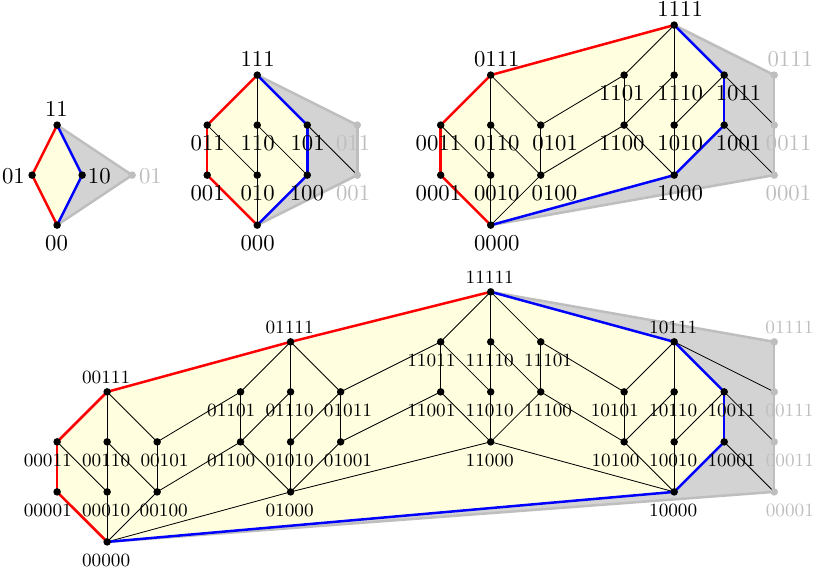}
\caption{The rhombic strips obtained from proof~1 of Theorem~\ref{thm:Q-RS} for $n=2,3,4,5$.}
\label{fig:Q-RS1b}
\end{figure}

For the base case $n=2$ of the construction we take $R_2\ass G(Q_2)$, $C_2\ass(00,01,11)$ and $D_2\ass(00,10,11)$, and in this case we have~$Z_2=\emptyset$.

For the induction step, let $R_n,C_n=(x_0,\ldots,x_n),D_n=(y_0,\ldots,y_n)$ be given for some~$n\geq 2$.
We construct the rhombic strip~$R_{n+1}\seq G(Q_{n+1})$ as follows:
We remove the zipper edges from~$R_n$, and consider the plane graph $R'\ass R_n\setminus Z_n$.
We take a copy of~$0R'$, i.e., the graph $R'$ in which all vertices are prefixed with a 0-bit, maintaining the same embedding.
In addition, we take a mirrored copy of $1R'$, i.e., the graph $R'$ in which all vertices are prefixed with a 1-bit and the embedding is mirrored along a vertical line.
The second (mirrored) graph is embedded to the right of the first one and one unit higher than the first one (note that prefixing with a 1-bit increases the ranks).
We connect the two graphs by the edges
\[ E\ass\{(0y_i,1y_i)\mid i=0,\ldots,n\} \]
and the zipper edges
\[ Z_{n+1}\ass\{(0x_i,1x_i)\mid i=1,\ldots,n-1\}, \]
which yields the rhombic strip
\[ R_{n+1}\ass0R'\cup 1R'\cup E\cup Z_{n+1}. \]
Furthermore, the new leftmost and rightmost chains are
\[ C_{n+1}\ass(0C_n,1^{n+1}) \text{  and  } D_{n+1}\ass(0^{n+1},1C_n), \]
respectively.

One can check that the vertices on the chains~$C_n$ and~$D_n$ have the explicit form
$x_i=0^{n-i}1^i$ for $i=0,\ldots,n$ and $y_i=10^{n-i}1^{i-1}$ for $i=1,\ldots,n$, respectively.
\end{proof}

\begin{proof}[Proof~2 of Theorem~\ref{thm:Q-RS}]
An \defi{$x$-monotone path} in a drawing of a graph is a path with the property that its vertices are embedded with strictly increasing abscissa values.
We construct a rhombic strip~$R_n\seq G(Q_n)$ by induction on~$n$.
The notations used in this construction are illustrated in Figures~\ref{fig:Q-RS2a} and~\ref{fig:Q-RS2b}.
The construction maintains a partition~$R_n=P_n\cup A_n\cup B_n\cup Z_n$ into four sets with the following properties:
\begin{itemize}[leftmargin=4mm]
\item $P_n$ is an $x$-monotone Hamiltonian path\footnote{It can be shown that this path follows the binary reflected Gray code, though with a reflected labeling of bits compared to the sequence defined in \cref{sec:cyclic}.} of~$G(Q_n)$ that starts at the vertex~$0^n$ and ends at a vertex of Hamming distance $1$ (red in the figures).
\item The sequence $C_n=(x_0,\ldots,x_n)$ of leftmost vertices across all ranks is a chain~$C_n\seq R_n$ (cyan in the figures), and the sequence $(0^n,D_n)=(0^n,y_1,\ldots,y_n)$ of rightmost vertices across all ranks is a chain~$(0^n,D_n)\seq R_n$ (pink in the figures), and the two chains are connected by the zipper edges $Z_n\ass\{(x_i,y_{i+1})\mid i=1,\ldots,n-2\}\seq R_n$ (in the gray-shaded area in the figures).
\item The edges $A_n$ and~$B_n$ are the edges of~$R_n$ strictly above and below~$P_n$, respectively (in the blue-shaded and yellow-shaded areas in the figures).
\end{itemize}

\begin{figure}[t!]
\centering
\includegraphics[page=3]{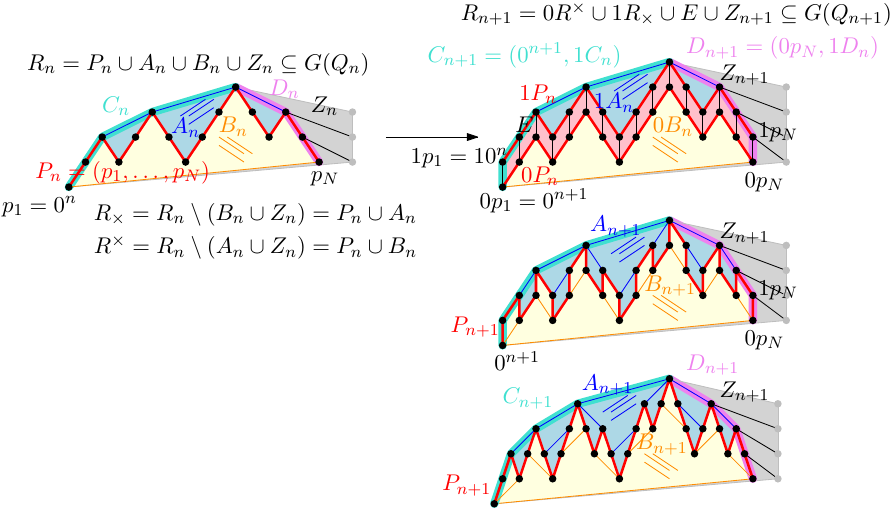}
\caption{Schematic illustration of proof~2 of Theorem~\ref{thm:Q-RS}.}
\label{fig:Q-RS2a}
\end{figure}

\begin{figure}[t!]
\centering
\includegraphics[page=4]{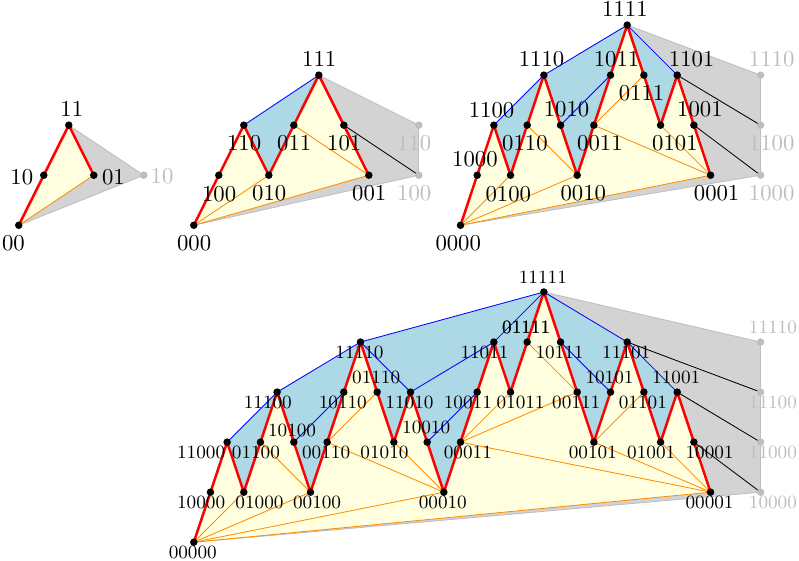}
\caption{The rhombic strips obtained from proof~2 of Theorem~\ref{thm:Q-RS} for $n=2,3,4,5$.}
\label{fig:Q-RS2b}
\end{figure}

For the base case $n=2$ of the construction we take $R_2\ass G(Q_2)$ with $P_2\ass (00,10,11,01)$, which gives $C_2=(00,10,11)$, $D_2=(01,11)$, $Z_2=\emptyset$, $A_2=\emptyset$ and~$B_2=\{(00,01)\}$.

For the induction step, let $P_n,C_n,D_n,A_n,B_n$ be given for some~$n\geq 2$.
We construct the rhombic strip~$R_{n+1}$ as follows:
We consider the plane graphs~$R^\times\ass R_n\setminus (A_n\cup Z_n)=P_n\cup B_n$ and $R_\times\ass R_n\setminus (B_n\cup Z_n)=P_n\cup A_n$.
We take a copy of~$0R^\times$ and a copy of~$1R_\times$, maintaining the same embedding.
Both graphs are embedded at the same positions horizontally, but the second graph is placed one unit higher than the first one, so that the corresponding vertices on the two copies of~$P_n$ line up vertically (note again that prefixing with~1 increases the ranks).
Let $P_n\assr(p_1,\ldots,p_N)$, $N=2^n$, be the sequence of vertices on the Hamiltonian path~$P_n$.
Then we connect the two graphs by the edges
\[ E\ass\{(0p_i,1p_i)\mid i=1,\ldots,N\}\]
and the zipper edges
\[ Z_{n+1}\ass1Z_n\cup\{(1p_N,10^n)\}, \]
which yields the rhombic strip
\[ R_{n+1}\ass0R^\times\cup 1R_\times\cup E\cup Z_{n+1}. \]
The new Hamiltonian path is
\begin{equation}
\label{eq:Pnp1}
P_{n+1}\ass(0p_1,1p_1,1p_2,0p_2,0p_3,1p_3,\ldots,0p_{N-1},1p_{N-1},1p_N,0p_N),
\end{equation}
which becomes $x$-monotone by slight perturbations of the vertical edges.
Furthermore, the new leftmost and rightmost chains are
\[ C_{n+1}\ass(0^{n+1},1C_n) \text{  and  } D_{n+1}\ass(0p_N,1D_n), \]
respectively (with the exception of $0^{n+1}$, which is not part of $D_{n+1}$), and the sets of edges above and below~$P_{n+1}$ are
\begin{equation}
\label{eq:ABp1}
\begin{split}
 A_{n+1}&\ass1A_n\cup\{(1p_i,1p_{i+1})\mid i=2,4,6,\ldots,N-2\} \text{ and } \\
 B_{n+1}&\ass0B_n\cup\{(0p_i,0p_{i+1})\mid i=1,3,5,\ldots,N-1\},
\end{split}
\end{equation}
respectively.
Note that in~\eqref{eq:Pnp1} and~\eqref{eq:ABp1} we have used that $N=2^n$ is even.

One can check that the vertices~$x_i$ and~$y_i$ on the chains~$C_n$ and~$D_n$ have the explicit form $x_i=1^i0^{n-i}$ for $i=0,\ldots,n$ and $y_i=1^{i-1}0^{n-i}1$ for $i=1,\ldots,n$, respectively.
\end{proof}

Maybe surprisingly, both aforementioned proofs actually yield the same rhombic strip, up to reversal of the vertex labels, which can be checked for $n=2,3,4,5$ by carefully comparing Figures~\ref{fig:Q-RS1b} and~\ref{fig:Q-RS2b}, and which can easily be shown in general using induction.

\begin{figure}[b!]
\makebox[0cm]{ 
\includegraphics{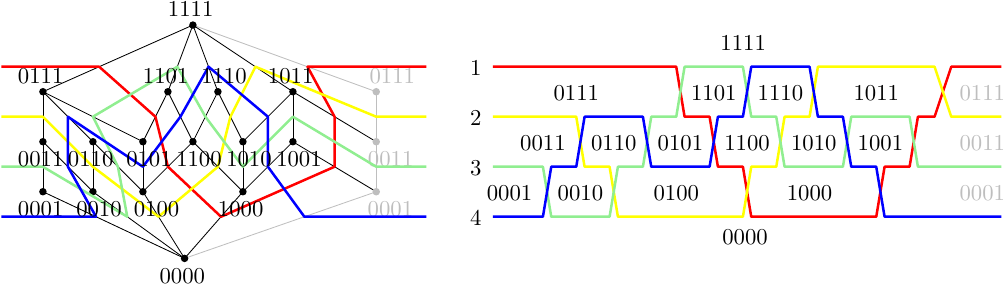}
}
\caption{The dual graph of a rhombic strip of~$G(Q_n)$ is a monotone $n$-Venn diagram.
Vertices and regions are encoded by their characteristic vectors of length~$n$.}
\label{fig:venn}
\end{figure}

There is an interesting connection between rhombic strips in~$G(Q_n)$ and Venn diagrams.
Specifically, the dual graph of a rhombic strip of~$G(Q_n)$ is an \defi{$n$-Venn diagram}~\cite{MR1668051}, i.e., a diagram of $n$ simple closed curves in the plane with only finitely many intersection points, such that each of the $2^n$ possible intersection patterns is represented by exactly one of the regions in the diagram; see Figure~\ref{fig:venn}.
Furthermore, those diagrams are \defi{monotone}, which means that for all $k=1,\ldots,n-1$ and every region inside of exactly $k$ of the curves, there is a neighboring region inside of exactly $k-1$ curves, and another neighboring region inside of exactly~$k+1$ curves (this is a consequence of Lemma~\ref{lem:strip-order}~(ii)).

The two aforementioned proofs can be dualized into the language of Venn diagrams, and amount to constructing a monotone $(n+1)$-Venn diagram inductively from a monotone $n$-Venn diagram.
Such constructions have already been provided in Venn's original paper~\cite{venn_1880} from~1880.

\section{Hypercubes}
\label{sec:cube}

\subsection{The polytope and its face lattice}

\begin{figure}[htb]
\makebox[0cm]{ 
\includegraphics[page=2]{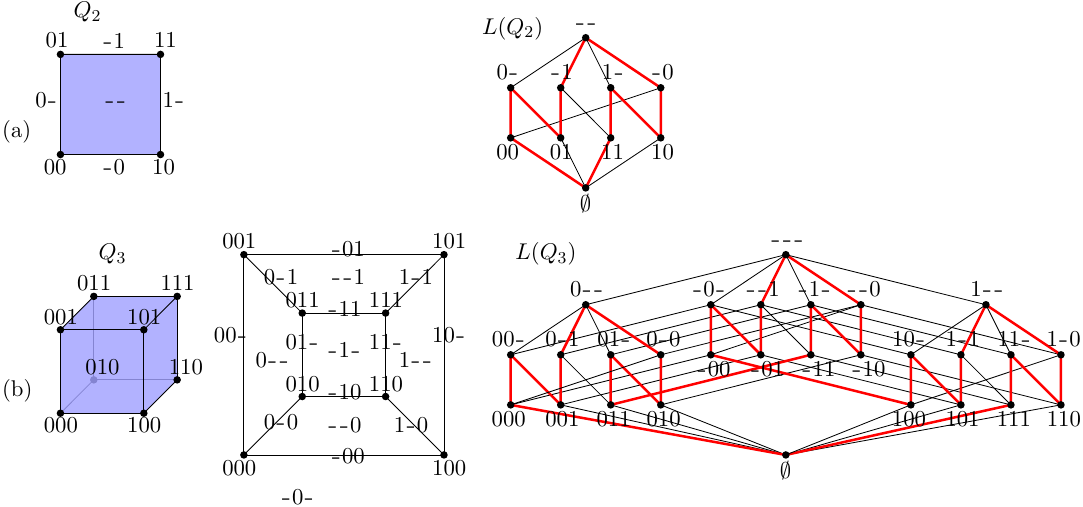}
}
\caption{The 2- and 3-cube and their face lattices, with a Hamiltonian cycle in the cover graph highlighted.
The cycles are the same as in Figure~\ref{fig:hc-cube-perm}~(a2) and~(a3).}
\label{fig:cube}
\end{figure}

The \defi{$n$-(hyper)cube} is the convex hull of the point set~$\{0,1\}^n$; see Figure~\ref{fig:cube}.
Equivalently, it is the $n$-fold Cartesian product of the unit interval~$[0,1]$.
The skeleton of this polytope is isomorphic to the cover graph of the Boolean lattice~$Q_n$, and so we use~$Q_n$ to also denote the polytope.
Each face~$F$ of~$Q_n$ can be encoded by a ternary string~$x_F\in\{0,1,\hyph\}^n$, where the set of vertices contained in~$F$ is given by replacing in~$x_F$ all occurrences of~$\hyph$ by either~0 or~1; see Figure~\ref{fig:cube}.
The dimension of the face~$F$ is given by the number of~$\hyph$s in the string~$x_F$.
For example, $x_F=01\hyph\hyph=\{0100,0101,0110,0111\}$ represents a 2-dimensional face~$F$ (a quadrilateral), and $\hyph^n=Q_n$ is the entire polytope.

The cover graph of the face lattice~$G(L(Q_n))$ has as vertex set~$\{0,1,\hyph\}^n\cup\{\emptyset\}$, and as edges all pairs of ternary strings that differ in replacing a single~0 or~1 by~$\hyph$, or vice versa, plus the edges~$(\emptyset,x)$ for all $x\in\{0,1\}^n$.

\subsection{Hamiltonian cycles and rhombic strips in the face lattice}

The following two theorems are the main results of this section.

\begin{theorem}
\label{thm:LQ-HC}
For any $n\geq 1$, the graph~$G(L(Q_n))$ has a Hamiltonian cycle.
\end{theorem}

The construction described in the following proof is illustrated in Figure~\ref{fig:cube} for $n=2,3$ and in Figure~\ref{fig:hc-cube-perm}~(a1)--(a4) for $n=1,2,3,4$.

\begin{proof}
We can list all $3^n$ strings from the set~$\{0,1,\hyph\}^n$ using the ternary reflected Gray code in such a way that any two consecutive strings differ in a flip~$0\leftrightarrow \hyph$ or $\hyph\leftrightarrow 1$.
This listing is defined inductively as~$\Gamma_1\ass(0,\hyph,1)$, and for $n\geq 2$ as
\[ \Gamma_n\ass(0\Gamma_{n-1},\hyph\rev(\Gamma_{n-1}),1\Gamma_{n-1}). \]
In words, the listing~$\Gamma_n$ consists of three copies of the previous listing~$\Gamma_{n-1}$, where all strings in the first copy are prefixed with~0, all strings in the second copy are prefixed with~$\hyph$, and all strings in the third copy are prefixed with~1, and the order of strings in the second copy is reversed.
It is easy to check that the sequence~$\Gamma_n$ starts with the string~$0^n$ and ends with the string~$1^n$, i.e., with two faces of rank~0 in the face lattice.
Consequently, $(\Gamma_n,\emptyset)$ is the desired Hamiltonian cycle in~$G(L(Q_n))$.
\end{proof}

\begin{theorem}
\label{thm:LQ-RS}
For any $n\geq 1$, the graph~$G(L(Q_n))$ has a rhombic strip.
\end{theorem}

For proving this result, we combine the techniques developed in the two proofs of Theorem~\ref{thm:Q-RS} presented in Section~\ref{sec:simplex}.
We use induction and embed three copies of the previous structures, one prefixed with~0, one prefixed with~$\hyph$, and one prefixed with~$1$.
The two structures prefixed with~0 and~$\hyph$ are connected as in the second proof via `stacking', and the two structures prefixed with~$\hyph$ and~1 are connected as in the first proof via `mirroring'.
An efficient Gray code algorithm implementing this construction in C++ can be found at~\cite{cos_sperm}.

\begin{proof}
We construct a rhombic strip~$R_n\seq G(L(Q_n))$ by induction on~$n$.
The notations used in this construction are illustrated in Figures~\ref{fig:LQ-RS1} and~\ref{fig:LQ-RS2}.
We define the abbreviation $B_n'\ass\{(\emptyset,x)\mid x\in\{0,1\}^n\}$ for the set of edges between the minimum~$\emptyset$ of~$L(Q_n)$ and all rank~0 faces.
The construction maintains a partition $R_n=P_n\cup A_n\cup B_n\cup B_n'\cup Z_n$ into five sets with the following properties:
\begin{itemize}[leftmargin=4mm]
\item $P_n$ is an $x$-monotone Hamiltonian path of~$G(L(Q_n))\setminus\{\emptyset\}$ that starts and ends at rank~0 vertices (red in the figures).
\item The sequence $(\emptyset,C_n)=(\emptyset,x_0,\ldots,x_n)$ of leftmost vertices across all ranks is a chain~$(\emptyset,C_n)\seq R_n$ (cyan in the figures), and the sequence $(\emptyset,D_n)=(\emptyset,y_0,\ldots,y_n)$ of rightmost vertices across all ranks is a chain~$(\emptyset,D_n)\seq R_n$ (pink in the figures), and the two chains are connected by the zipper edges $Z_n\ass\{(x_i,y_{i+1})\mid i=0,\ldots,n-2\}\seq R_n$ (in the gray-shaded area in the figures).
\item The edges $A_n$ and~$B_n\cup B_n'$ are the edges of~$R_n$ strictly above and below~$P_n$, respectively (in the blue-shaded and yellow-shaded areas in the figures).
\end{itemize}

\begin{figure}[t!]
\centering
\includegraphics[page=5]{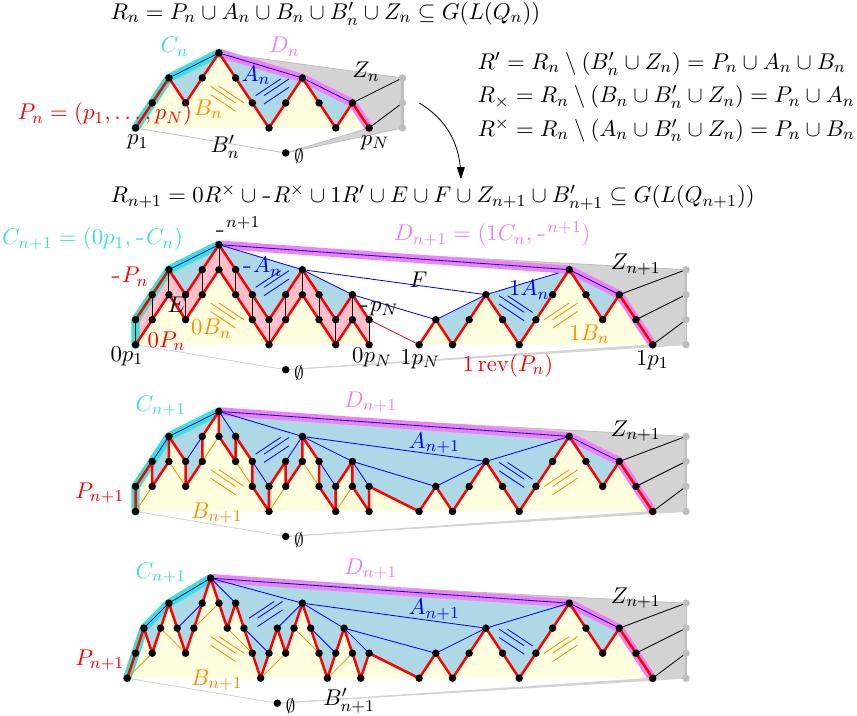}
\caption{Schematic illustration of the proof of Theorem~\ref{thm:LQ-RS}.}
\label{fig:LQ-RS1}
\end{figure}

\begin{figure}[htb]
\centering
\makebox[0cm]{ 
\includegraphics[page=6]{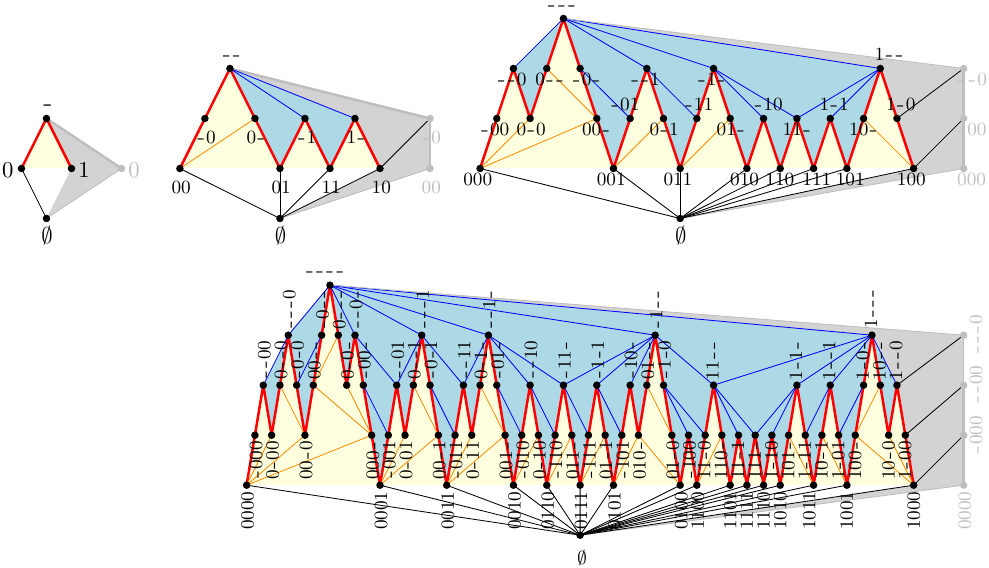}
}
\caption{The rhombic strips obtained from the proof of Theorem~\ref{thm:LQ-RS} for $n=1,2,3,4$.}
\label{fig:LQ-RS2}
\end{figure}

For the base case $n=1$ of the construction we take $R_1\ass G(L(Q_1))$ with $P_1\ass (0,\hyph,1)$, which gives $C_1=(0,\hyph)$, $D_1=(1,\hyph)$, $Z_1=\emptyset$, $A_1=\emptyset$ and~$B_1=\emptyset$.

For the induction step, let $P_n,C_n=(x_0,\ldots,x_n),D_n=(y_0,\ldots,y_n),A_n,B_n$ be given for some~$n\geq 1$.
We construct the rhombic strip~$R_{n+1}$ as follows:
We consider the plane graphs:
\begin{align*}
R^\times&\ass R_n\setminus (A_n\cup B_n'\cup Z_n)=P_n\cup B_n\\
R_\times&\ass R_n\setminus (B_n\cup B_n'\cup Z_n)=P_n\cup A_n\\
R'&\ass R_n\setminus (B_n'\cup Z_n)=P_n\cup A_n\cup B_n.
\end{align*}
We take a copy of~$0R^\times$ and a copy of~$\hyph R_\times$, maintaining the same embedding, and a vertically mirrored copy of~$1R'$.
The first two graphs are embedded at the same positions horizontally, but the second graph is placed one unit higher than the first one, so that the corresponding vertices on the two copies of~$P_n$ line up vertically (note that prefixing with $\hyph$ increases the ranks).
The third (mirrored) graph is embedded to the right of the first two and at the same height as the first one.
Let $P_n\assr(p_1,\ldots,p_N)$, $N=3^n$, be the sequence of vertices on the Hamiltonian path~$P_n$.
Then we connect the first and second graph by the edges
\[ E\ass\{(0p_i,\hyph p_i)\mid i=1,\ldots,N\},\]
the second and third graph by the edges
\[ F\ass\{(1y_i,\hyph y_i)\mid i=0,\ldots,n\}, \]
and the third and first graph by the zipper edges
\[ Z_{n+1}\ass\{(1x_i,\hyph x_i)\mid i=0,\ldots,n-1\}, \]
which yields the rhombic strip
\[ R_{n+1}\ass0R^\times\cup \hyph R_\times\cup 1R' \cup E\cup F\cup Z_{n+1}\cup B_{n+1}'. \]
The new Hamiltonian path (for $G(L(Q_{n+1}))\setminus\{\emptyset\}$) is
\begin{equation}
\label{eq:Pnp1n}
P_{n+1}\ass(0p_1,\hyph p_1,\hyph p_2,0p_2,0p_3,\hyph p_3,\ldots,\hyph p_{N-1},0p_{N-1},0p_N,\hyph p_N,1\rev(P_n)),
\end{equation}
which becomes $x$-monotone, again, by slight perturbations of the vertical edges in the first part.
Furthermore, the new leftmost and rightmost chains are
\[ C_{n+1}\ass(0p_1,\hyph C_n) \text{  and  } D_{n+1}\ass(1C_n,\hyph^{n+1}), \]
respectively (with the exception of $\emptyset$, which is neither part of $C_{n+1}$ nor~$D_{n+1}$), and the sets of edges above and below~$P_{n+1}$ (in addition to~$B_{n+1}'$) are
\begin{equation}
\label{eq:ABp1n}
\begin{split}
A_{n+1}&\ass\hyph A_n\cup 1A_n\cup\{(\hyph p_i,\hyph p_{i+1})\mid i=2,4,6,\ldots,N-1\}\cup (F\setminus\{(1p_N,\hyph p_N)\}) \text{ and } \\
B_{n+1}&\ass0B_n\cup 1B_n\cup\{(0p_i,0p_{i+1})\mid i=1,3,5,\ldots,N-2\},
\end{split}
\end{equation}
respectively.
Note that in~\eqref{eq:Pnp1n} and~\eqref{eq:ABp1n} we have used that $N=3^n$ is odd.

One can check that the vertices~$x_i$ and~$y_i$ on the chains~$C_n$ and~$D_n$, respectively, have the explicit form $x_i=\hyph^i0^{n-i}$ for $i=0,\ldots,n$ and $y_i=1\hyph^i0^{n-1-i}$ for $i=0,\ldots,n-1$.
\end{proof}

\section{Permutahedra}
\label{sec:perm}

\subsection{The polytope and its face lattice}

\begin{figure}[htb]
\centering
\includegraphics[page=1]{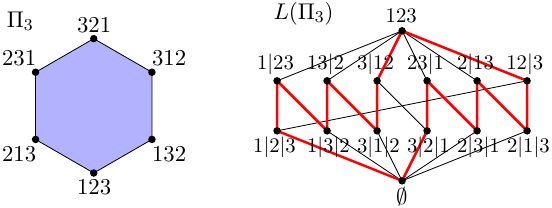}
\caption{The 2-dimensional permutahedron~$\Pi_3$ and its face lattice, with a Hamiltonian cycle in the cover graph highlighted.
The cycle is the same as in Figure~\ref{fig:hc-cube-perm}~(b3).}
\label{fig:perm}
\end{figure}

The \defi{permutahedron~$\Pi_n$} is the convex hull of $\{(\pi(1),\ldots,\pi(n))\mid \pi\in S_n\}$, where $S_n$ denotes the set of all permutations on~$[n]$; see Figure~\ref{fig:perm}.
The permutahedron is~$(n-1)$-dimensional, and its edges connect pairs of permutations that differ in an adjacent transposition.\footnote{By this we mean a transposition of adjacent positions, i.e., $\pi(i)\leftrightarrow\pi(i+1)$, which is true if vertices are labeled by the inverse of the permutation that defines the vertex coordinates in our definition via the convex hull.}
More generally, faces of the permutahedron can be described as follows:
an (ordered) \defi{set partition} of~$[n]$ is a sequence~$A_1,\ldots,A_k$ of nonempty disjoint subsets~$A_i\seq [n]$ whose union is~$[n]$, where $k\in[n]$.
We write such a partition as a string $A_1|A_2|\cdots|A_k$, and we refer to the sets~$A_i$, $i\in[k]$, as \defi{blocks}, and to the characters~$|$ as \defi{bars}.
Each $k$-dimensional face~$F$ of~$\Pi_n$ corresponds to a set partition~$x_F=A_1|\cdots|A_{n-k}$ with $n-k$ blocks, where the set of permutations contained in~$F$ is given by permuting the elements in each block of the partition~$x_F$ arbitrarily.
For example~$x_F=25|1|34$ encodes the 2-dimensional face~$F$ of~$\Pi_5$ that contains the permutations $\{25134,52134,25143,52143\}$ (a quadrilateral).
The facets of~$\Pi_n$ are set partitions with two blocks~$A_1|A_2$, i.e., these can be identified with nonempty proper subsets~$A_1\seq [n]$, and there are~$2^n-2$ of them.
The total number of faces of the permutahedron (excluding~$\emptyset$) is given by the Fubini numbers or ordered Bell numbers (OEIS~A000670).

The cover graph~$G(L(\Pi_n))$ has as vertex set all set partitions of~$[n]$, and edges between pairs of partitions that differ in removing a bar or adding a bar, which results in joining two blocks to one or splitting a block into two, respectively, plus the edges~$(\emptyset,\pi)$ for all $\pi\in S_n$.

\subsection{Hamiltonian cycles in the face lattice}

The following theorem and the analogous result for $B$-permutahedra (stated in Theorem~\ref{thm:LPib-HC} below) are the two main results of this section.

\begin{theorem}
\label{thm:LPi-HC}
For any $n\geq 2$, the graph $G(L(\Pi_n))$ has a Hamiltonian cycle.
\end{theorem}

The Hamiltonian cycles constructed in this proof are shown in Figure~\ref{fig:hc-cube-perm}~(b2)--(b4) for $n=2,3,4$.
An efficient Gray code algorithm implementing this construction in C++ can be found at~\cite{cos_part}.

\begin{proof}
For any set partition~$x=A_1|A_2|\cdots|A_k$ of~$[n-1]$ and for $i=0,\ldots,k$, we let $\wc{c}_i(x)$ be the set partition of~$[n]$ obtained by adding the singleton set~$\{n\}$ as a new block after~$A_i$ and before~$A_{i+1}$ (at the beginning if~$i=0$ and the end if~$i=k$), i.e.,
\[ \wc{c}_i(x)\ass A_1|\cdots|A_i|\{n\}|A_{i+1}|\cdots|A_k. \]
Furthermore, for $i=1,\ldots,k$, we let $\wh{c}_i(x)$ be the set partition of~$[n]$ obtained by adding the element~$n$ to the set~$A_i$, i.e.,
\[ \wh{c}_i(x)\ass A_1|\cdots|A_{i-1}|A_i\cup\{n\}|A_{i+1}|\cdots|A_k. \]
We then define a sequence~$\rvec{c}(x)$ of set partitions of~$[n]$ by
\[ \rvec{c}(x)\ass\big(\wc{c}_0(x),\wh{c}_1(x),\wc{c}_1(x),\wh{c}_2(x),\wc{c}_2(x),\ldots,\wh{c}_{k-1}(x),\wc{c}_{k-1}(x),\wh{c}_k(x),\wc{c}_k(x)\big).\]
In words, $\rvec{c}(x)$ is obtained from the set partition~$x$ by either adding~$\{n\}$ as a new singleton block between two existing blocks, or by joining the element~$n$ with one of the existing blocks in~$x$, alternatingly from left to right.
Consequently, the sequence~$\rvec{c}(x)$ has length~$2k+1$, and consists alternatingly of $k$-faces and $(k+1)$-faces of~$\Pi_{n}$.
In fact, any two consecutive faces in the sequence~$\rvec{c}(x)$ form an edge in~$G(L(\Pi_{n+1}))$, i.e., the sequence describes a path in~$G(L(\Pi_n))$ that alternates between rank~$r$ and~$r+1$, where $r=n-1-k$ is the rank of~$x$ in~$L(\Pi_{n-1})$.
We write $\lvec{c}(x)\ass\rev(\rvec{c}(x))$ for the reverse sequence/path.
For example, for $x=25|1|34\in L(\Pi_5)$ we have the paths
\begin{align*}
\rvec{c}(x)&=(\redsix|25|1|34,25\redsix|1|34,25|\redsix|1|34,25|1\redsix|34,25|1|\redsix|34,25|1|34\redsix,25|1|34|\redsix) \text{ and} \\
\lvec{c}(x)&=(25|1|34|\redsix,25|1|34\redsix,25|1|\redsix|34,25|1\redsix|34,25|\redsix|1|34,25\redsix|1|34,\redsix|25|1|34)
\end{align*}
in~$G(L(\Pi_6))$.

To prove the theorem, we construct a path~$P_n$ in~$G(L(\Pi_n))$ that visits all faces of~$\Pi_n$ except~$\emptyset$ and that starts and ends at rank~0 faces (i.e., two permutations), and therefore $(P_n,\emptyset)$ is the desired Hamiltonian cycle in~$G(L(\Pi_n))$.

The path~$P_n$ is constructed inductively as follows:
For the base case $n=2$ of the construction we take $P_2\ass1|2,12,2|1$.

For the induction step, let $P_{n-1}\assr(x_1,\ldots,x_N)$ be the path in~$G(L(\Pi_{n-1}))$.
Recall that $L(\Pi_{n-1})$ is bipartite, with $P_{n-1}$ starting and ending in the same part, so $N$ is always odd.
We define
\begin{equation}
\label{eq:PPi}
P_n\ass(\lvec{c}(x_1),\rvec{c}(x_2),\lvec{c}(x_3),\rvec{c}(x_4),\ldots,\lvec{c}(x_N)).
\end{equation}
Note that both the first and last entry in each of the subsequences~$\lvec{c}(x_1)$ and~$\lvec{c}(x_N)$ have rank~0, and so~$P_n$ starts and ends with rank~0 faces.
We observe that if $x$ and~$y$ are adjacent in~$G(L(\Pi_{n-1}))$, then the last entry of~$\lvec{c}(x)$ and the first entry of~$\rvec{c}(y)$ (which are $\{n\}|x$ and $\{n\}|y$, respectively), are adjacent in~$G(L(\Pi_n))$.
Similarly the last entry of~$\rvec{c}(x)$ and the first entry of~$\lvec{c}(y)$ (which are $x|\{n\}$ and~$y|\{n\}$, respectively) are adjacent.
From this it follows inductively that $P_n$ as defined in~\eqref{eq:PPi} is indeed a path in~$G(L(\Pi_n))$ that visits all faces of~$\Pi_n$ except~$\emptyset$.
This completes the proof.
\end{proof}

One can easily prove by induction the following two noteworthy properties of the path~$P_n$ in~$G(L(\Pi_n))$ constructed in the proof of Theorem~\ref{thm:LPi-HC}.
Firstly, for $P_n\assr(x_1,\ldots,x_N)$ and all $i=1,\ldots,N$, if $x_i=A_1|A_2|\cdots|A_k$, then we have $x_{N+1-i}=A_k|\cdots|A_2|A_1$, i.e., the set partitions~$x_i$ and~$x_{N+1-i}$ differ only in reversing the order of their blocks.
Secondly, let us consider the subsequence~$P_n'$ of~$P_n$ given by all permutations, i.e., set partitions with $n$ blocks.
For a permutation~$x\in S_{n-1}$, let $\lvec{c}'(x)$ denote the subsequence of~$\lvec{c}(x)$ of permutations from~$S_n$, i.e., $\lvec{c}'(x)=(\wc{c}_{n-1}(x),\wc{c}_{n-2}(x),\ldots,\wc{c}_1(x),\wc{c}_0(x))$.
Then for $P_{n-1}'\assr(x_1,\ldots,x_N)$ we have $P_n'=(\lvec{c}'(x_1),\lvec{c}'(x_2),\lvec{c}'(x_3),\ldots,\lvec{c}'(x_N))$, i.e., the next sequence of permutations is obtained by repeatedly inserting the new largest value~$n$ from right to left in the previous list of permutations.
This of course does not yield a Hamiltonian path in~$G(\Pi_n)$.

\subsection{$B$-permutahedra}

\begin{figure}[b!]
\centering
\includegraphics[page=2]{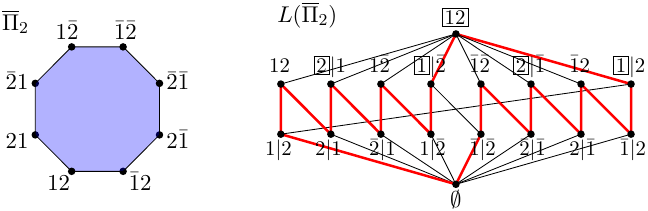}
\caption{The 2-dimensional $B$-permutahedron~$\overline{\Pi}_2$ and its face lattice, with a Hamiltonian cycle in the cover graph highlighted.
Negative signs are indicated by overlining, and boxes indicate type~2 faces, i.e., both signs for the corresponding entries.
The cycle is the same as the one constructed in the proof of Theorem~\ref{thm:LPib-HC}.}
\label{fig:bperm2}
\end{figure}

The \defi{$B$-permutahedron~$\overline{\Pi}_n$} is the convex hull of all signed permutations, i.e., the convex hull of $\{(\pm \pi(1),\ldots,\pm\pi(n))\mid \pi\in S_n\}$; see Figure~\ref{fig:bperm} and Figure~\ref{fig:bperm2}.

The $B$-permutahedron is \mbox{$n$-dimensional}.
Its edges connect pairs of permutations that either differ in an adjacent transposition, preserving all signs, or in a complementation of the sign of the first entry.
Each $k$-dimensional face~$F$ of~$\overline{\Pi}_n$, $k\geq 1$, corresponds to a signed set partition of one of the following two types:
\begin{enumerate}[leftmargin=15mm]
\item[(type 1)] a set partition $x_F=A_1|A_2|\cdots|A_{n-k}$ of~$[n]$ with $n-k$ blocks where each $i\in[n]$ carries a positive or negative sign;
\item[(type 2)] a set partition $x_F=\boxed{A_1}|A_2|\cdots|A_{n-k+1}$ of~$[n]$ with $n-k+1$ blocks where each $i\in[n]\setminus A_1$ carries a positive or negative sign, and all $i\in A_1$ carry both signs (positive and negative).
\end{enumerate}
For type~1, the set of signed permutations contained in~$F$ is given by permuting all elements with their corresponding signs in each block of~$x_F$ arbitrarily.
For type~2, one also iterates over all ways of selecting one of the two possible signs for all elements in the first block~$A_1$.

A type~1 face $A_1|A_2|\cdots|A_k$ has the cover relations
\begin{align*}
A_1|A_2|\cdots|A_k &\;\cov\; \boxed{A_1}|A_2|\cdots|A_k, \quad \text{and} \\
A_1|A_2|\cdots|A_k &\;\cov\; A_1|\cdots|A_i\cup A_{i+1}|\cdots|A_k \quad \text{ for } i=1,\ldots,k-1,
\end{align*}
whereas a type~2 face $\boxed{A_1}|A_2|\cdots|A_{k}$ has the cover relations
\begin{align*}
\boxed{A_1}|A_2|\cdots|A_{k} &\;\cov\; \boxed{A_1\cup A_2}|\cdots|A_{k}, \quad \text{and} \\ \boxed{A_1}|A_2|\cdots|A_{k} &\;\cov\; \boxed{A_1}|A_2|\cdots|A_i\cup A_{i+1}|\cdots|A_k \quad \text{ for } i=2,\cdots,k-1.
\end{align*}

\begin{theorem}
\label{thm:LPib-HC}
For any $n\geq 1$, the graph $G(L(\overline{\Pi}_n))$ has a Hamiltonian cycle.
\end{theorem}

The Hamiltonian cycle constructed in this proof is shown in Figure~\ref{fig:bperm2} for~$n=2$.

\begin{proof}
For a type~1 signed set partition $x=A_1|\cdots|A_k$ of~$[n-1]$ we define
\begin{align*}
\wc{c}_i^+(x)&\ass A_1|\cdots|A_i|\{n\}|A_{i+1}|\cdots|A_k, \\
\wc{c}_i^-(x)&\ass A_1|\cdots|A_i|\{\overline{n}\}|A_{i+1}|\cdots|A_k, \quad i=0,\ldots,k, \\
\wh{c}_i^+(x)&\ass A_1|\cdots|A_{i-1}|A_i\cup\{n\}|A_{i+1}|\cdots|A_k, \\
\wh{c}_i^-(x)&\ass A_1|\cdots|A_{i-1}|A_i\cup\{\overline{n}\}|A_{i+1}|\cdots|A_k, \quad i=1,\ldots,k,
\end{align*}
where the overline indicates a negative sign, and we further define
\[ \wh{c}_0^\pm(x)\ass\boxed{\{n\}}|A_1|\cdots|A_k. \]
For a type~2 signed set partition $y=\boxed{A_1}|A_2|\cdots|A_k$ of~$[n-1]$ we define $\wc{c}_i^+(y)$, $\wc{c}_i^-(y)$ for $i=1,\ldots,k$ as before, and $\wh{c}_i^+(y)$, $\wh{c}_i^-(y)$ for $i=2,\ldots,k$ as before, plus the additional
\[ \wh{c}_1^\pm(y)\ass \boxed{A_1\cup\{n\}}|A_2|\cdots|A_{k}. \]
In each of these two cases, we define a sequence~$c(x)$ and~$c(y)$ of signed set partitions of~$[n]$ by
\begin{align*}
c(x)&\ass(\wc{c}_k^+(x),\wh{c}_k^+(x),\ldots,\wh{c}_1^+(x),\wc{c}_0^+(x),\wh{c}_0^\pm(x),\wc{c}_0^-(x),\wh{c}_1^-(x),\ldots,\wh{c}_k^-(x),\wc{c}_k^-(x)) \text{ and} \\
c(y)&\ass(\wc{c}_k^+(y),\wh{c}_k^+(y),\ldots,\wh{c}_2^+(y),\wc{c}_1^+(y),\wh{c}_1^\pm(y),\wc{c}_1^-(y),\wh{c}_2^-(y),\ldots,\wh{c}_k^-(y),\wc{c}_k^-(y)),
\end{align*}
respectively.
Note that both sequences define paths in~$G(L(\overline{\Pi}_{n}))$.

To prove the theorem, we construct a path~$P_n$ in~$G(L(\overline{\Pi}_n))$ that visits all faces of~$\overline{\Pi}_n$ except~$\emptyset$ and that starts and ends at rank~0 faces (i.e., two signed permutations), and therefore~$(P_n,\emptyset)$ is the desired Hamiltonian cycle in~$G(L(\overline{\Pi}_n))$.
The path~$P_n$ is constructed inductively as follows: For the base case~$n=1$ of the construction we take $P_1\ass 1,\boxed{1},\overline{1}$.
For the induction step, let $P_{n-1}\assr (x_1,\ldots,x_N)$ be the path in~$G(L(\overline{\Pi}_{n-1}))$.
Then we define
\[ P_n\ass (c(x_1),\rev(c(x_2)),c(x_3),\ldots,c(x_N)). \]
It can be checked straightforwardly that the path~$P_n$ has the required properties.
\end{proof}

\begin{figure}[b!]
\makebox[0cm]{ 
\includegraphics{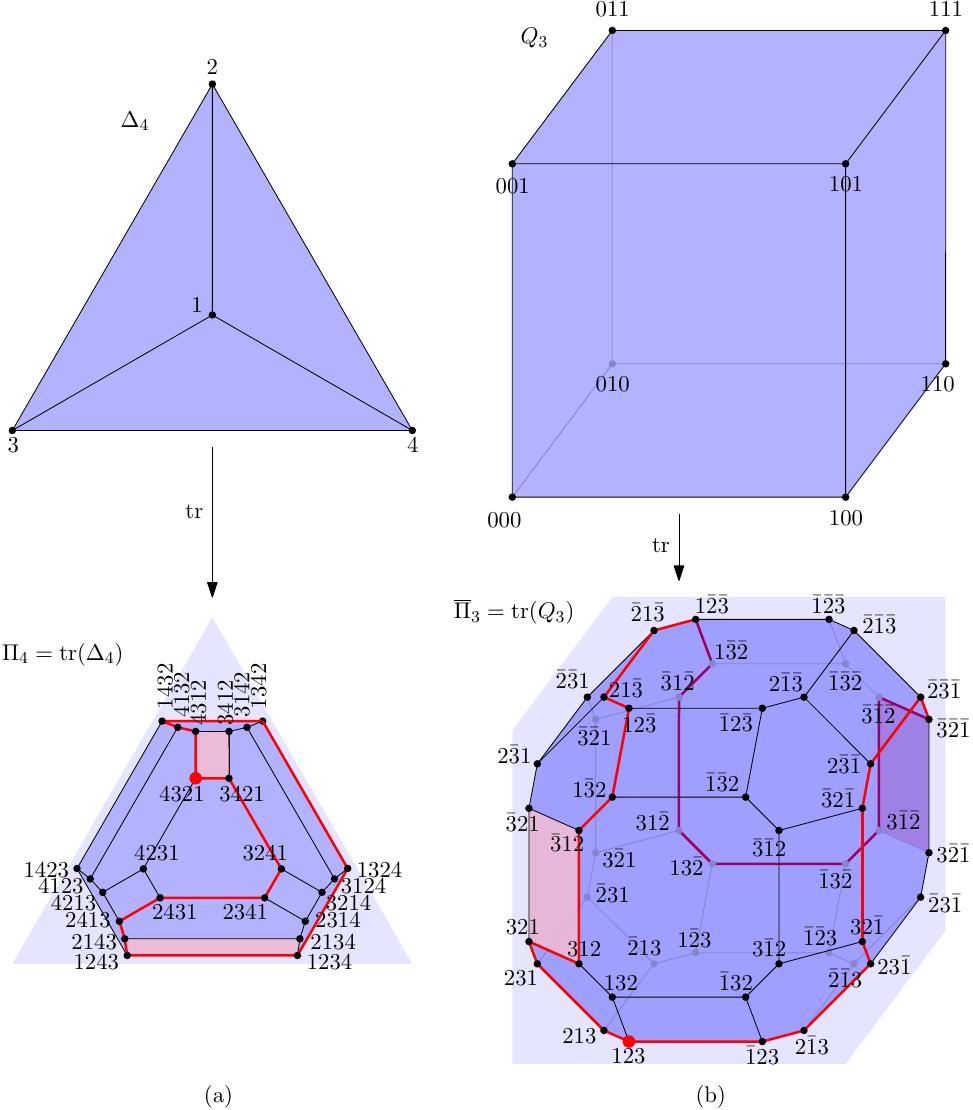}
}
\caption{(a) Truncating the simplex yields the permutahedron. (b) Truncating the hypercube yields the $B$-permutahedron.
The facet-Hamiltonian cycles in the permutahedron and B-permutahedron are obtained via Lemma~\ref{lem:truncate} from the rhombic strips shown in Figures~\ref{fig:Q-RS1b} and~\ref{fig:LQ-RS2}, respectively.
The highlighted vertex corresponds to the leftmost chain in the rhombic strip.
Each of the highlighted four-cycles yields two possibilities for the facet-Hamiltonian cycle to walk around this face, corresponding to two options of sweeping the chain through the rhombic strip.}
\label{fig:truncate}
\end{figure}

\section{Truncation}
\label{sec:truncate}

Given a polytope~$P$, the \defi{omnitruncation} of~$P$, denoted~$\tr(P)$, is the simple polytope of the same dimension as~$P$ that has one vertex for each flag of~$P$.
Geometrically, it is obtained by truncating all the (non-trivial) faces of~$P$; see Figure~\ref{fig:truncate}.

It is well-known that truncating the simplex yields the permutahedron, i.e., we have $\tr(\Delta_n)=\Pi_n$, and truncating the hypercube yields the $B$-permutahedron, i.e., we have $\tr(Q_n)=\overline{\Pi}_n$.

In their paper~\cite{MR4863586}, the authors established the following connection between rhombic strips in the face lattice~$L(P)$ of a polytope~$P$ and facet-Hamiltonian cycles in the omnitruncated polytope~$\tr(P)$.

\begin{lemma}[\cite{MR4863586}]
\label{lem:truncate}
The graph~$G(\tr(P))$ has a facet-Hamiltonian cycle if and only if $G(L(P))$ has a rhombic strip.
\end{lemma}

The facet-Hamiltonian cycle in~$G(\tr(P))$ is obtained by sweeping a chain from left-to-right through the rhombic strip of~$G(L(P))$, changing one element of the chain at a time, moving it across a diamond.
If several elements in the chain can be changed, then we have freedom in constructing several different facet-Hamiltonian cycles; see Figure~\ref{fig:truncate}.

Applying Lemma~\ref{lem:truncate} for $P\ass\Delta_n$ and using Theorem~\ref{thm:Q-RS} thus proves that the permutahedron~$\Pi_n$ has a facet-Hamiltonian cycle; see Figure~\ref{fig:truncate}~(a).
In fact, we obtain many different such cycles, corresponding to the different constructions discussed in Section~\ref{sec:simplex} that establish Theorem~\ref{thm:Q-RS}.

Similarly, applying Lemma~\ref{lem:truncate} for $P\ass Q_n$ and using Theorem~\ref{thm:LQ-RS} yields the following result, which affirmatively resolves Conjecture~1 raised by Akitaya, Cardinal, Felsner, Kleist and Lauff~\cite{MR4863586}; see Figure~\ref{fig:truncate}~(b).

\begin{theorem}
\label{thm:BPi-fHC}
For any $n\geq 2$, the $B$-permutahedron $G(\overline{\Pi}_n)$ has a facet-Hamiltonian cycle.
\end{theorem}

\section{Associahedra}
\label{sec:asso}

\subsection{The polytope and its face lattice}

\begin{figure}[htb]
	\centering
	\includegraphics{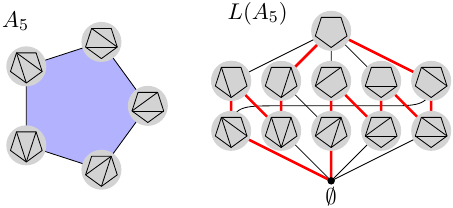}
	\caption{The 2-dimensional associahedron~$A_5$ and its face lattice, with a Hamiltonian cycle in the cover graph highlighted.
		The cycle is the same as in Figure~\ref{fig:hc-asso} for $n=5$.}
	\label{fig:asso}
\end{figure}

The combinatorial model of the associahedron that we will use in this paper does not use binary trees, as mentioned in the introduction, but instead the corresponding dual graphs, namely triangulations of a convex $n$-gon; see Figure~\ref{fig:3poly}~(c).
Specifically, the \defi{associahedron $A_n$} is the polytope whose face lattice~$L(A_n)$ is the refinement order of dissections of a convex $n$-gon; see Figure~\ref{fig:asso}.
This polytope is $(n-3)$-dimensional.
The vertices of~$A_n$ are in bijection with triangulations of a convex $n$-gon, and edges correspond to flips in a triangulation, where a flip removes an inner edge between two triangles, and replaces it by the other diagonal of the resulting quadrilateral.
More generally, the $k$-dimensional faces are given by all dissections of the $n$-gon with exactly $(n-k-3)$ inner edges.
In the face lattice~$L(A_n)$, the cover relations are between dissections that differ in adding or removing a single (inner) edge.
The total number of faces of the associahedron (excluding $\emptyset$) is given by the little Schr\"oder numbers (OEIS A001003).

\subsection{Hamiltonian cycles in the face lattice}

The following is the main result of this section.

\begin{theorem}
\label{thm:LAsso-HC}
For any $n\geq 4$, the graph $G(L(A_n))$ has a Hamiltonian cycle.
\end{theorem}

The Hamiltonian cycles constructed in this proof are shown in Figure~\ref{fig:hc-asso} for $n=3,4,5,6$.
An efficient Gray code algorithm implementing this construction in C++ can be found at~\cite{cos_kary}.

\begin{proof}
We label the points $1,\ldots,n$ in counterclockwise order.
For a given dissection~$X$, we consider the set of edges incident with the point~$n$, and we let $v_1,\ldots,v_k$ be the neighbors of~$n$ in increasing order.
Clearly, we have $k\geq 2$ and $v_1=1$ and~$v_k=n-1$.
For $i=2,\ldots,k$, we let $\wh{c}_i(X)$ be the dissection of the $(n+1)$-gon obtained from~$X$ by expanding the vertex~$n$ into two vertices~$n$ and~$n+1$ connected by an edge, where $n+1$ retains $v_1,\ldots,v_{i-1}$ as neighbors, and $n$ retains $v_i,\ldots,v_k$ as neighbors; see Figure~\ref{fig:asso-ins}~(a).
Furthermore, for $i=2,\ldots,k$, we let $\wc{c}_i(X)$ be the dissection obtained from~$\wh{c}_i(X)$ by adding the edge~$(v_i,n+1)$.
Note that $\wc{c}_i(X)$ for $i=2,\ldots,k-1$ is also obtained from~$\wh{c}_{i+1}(X)$ by adding the edge~$(v_i,n)$.
Lastly, we let $\wc{c}_1(X)$ be the dissection obtained from~$\wh{c}_2(X)$ by adding the edge~$(v_1,n)=(1,n)$.

\begin{figure}[b!]
\includegraphics{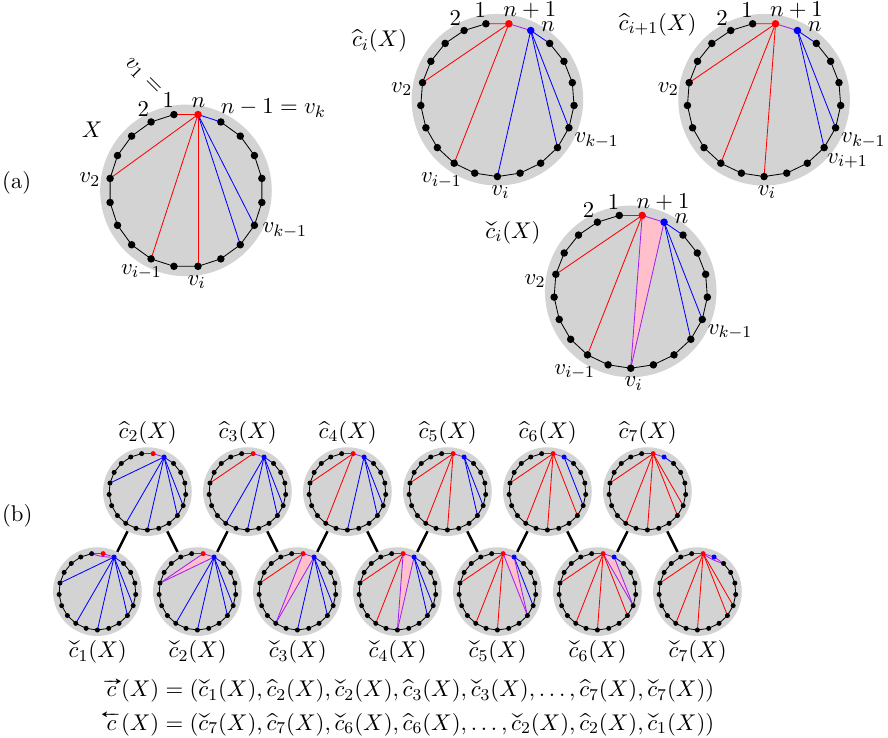}
\caption{Illustration of the proof of Theorem~\ref{thm:LAsso-HC}.
Only edges incident with~$n$ or~$n+1$ are shown in the figure, whereas all others are omitted for clarity.}
\label{fig:asso-ins}
\end{figure}

We define the sequence~$\rvec{c}(X)$ of dissections of an $(n+1)$-gon by
\[ \rvec{c}(X)\ass\big(\wc{c}_1(X),\wh{c}_2(X),\wc{c}_2(X),\wh{c}_3(X),\wc{c}_3(X),\ldots,\wh{c}_{k-1}(X),\wc{c}_{k-1}(X),\wh{c}_k(X),\wc{c}_k(X)\big), \]
and we note that it has length~$2k-1$ and describes a path in~$G(L(A_{n+1}))$ that alternates between ranks~$r$ and~$r+1$, where $r$ is the rank of~$X$ in~$L(A_n)$; see Figure~\ref{fig:asso-ins}~(b).
We write $\lvec{c}(X)\ass\rev(\lvec{c}(X))$ for the reverse sequence/path.

To prove the theorem, we construct a path~$P_n$ in~$G(L(A_n))$ that visits all faces of~$A_n$ except~$\emptyset$ and that starts and ends at rank~0 faces (i.e., two triangulations), and therefore~$(P_n,\emptyset)$ is the desired Hamiltonian cycle in~$G(L(A_n))$.

The path~$P_n$ is constructed inductively as follows:
For the base case~$n=4$ of the construction we take~$P_4$ consisting of the three dissections of a 4-gon shown in Figure~\ref{fig:hc-asso}.

For the induction step, let $P_n\assr(X_1,\ldots,X_N)$ be the path in~$G(L(A_n))$, i.e., each $X_i$ is a dissection of an $n$-gon.
Then we define
\begin{equation}
\label{eq:PAsso}
P_{n+1}\ass(\lvec{c}(X_1),\rvec{c}(X_2),\lvec{c}(X_3),\rvec{c}(X_4),\ldots,\lvec{c}(X_N)).
\end{equation}
Note that both the first and last entry in each of the subsequences~$\lvec{c}(X_1)$ and~$\lvec{c}(X_N)$ have rank~0, and so~$P_{n+1}$ starts and ends with rank~0 faces.
We observe that if $X$ and~$Y$ are adjacent in~$G(L(A_n))$, then the last entry of~$\lvec{c}(X)$ and the first entry of~$\rvec{c}(Y)$ are adjacent in~$G(L(A_{n+1}))$, and similarly the last entry of~$\rvec{c}(X)$ and the first entry of~$\lvec{c}(Y)$ are adjacent.
From this it follows inductively that $P_{n+1}$ as defined in~\eqref{eq:PAsso} is indeed a path in~$G(L(A_{n+1}))$ that visits all faces of~$A_{n+1}$ except~$\emptyset$.
This completes the proof.
\end{proof}

\section{Cyclic polytopes}
\label{sec:cyclic}

Cyclic polytopes play an important role in polyhedral combinatorics.
According to the upper bound theorem, they maximize the number of faces for a given dimension and number of vertices.
Let us recall their definition.

\subsection{The polytope and its face lattice}

The \defi{moment curve} in~$\mathbb{R}^d$ is defined as $f:\mathbb{R}\rightarrow \mathbb{R}^d$, $t\mapsto f(t)=(t,t^2,\ldots,t^d)$.
The \defi{cyclic polytope} in~$\mathbb{R}^d$ on $n$ vertices, denoted $C(n,d)$ for $n>d\geq 2$, is defined as the convex hull of $n$ consecutive points $f(1),\ldots,f(n)$ on the moment curve.
The polytope~$C_{n,d}$ is $d$-dimensional, and for $n=d+1$ it is isomorphic to the simplex~$\Delta_n$.
Each face of~$C_{n,d}$ corresponds to a subset~$F\seq[n]$, namely $i\in F$ if and only if the corresponding point~$f(i)$ belongs to the face.
It is known that $C_{n,d}$ is a \defi{simplicial} polytope, i.e., all of its facets are $(d-1)$-dimensional simplices (on $d$ points).

To describe the face lattice of cyclic polytopes, it essential to describe the facets.
As the polytopes are simplicial, all other faces are then obtained simply by taking all possible subsets of the facets.
The facets of cyclic polytopes are given by \defi{Gale's evenness conditions} \cite[Thm.~0.7]{MR1311028}.
Specifically, a subset $F=\{i_1<i_2<\cdots<i_d\}\seq [n]$ is a facet of~$C_{n,d}$ if and only if any two consecutive elements in~$[n]\setminus F$ are separated by an even number of elements in~$F$.
In other words, any maximal subsequence of consecutive elements in~$F$ that excludes~1 and~$n$ has even length.

\subsection{Description as a subposet of the Boolean lattice}

\begin{figure}[h!]
\includegraphics[page=1]{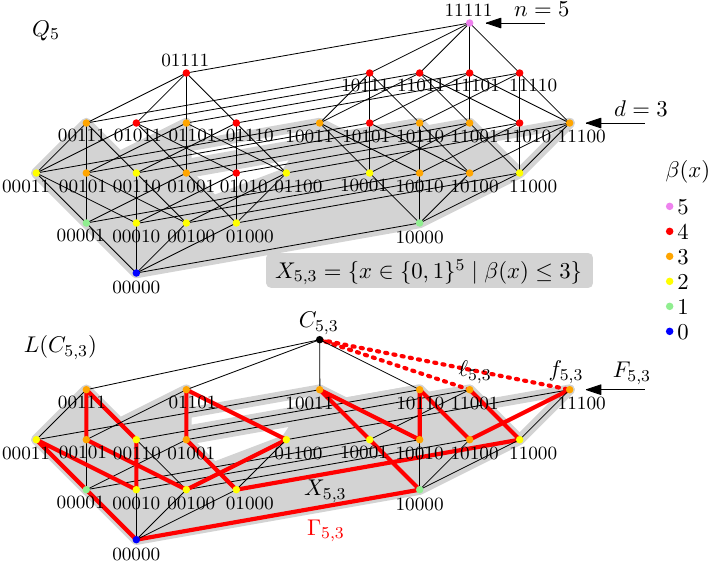}
\caption{Illustration of Lemma~\ref{lem:cyclic-downset}.
The Hamiltonian cycle highlighted in the face lattice~$L(C_{5,3})$ is obtained from Lemma~\ref{lem:brgc-cyclic}.}
\label{fig:downset}
\end{figure}

In the following, we describe each face~$F\seq [n]$ of~$C_{n,d}$ by its corresponding characteristic vector, i.e., a bitstring~$x_F$ of length~$n$ with a 1-bit at position~$i$ if $i\in F$ and a 0-bit at position~$i$ if $i\notin F$.
For any $x\in\{0,1\}^n$ we let $w(x)$ be the Hamming weight of~$x$.
We refer to a sequence of~1s in~$x$ that is preceded and followed by a 0-bit as a \defi{run}, and to the number of 1s in the sequence as the \defi{length} of the run.
Note that a sequence of 1s starting at the first position or ending at the last position is \emph{not} a run.
We say that a run is \defi{even} or \defi{odd}, depending on the parity of its length, and we write $o(x)$ for the number of odd runs in~$x$.
For example, $x=110001111001000$ contains two runs, one even (of length~4) and one odd (of length~1), and therefore $o(x)=1$.
By Gale's evenness condition,
\[ F_{n,d}\ass\{x\in \{0,1\}^n\mid w(x)=d\wedge o(x)=0\} \]
is the set of facets of~$C_{n,d}$.
As $C_{n,d}$ is simplicial, the face lattice~$L(C_{n,d})$ (without the trivial face at the top) is isomorphic to the downset of~$F_{n,d}$ in the Boolean lattice~$Q_n$.
We can describe this downset combinatorially, as follows.
For any $x\in\{0,1\}^n$ we introduce the parameter $\beta(x)\ass w(x)+o(x)$, which is the sum of the Hamming weight of~$x$ and the number of odd runs in~$x$.
For example, for $x=110001111001000$ from before we have $w(x)=7$ and $o(x)=1$, and therefore $\beta(x)=7+1=8$.
The following lemma is illustrated in Figure~\ref{fig:downset}.

\begin{lemma}
\label{lem:cyclic-downset}
The face lattice~$L(C_{n,d})$ without the trivial face at the top is isomorphic to the subposet of~$Q_n$ induced by the set
\begin{equation}
\label{eq:Xnd}
X_{n,d}\ass\{x\in\{0,1\}^n\mid \beta(x)\leq d\}.
\end{equation}
\end{lemma}

In words, characteristic vectors of faces of~$C_{n,d}$ are precisely the bitstrings whose $\beta$-value is bounded from above by~$d$.
Note that for all $x\in X_{n,d}$ with $w(x)=d$ we have $o(x)=0$, and more generally if $w(x)=d-s$ then we have $o(x)\leq s$.

\begin{proof}
Firstly, we need to show that if $x\in F_{n,d}$, then any $y$ in the downset of~$x$ satisfies $\beta(y)\leq d$.
For this it is enough to prove the following claim:
Let $x'\in\{0,1\}^n$ and suppose that $y'$ is obtained from~$x'$ by flipping a single 1-bit to a 0-bit, then we have $\beta(y')\leq \beta(x')$.
Clearly, we have $w(y')=w(x')-1$.
If the flipped 1-bit is contained in an even run in~$x'$, then in~$y'$ this will create one odd run and at most one even run.
On the other hand, if the 1-bit is contained in an odd run in~$x'$, then in~$y'$ this will create either at most two odd runs or at most two even runs.
Lastly, if the 1-bit is part of a maximal leading or terminal sequence of 1s in~$x'$, then this will create at most one even or odd run.
In each case, $y'$ has at most one odd run more than~$x'$, implying $o(y')\leq o(x')+1$.
Overall, we have $\beta(y')=w(y')+o(y')\leq (w(x')-1)+(o(x')+1)=w(x')+o(x')=\beta(x')$, as desired.

Secondly, we need to argue that if~$y\in\{0,1\}^n$ satisfies $\beta(y)\leq d$, then $y$ is contained in the downset of some $x\in F_{n,d}$.
If $w(y)=d$ then $\beta(y)\leq d$ implies $o(y)=0$ and therefore $y\in F_{n,d}$.
Otherwise, we have $w(y)<d$.
To proceed, it is enough to prove the following claim:
For any $y'\in\{0,1\}^n\setminus\{1^n\}$, one can flip a 0-bit in~$y'$ to a 1-bit to obtain a bitstring~$x'$ such that $\beta(x')=\beta(y')+1$ if $o(y')=0$ and $\beta(x')\leq \beta(y')$ if $o(y')>0$.
Repeatedly applying such flips to some~$y'\in\{0,1\}^n$ with $\beta(y')\leq d$ and $w(y')<d$ increases the weight and thus leads to a $y$ in the upset of~$y'$ with $\beta(y)\leq d$ and $w(y)=d$, as covered in the first case.
The proof of the claim goes as follows.
If $o(y')=0$, then flipping the first 0-bit in $y'$ gives a bitstring~$x'$ without any odd runs, so $o(x')=o(y')=0$ and $w(x')=w(y')+1$, yielding $\beta(x')=w(x')+o(x')=(w(y')+1)+o(y')=\beta(y')+1$.
On the other hand, if $o(y')>0$, then flipping the 0-bit preceding one of the odd runs in~$y'$ gives a bitstring~$x'$ that has at least one odd run less than~$y'$, so $o(x')\leq o(y')-1$ and $w(x')=w(y')+1$, yielding $\beta(x')=w(x')+o(x')\leq (w(y')+1)+(o(y')-1)=\beta(y')$, as claimed.

This completes the proof of the lemma.
\end{proof}

\subsection{Hamiltonian cycles in the face lattice}

The following is the main result of this section.

\begin{theorem}
\label{thm:cyclic}
For any $n>d\geq 2$, the graph $G(L(C_{n,d}))$ has a Hamiltonian cycle.
\end{theorem}

By Lemma~\ref{lem:cyclic-downset}, our goal is to list the bitstrings in the set~$X_{n,d}$ defined in~\eqref{eq:Xnd} by flipping a single bit in each step (also, the trivial face at the top of the face lattice still needs to be inserted).
We will do this by restricting the binary reflected Gray code in~$Q_n$ to the subset~$X_{n,d}$.
This sublist technique is very powerful and has been exploited in several previous papers (see \cite[Sec.~3.16--3.18+4.1]{MR4649606}).
However, the fact that it works for the set~$X_{n,d}$ is new and somewhat surprising, as in all previous applications it only gives less restrictive Gray codes (allowing up to~2 or~3 flipped bits between consecutive bitstrings, whereas we insist on single flips).

\begin{figure}[h!]
\includegraphics[page=2]{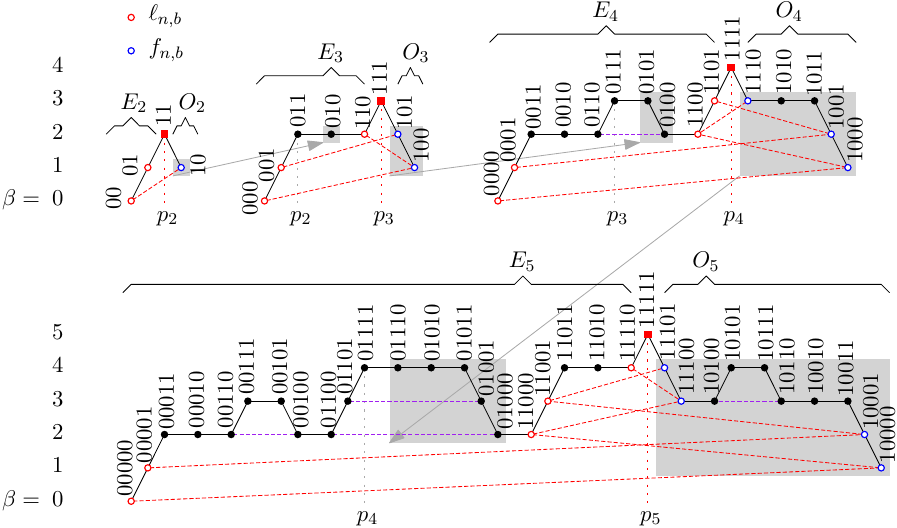}
\caption{Illustration of the binary reflected Gray code, with the vertices aligned vertically according to their~$\beta$-values.
Note that the $\beta$-value changes by at most~1 in every step along the Gray code, as described by Lemma~\ref{lem:brgc-aux}~(iii).
The gray boxes highlight all bitstrings in~$O_n$, located after position~$p_n$ in the listing~$\Gamma_n$ by Lemma~\ref{lem:brgc-aux}~(i).
Their $\beta$-values increase by~$+1$ according to Lemma~\ref{lem:brgc-aux}~(ii) in the induction step, i.e., the content of the gray boxes is shifted up by one unit when prefixing with~0.
The vertices $f_{n,b}$ and~$\ell_{n,b}$, whose explicit form is given by Lemma~\ref{lem:brgc-aux}~(iv), are highlighted.
Dashed lines connect certain pairs of vertices differing in a single bit.
}
\label{fig:brgc}
\end{figure}

The following definitions and properties are illustrated in Figure~\ref{fig:brgc}.
We define
\begin{align*}
E_n &\ass \{x=(x_1,\ldots,x_n)\in\{0,1\}^n\mid \exists k\geq 0 \text{ even s.t.\ } (x_1,\ldots,x_{k+1})=1^k0\}, \text{ and } \\
O_n &\ass \{x=(x_1,\ldots,x_n)\in\{0,1\}^n\mid \exists k\geq 1 \text{ odd s.t.\ }(x_1,\ldots,x_{k+1})=1^k0\},
\end{align*}
i.e., these are the sets of bitstrings of length~$n$ with an even or odd number of leading 1s, respectively.
For example, we have $E_3=\{000,001,011,010,110\}$ and $O_3=\{100,101\}$.
Clearly, $E_n\cup O_n\cup\{1^n\}$ is a partition of $\{0,1\}^n$.

For any $x\in\{0,1\}^{n-1}$ and $a\in\{0,1\}$ we have
\begin{equation}
\label{eq:beta-bx}
\beta(ax)=
\begin{cases}
\beta(x)   & \text{if } a=0 \text{ and } x\in E_{n-1}\cup\{1^{n-1}\}, \\
\beta(x)+1 & \text{if } a=1 \text{ or } x\in O_{n-1}.
\end{cases}
\end{equation}
Note that the condition stated in the second case of~\eqref{eq:beta-bx} is precisely the negation of the condition stated in the first case.

The binary reflected Gray code is defined as $\Gamma_1\ass (0,1)$ and
\begin{equation}
\label{eq:brgc}
\Gamma_n\ass 0\Gamma_{n-1},1\rev(\Gamma_{n-1})
\end{equation}
for $n\geq 1$.
In words, the listing~$\Gamma_n$ consists of two copies of the previous listing~$\Gamma_{n-1}$, where all bitstrings in the first copy are prefixed with~0 and all bitstrings in the second copy are prefixed with~1, and the order of bitstrings in the second copy is reversed.
We let $p_n$ be the position of the string~$1^n$ in~$\Gamma_n$.
For $b=0,\ldots,n-1$ we write $\ell_{n,b}$ for the last bitstring~$x$ in~$\Gamma_n$ before~$1^n$ with $\beta(x)=b$.
Similarly, for $b=1,\ldots,n-1$ we let $f_{n,b}$ denote the first bitstring~$x$ in~$\Gamma_n$ after~$1^n$ with~$\beta(x)=b$.

For any bitstring$~x\in\{0,1\}^n$ we write~$x^-\in\{0,1\}^{n-1}$ for the string obtained from~$x$ by removing the first bit.
The following lemma captures four crucial properties of the binary reflected Gray code that we will exploit for our construction later; see Figure~\ref{fig:brgc}.

\begin{lemma}
\label{lem:brgc-aux}
The binary reflected Gray code $\Gamma_n\assr(x_1,\ldots,x_{2^n})$ has the following properties:
\begin{enumerate}[label=(\roman*),leftmargin=8mm,topsep=1mm]
\item For any $i<p_n$ we have $x_i\in E_n$ and for any $i>p_n$ we have $x_i\in O_n$.
In words, all bitstrings before~$1^n$ have an even number of leading~1s and all bitstrings after~$1^n$ have an odd number of leading~1s.
\item For any bitstring~$x_i$ we have
\[
\beta(x_i)=
\begin{cases}
\beta(x_i^-)   & \text{if } i\leq p_{n-1}, \\
\beta(x_i^-)+1 & \text{if } i>p_{n-1}.
\end{cases}
\]
\item For $i=1,\ldots,2^n-1$ we have $\beta(x_{i+1})-\beta(x_i)\in\{-1,0,+1\}$.
In words, the $\beta$-values of any two consecutive bitstrings in~$\Gamma_n$ differ by at most~1.
\item For $b=0,\ldots,n-1$ we have
\begin{subequations}
\label{eq:ell-f}
\begin{equation}
\ell_{n,b}=
\begin{cases}
1^b0^{n-b} & \text{if $b$ is even}, \\
1^{b-1}0^{n-b}1 & \text{if $b$ is odd},
\end{cases}
\end{equation}
and for $b=1,\ldots,n-1$ we have
\begin{equation}
\label{eq:fnb}
f_{n,b}=
\begin{cases}
1^b0^{n-b} & \text{if $b$ is odd}, \\
1^{b-1}0^{n-b}1 & \text{if $b$ is even}.
\end{cases}
\end{equation}
\end{subequations}
\end{enumerate}
\end{lemma}

We emphasize that in part~(ii) of the lemma, the bitstring~$x_i$ is part of the Gray code~$\Gamma_n$, but the index~$p_{n-1}$ is (by definition) determined with respect to the Gray code~$\Gamma_{n-1}$.
Consequently, by~\eqref{eq:brgc} the bitstring at position~$p_{n-1}$ in~$\Gamma_n$ is the bitstring~$01^{n-1}$.

\begin{proof}
To prove part~(i), we argue by induction on~$n$.
The induction basis $n=1$ is trivial.
For the induction step, consider the definition~\eqref{eq:brgc}.
Clearly, all bitstrings in $0\Gamma_{n-1}$ belong to the set~$E_n$, as they have zero leading 1s, which is an even number.
By induction, we have $\Gamma_{n-1}=E'_{n-1},1^{n-1},O'_{n-1}$, where $E'_{n-1}$ and $O'_{n-1}$ are certain orderings of the elements of $E_{n-1}$ and~$O_{n-1}$, respectively.
Consequently, we have $1\rev(\Gamma_{n-1})=1\rev(O'_{n-1}),1^n,1\rev(E'_{n-1})$, i.e., all bitstrings before~$1^n$ belong to the set~$E_n$, and all bitstrings after~$1^n$ belong to the set~$O_n$.
This proves part~(i).

Part~(ii) follows from~\eqref{eq:brgc}, part~(i) of the lemma (applied to~$\Gamma_{n-1}$) and~\eqref{eq:beta-bx}.

To prove part~(iii), we argue by induction on~$n$.
Using~\eqref{eq:brgc} and part~(ii) settles the claim for all $i=1,\ldots,2^n-1$ except for $i=p_{n-1}$ and $i=2^{n-1}$.
For the transitions at these two indices we continue the computation as follows.
For $i=p_{n-1}$ we obtain $\beta(x_i)=\beta(01^{n-1})=n-1$ and $\beta(x_{i+1})=\beta(0f_{n-1,n-2})=n-1$ with the help of~\eqref{eq:fnb}, so $\beta(x_{i+1})-\beta(x_i)=0$.
For $i=2^{n-1}$ we obtain $\beta(x_i)=\beta(010^{n-2})=2$ and $\beta(x_{i+1})=\beta(110^{n-2})=2$, so $\beta(x_{i+1})-\beta(x_i)=0$.
This proves part~(iii).

To prove~(iv), first observe that the definition~\eqref{eq:brgc} implies $\ell_{n,0}=0^n$, $\ell_{n,1}=0^{n-1}1$ and
\begin{equation*}
\begin{split}
\ell_{n,b}&=1f_{n-1,b-1} \text{  for } b=2,\ldots,n-1; \\
f_{n,b}&=1\ell_{n-1,b-1} \text{  for } b=1,\ldots,n-1.
\end{split}
\end{equation*}
From those recursive relations the explicit formulas follow easily by induction.
\end{proof}

\begin{lemma}
\label{lem:brgc-cyclic}
For any $d\geq 1$ and~$n\geq d+1$, let $\Gamma_{n,d}$ be the sublist of~$\Gamma_n$ induced by the set~$X_{n,d}$ defined in~\eqref{eq:Xnd}, i.e., the list obtained from~$\Gamma_n$ by deleting all bitstrings~$x$ with $\beta(x)>d$.
Then in~$\Gamma_{n,d}$ any two cyclically consecutive bitstrings differ in a single bit, except the consecutive strings~$\ell_{n,d}$ and~$f_{n,d}$, which differ in a transposition $0\leftrightarrow 1$.
\end{lemma}

In Figure~\ref{fig:brgc}, the bitstrings are arranged vertically according to their $\beta$-values.
Consequently, the sublist~$\Gamma_{n,d}$ of~$\Gamma_n$ is obtained by taking all bitstrings from left to right whose height is at most~$d$, or equivalently, by cutting out subsequences of bitstrings at height more than~$d$, taking shortcuts along the dashed edges.
The resulting sublist~$\Gamma_{5,3}$ is highlighted in Figure~\ref{fig:downset}.

As mentioned in Section~\ref{sec:algo}, the sublist~$\Gamma_{n,d}$ of~$\Gamma_n$ can be computed in constant time per vertex, despite the fact that sometimes many intermediate vertices are skipped when taking a shortcut.
This is achieved by establishing and exploiting explicit successor rules between shortcut vertices, analogous to what has been done for similar generation algorithms based on sublists of the binary reflected Gray code; see for example~\cite{MR0349274}.
For details, see our C++ implementation~\cite{cos_cyclic}.

\begin{proof}
Let $\Gamma_n\assr(x_1,\ldots,x_{2^n})$, and consider a maximal subsequence of consecutive indices~$k\in\{1,\ldots,2^n\}$ with $\beta(x_k)>d$.
By Lemma~\ref{lem:brgc-aux}~(iii), the preceding index~$i$ and the following index~$j$ of this subsequence satisfy $\beta(x_i)=\beta(x_j)=d$, i.e., we have $\beta(x_k)>d$ for all $i<k<j$.
If $x_i=\ell_{n,d}$ and~$x_j=f_{n,d}$, then $x_i$ and~$x_j$ differ in a transposition $0\leftrightarrow 1$ by Lemma~\ref{lem:brgc-aux}~(iv).
If not, then we have $d\geq 2$ and we consider two cases.
If $i<j<p_{n-1}$ or $p_{n-1}<i<j$, then the claim follows by induction on~$n$, using Lemma~\ref{lem:brgc-aux}~(ii).
It remains to consider the case $i<p_{n-1}<j$.
The assumption that $\beta(x_k)>d$ for all $i<k<j$, Lemma~\ref{lem:brgc-aux}~(ii), and the fact that $0\Gamma_{n-1}$ ends with $y\ass 010^{n-2}$ which satisfies~$\beta(y)=2$, imply that $x_i=0\ell_{n-1,d}$ and $x_j=0f_{n-1,d-1}$.
Therefore, if $d$ is even we obtain from Lemma~\ref{lem:brgc-aux}~(iv) that $x_i=01^d0^{n-d-1}$ and $x_j=01^{d-1}0^{n-d}$, i.e., $x_i$ and~$x_j$ differ only in the bit at position~$d+1$.
On the other hand, if $d$ is odd, we obtain $x_i=01^{d-1}0^{n-d-1}1$ and $x_j=01^{d-2}0^{n-d}1$, i.e., $x_i$ and~$x_j$ differ only in the bit at position~$d$.
\end{proof}

We are now in position to prove Theorem~\ref{thm:cyclic}.

\begin{proof}[Proof of Theorem~\ref{thm:cyclic}]
Combine Lemmas~\ref{lem:cyclic-downset} and \ref{lem:brgc-cyclic}.
We claim that the two strings~$\ell_{n,d}$ and $f_{n,d}$ that are consecutive in the listing~$\Gamma_{n,d}$ both correspond to facets of~$C_{n,d}$, so we can insert the trivial face~$C_{n,d}$ at the top of the face lattice~$L(C_{n,d})$ between these two in order to obtain a Hamiltonian cycle in~$G(L(C_{n,d}))$; see the bottom part of Figure~\ref{fig:downset}.
Indeed, $x\in\{\ell_{n,d},f_{n,d}\}$ satisfies $\beta(x)=d$ by definition, and $o(x)=0$ by~\eqref{eq:ell-f}, implying that $w(x)=d$.
\end{proof}

\section{3-dimensional polytopes}
\label{sec:3dim}

\subsection{Hamiltonian cycles in the face lattice}

We consider 3-dimensional polytopes via their Schlegel diagrams, which are the 3-connected plane graphs.
We establish the following auxiliary statement, which actually holds even for 2-connected plane graphs.

\begin{lemma}
\label{lem:plane-path}
Let $H=(V,E)$ be a 2-connected plane graph, let $(e_1,v_1,e_2,v_2,\ldots,v_{k-1},e_k,v_k)$ be the sequence of edges and vertices in counterclockwise order along the outer face~$\ol{f}$, and let $f$ be the inner face incident with~$e_k$.
Then $G(L(H))\setminus\{\ol{f},\emptyset,H\}$ has a Hamiltonian path that starts at~$f$, visits each of the pairs~$e_i,v_i$ for $i=1,\ldots,k-1$ consecutively, and ends with the triple~$v_{k-1},e_k,v_k$.
\end{lemma}

\begin{figure}[htb]
\includegraphics[page=2]{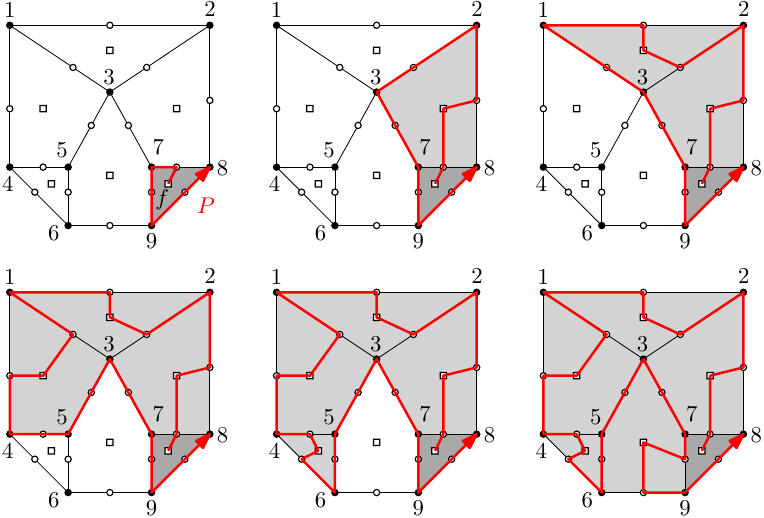}
\caption{Illustration of the path constructed inductively as described in the proof of Lemma~\ref{lem:plane-path} for the polytope~$P$ from Figure~\ref{fig:example}.
The resulting Hamiltonian cycle in~$G(L(P))$ is shown in parts~(d1) and~(d2) of that figure.}
\label{fig:path-example}
\end{figure}

\begin{proof}
This proof is illustrated schematically in Figure~\ref{fig:path}, and for a concrete example in Figure~\ref{fig:path-example}.
The proof uses the notion of the~\defi{weak dual graph}, which is obtained from the dual graph by removing the vertex that corresponds to the outer face~$\ol{f}$.
Note that if a plane graph is 2-connected, then its weak dual graph is connected.

\begin{figure}[htb]
\makebox[0cm]{ 
\includegraphics[page=1]{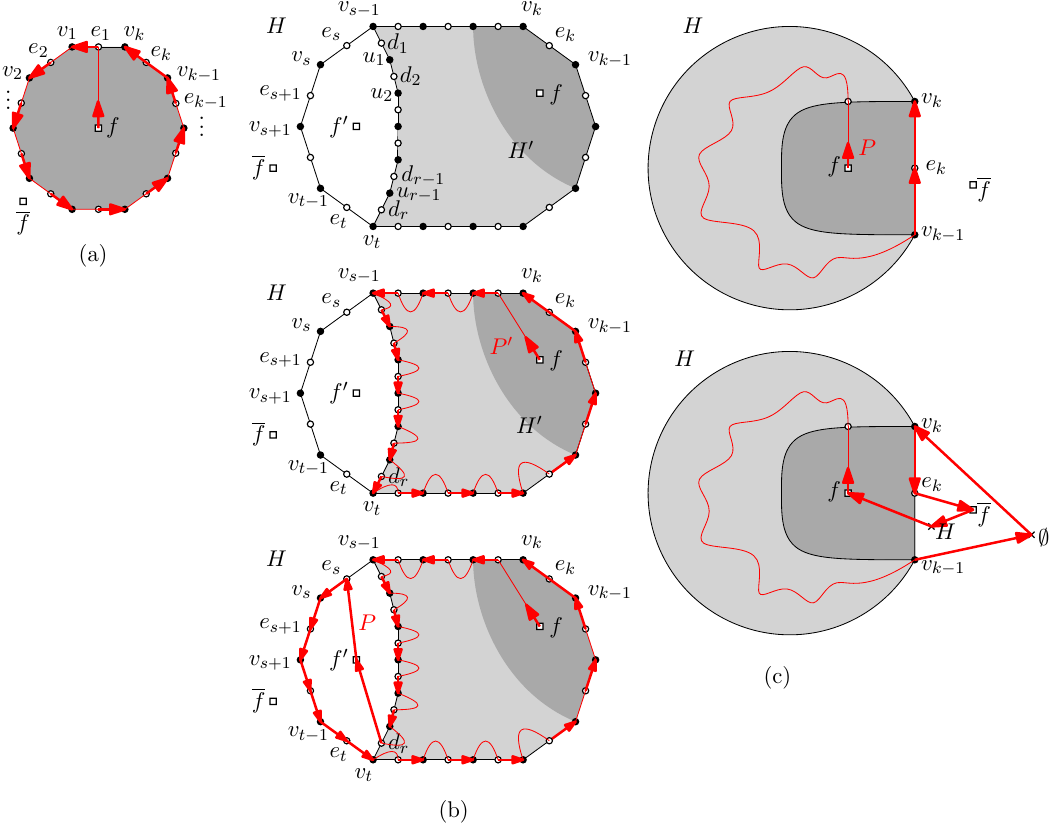}
}
\caption{Illustration of the proofs of Lemma~\ref{lem:plane-path} and Theorem~\ref{thm:3d-HC}.
Vertices, edges and faces are drawn as bullets, circles and squares, respectively.}
\label{fig:path}
\end{figure}

We prove the statement by induction on the number of faces of~$H$.
If $H$ has only one face, then $H$ is a cycle $(e_1,v_1,e_2,v_2,\ldots,v_{k-1},e_k,v_k)$, bounding the only inner face~$f$, and then $P\ass(f,e_1,v_1,e_2,v_2,\ldots,v_{k-1},e_k,v_k)$ is the desired path; see Figure~\ref{fig:path}~(a).

For the induction step, suppose that~$H$ has more than one face, and let $(e_1,v_1,e_2,v_2,\ldots,e_k,v_k)$ and $f$ be as in the lemma.
In the weak dual graph~$D$ of~$H$, we pick a vertex~$f'$ different from~$f$ that is also adjacent to~$\ol{f}$ in the dual graph; see Figure~\ref{fig:path}~(b).
Note that $f'$ is a face in the primal graph~$H$, and this choice ensures that removing the face~$f'$ from~$H$ yields again a 2-connected graph.
Let $(v_{s-1},e_s,v_s,e_{s+1},v_{s+1},\ldots,v_{t-1},e_t,v_t)$ for $1\leq s\leq t<k$ be the sequence of vertices and edges between~$f'$ and the outer face.
Furthermore, let $(d_1,u_1,d_2,u_2,\ldots,d_{r-1},u_{r-1},d_r)$, $r\geq 1$, be the sequence of edges and vertices between~$v_{s-1}$ and~$v_t$ in clockwise order around~$f'$.
Consider the 2-connected graph
\[ H'\ass H\setminus\{e_s,v_s,e_{s+1},v_{s+1},\ldots,v_{t-1},e_t\}, \]
which has the same faces as~$H$ except~$f'$, and which has the sequence
\[
(e_1,v_1,e_2,v_2,\ldots,v_{s-2},e_{s-1},v_{s-1},d_1,u_1,d_2,u_2,\ldots,d_{r-1},u_{r-1},d_r,v_t,e_{t+1},v_{t+1},\ldots,v_{k-1},e_k,v_k) \]
of edges and vertices in counterclockwise order along the outer face.
Moreover, we denote by~$\ol{d}$ the outer face of~$H'$.
By induction, $G(L(H'))\setminus\{\ol{d},\emptyset,H'\}$ has a Hamiltonian path~$P'$ that starts at~$f$, visits the pairs $e_i,v_i$ for $i=1,\ldots,s-1$ and $i=t+1,\ldots,k$ consecutively, as well as the pairs $d_i,u_i$ for $i=1,\ldots,r-1$ and the pair $d_r,v_t$, and that ends with the triple~$v_{k-1},e_k,v_k$.
To obtain the desired Hamiltonian path~$P$ for~$G(L(H))\setminus\{\ol{f}\cup f',\emptyset,H\}$, we insert the sequence
\[ (f',e_s,v_s,e_{s+1},v_{s+1},\ldots,v_{t-1},e_t) \]
between $d_r$ and~$v_t$ in~$P'$.
It can be checked directly that~$P$ has the required properties.

This completes the proof of the lemma.
\end{proof}

\begin{theorem}
\label{thm:3d-HC}
For any 2-connected plane graph~$H$, the graph~$G(L(H))$ has a Hamiltonian cycle.
Consequently, for any 3-dimensional polytope~$P$, the graph~$G(L(P))$ has a Hamiltonian cycle.
\end{theorem}

\begin{proof}
We take the path~$P=(f,\ldots,v_{k-1},e_k,v_k)$ in~$G(L(H))$ guaranteed by Lemma~\ref{lem:plane-path}.
We turn the path into a Hamiltonian cycle of~$G(L(H))$ by adding the outer face~$\ol{f}$ and the trivial cells~$\emptyset$ and~$H$, by replacing the last two entries~$e_k,v_k$ of~$P$ by the sequence~$(\emptyset,v_k,e_k,\ol{f},H)$, which makes the resulting sequence cyclic, as $H$ and~$f$ are adjacent in~$G(L(H))$; see Figure~\ref{fig:path}~(c).
\end{proof}

\subsection{Rhombic strips in the face lattice}

\begin{figure}[htb]
\makebox[0cm]{ 
\includegraphics[page=1]{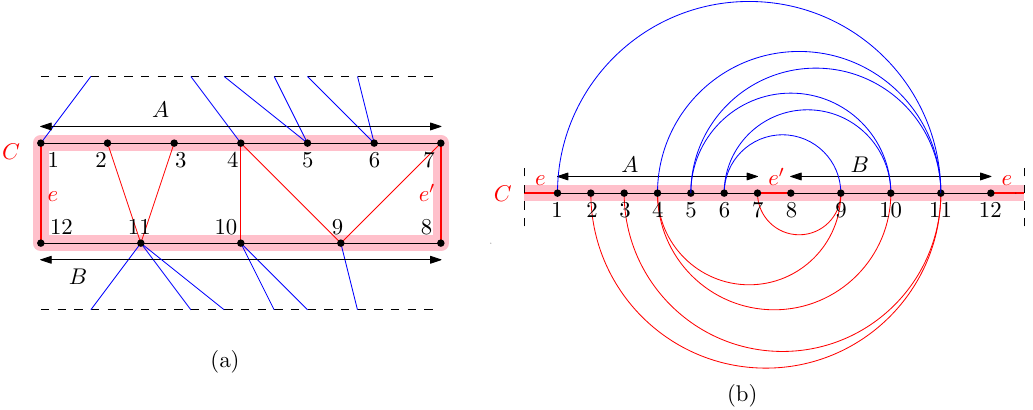}
}
\caption{Two illustrations of the same graph and Hamiltonian cycle that satisfies condition~(i) in Theorem~\ref{thm:3d-RS}.
In~(a), the $C$-chords outside of the cycle~$C$ wrap around the top and bottom boundary (dashed).
In~(b), the cycle~$C$ wraps around the left and right boundary, and the $C$-chords inside and outside of the cycle are shown below and above it, respectively.}
\label{fig:3-RS}
\end{figure}

The next result, illustrated in Figure~\ref{fig:3-RS}, gives a characterization for when the inclusion order~$L(H)$ of the cells of a plane graph~$H$ admits a rhombic strip.
In particular, this characterizes the 3-dimensional polytopes whose face lattice admits a rhombic strip.
Recall from Lemma~\ref{lem:strip-order}~(v) that a necessary condition for~$L(H)$ to have a rhombic strip is that~$H$ and its dual graph both admit a Hamiltonian cycle.
However, it turns out that these necessary conditions are not sufficient.

Given a graph~$H$ and Hamiltonian cycle~$C$ in~$H$, we refer to the edges of~$H\setminus C$ as \defi{$C$-chords}.

\begin{theorem}
\label{thm:3d-RS}
Let $H$ be a 2-connected plane graph.
Then $G(L(H))$ has a rhombic strip if and only if $H$ has a Hamiltonian cycle~$C=(v_1,\ldots,v_n)$ that satisfies one of the following two equivalent conditions:
\begin{enumerate}[label=(\roman*),leftmargin=8mm,topsep=1mm]
\item The cycle~$C$ can be (vertex-)partitioned into two paths~$A$ and~$B$ such that every $C$-chord has one endpoint on~$A$ and the other on~$B$.
\item There are no three $C$-chords $(v_{i_1},v_{i_2}),(v_{i_3},v_{i_4}),(v_{i_5},v_{i_6})$ with $1\leq i_1<i_2\leq i_3<i_4\leq i_5<i_6\leq n+1$, where $v_{n+1}=v_1$.
\end{enumerate}
\end{theorem}

It is easy to check that the existence of a Hamiltonian cycle~$C$ as in Theorem~\ref{thm:3d-RS}~(i) implies that the dual graph of~$H$ also has a Hamiltonian cycle.

\begin{proof}
In the first part of the proof, we show that $G(L(H))$ has a rhombic strip if and only if $H$ has a Hamiltonian cycle~$C$ that satisfies~(i).

`$\Rightarrow$':
We assume that~$L(H)$ has a rhombic strip~$R$.
Let $C$ be the Hamiltonian cycle in~$H$ defined by the cyclic ordering of all vertices of~$H$ in~$R$ (recall Lemma~\ref{lem:strip-order}~(i)+(v)).
By Lemma~\ref{lem:strip-order}~(i), for any edge~$e=(u,v)\in C$, both cover relations~$u\cov e$ and~$v\cov e$ are present as edges in~$R$.
The sequence of vertices and edges along~$C$ and the cover relations between them thus forms a cycle in the rhombic strip~$R$ that alternates between ranks~0 and~1.
By planarity and Lemma~\ref{lem:strip-order}~(ii), for any $C$-chord~$e=(u,v)$, exactly one of the cover relations~$u\cov e$ or~$v\cov e$ is present as an edge in~$R$.
Let $f$ and~$f'$ be the two faces incident with~$e$ in~$H$.
By Lemma~\ref{lem:strip-order}~(iii), both cover relations~$e\cov f$ and~$e\cov f'$ must be present as edges in~$R$, and therefore $f$ and~$f'$ appear consecutively in the cyclic ordering of all faces of~$H$ in~$R$ (recall Lemma~\ref{lem:strip-order}~(i)), i.e., there is a diamond~$(f,e,f',H)$ in~$R$.
As this observation applies to every $C$-chord~$e$, we obtain that no face of~$H$ is bounded by more than two $C$-chords, and thus the cyclic ordering of all faces of~$H$ given by~$R$ can be partitioned into two contiguous parts~$A'$ and~$B'$, where~$A'$ consists of the faces inside of~$C$ and $B'$ of the faces outside of~$C$, such that any two consecutive faces in~$A'$ and~$B'$ cover the one $C$-chord that separates them, and the two pairs of boundary entries of~$A'$ and~$B'$ each cover one edge of~$C$ that separates them.
Removing these two separating edges from~$C$  (keeping their end vertices) yields the partition into paths~$A$ and~$B$ with the desired properties.

\begin{figure}[htb]
\includegraphics[page=2]{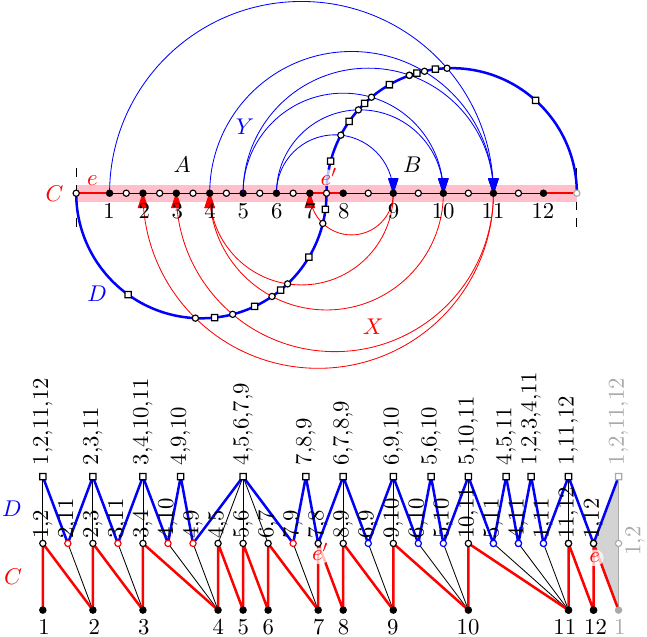}
\caption{Illustration of the proof of Theorem~\ref{thm:3d-RS}.
Every $C$-chord~$c$ has an arrow pointing to the vertex~$v$ such $(v,c)$ is an edge in the rhombic strip.
The top part shows the same graph as Figure~\ref{fig:3-RS}, and the bottom part shows the resulting rhombic strip.
Vertices, edges and faces are drawn as bullets, circles and squares, respectively.
For clarity, the trivial cells at the bottom and top are not shown in the rhombic strip.}
\label{fig:proof-3-RS}
\end{figure}

\begin{figure}[htb]
\includegraphics[page=3]{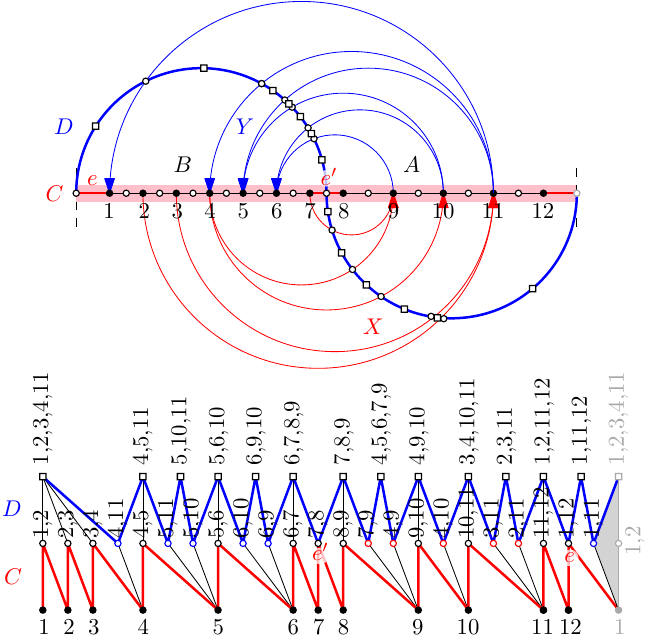}
\caption{A second rhombic strip constructed for the graph~$H$ from Figure~\ref{fig:proof-3-RS} by swapping the roles of~$A$ and~$B$.}
\label{fig:swap}
\end{figure}

`$\Leftarrow$':
The notations used in this part of the proof are illustrated in Figure~\ref{fig:proof-3-RS}.
We assume that~$H$ has a Hamiltonian cycle~$C$ satisfying condition~(i).
Let $e$ and~$e'$ be the two edges connecting the paths~$A$ and~$B$ to the cycle~$C$.

We specify the rhombic strip of~$L(H)$ by defining a cyclic ordering of the vertices, edges and faces of~$H$, which determines the order in which they are embedded on the sphere, at three different heights according to their rank.
Furthermore, we specify which cover relations of~$L(H)$ are added as (straight) edges to the rhombic strip, ensuring that there are no edge crossings and that all faces are diamonds.

Specifically, vertices of~$H$ and edges on~$C$ are ordered according to the cycle~$C$, and all cover relations between them are added to the rhombic strip, yielding a cycle that alternates between ranks~0 and~1 (the red cycle in Figure~\ref{fig:proof-3-RS}).
Faces of~$H$ and edges of~$H\setminus C$ are ordered according to the cycle~$D$ in the dual graph given by the duals of the edges in~$H\setminus C\cup\{e,e'\}$.
All cover relations between them are added to the rhombic strip, yielding a cycle that alternates between ranks~1 and~2 (the blue cycle in Figure~\ref{fig:proof-3-RS}).

The orderings of the edges of~$C\setminus\{e,e'\}$ and~$H\setminus C$ are interleaved as follows:
Let $X$ and~$Y$ denote the sets of $C$-chords inside or outside of~$C$, respectively.

Any edge~$\wh{e}\in A$ is inserted between the two edges of~$X\cup\{e,e'\}$ that bound the same face~$f$ inside of~$C$, whereas any edge~$\wh{e}\in B$ is inserted between the two edges of~$Y\cup\{e,e'\}$ that bound the same face~$f$ outside of~$C$.
In both cases, the cover relation $\wh{e}\lessdot f$ is added as an edge to the rhombic strip.
In the example of Figure~\ref{fig:proof-3-RS}, the edge~$\wh{e}=(4,5)\in A$ is inserted between the two edges $(4,9),(7,9)\in X\cup\{e,e'\}$, which all bound the face~$f=(4,5,6,7,9)$, so the cover relation $(4,5)=\wh{e}\lessdot f=(4,5,6,7,9)$ is added.
The edge~$\wh{e}=(10,11)\in B$ is inserted between the two edges $(5,10),(5,11)\in Y\cup\{e,e'\}$, which all bound the face~$f=(5,10,11)$, so the cover relation $(10,11)=\wh{e}\lessdot f=(5,10,11)$ is added.

Any edge~$\wh{e}\in X$ is inserted between the two edges of~$A\cup\{e,e'\}$ that have the same endpoint~$v$ as~$\wh{e}$, whereas any edge~$\wh{e}\in Y$ is inserted between the two edges of~$B\cup\{e,e'\}$ that have the same endpoint~$v$ as~$\wh{e}$.
In both cases, the cover relation $v\lessdot \wh{e}$ is added as an edge to the rhombic strip.
In the example of Figure~\ref{fig:proof-3-RS}, the edge~$\wh{e}=(7,9)\in X$ is inserted between the two edges~$(6,7),(7,8)\in A\cup\{e,e'\}$, which all share the endpoint~$v=7$, so the cover relation $7=v\lessdot\wh{e}=(7,9)$ is added.
The edge~$\wh{e}=(1,11)\in Y$ is inserted between the two edges~$(10,11),(11,12)\in B\cup\{e,e'\}$, which all share the endpoint~$v=11$, so the cover relation $11=v\lessdot \wh{e}=(1,11)$ is added.

In the second part of the proof, we show the equivalence between conditions~(i) and~(ii).

(i)$\Rightarrow$(ii):
Suppose that the cycle~$C$ satisfies~(i), and let $A$ and~$B$ be the paths that partition~$C$ as stated in~(i).
Every $C$-chord requires that one of its endpoints belongs to~$A$ and the other to~$B$.
Consequently, if there were three $C$-chords $(v_{i_1},v_{i_2}),(v_{i_3},v_{i_4}),(v_{i_5},v_{i_6})$ with $1\leq i_1<i_2\leq i_3<i_4\leq i_5<i_6\leq n+1$, then the membership of vertices with respect to the paths~$A$ and~$B$ along the cycle~$C$ would alternate at least three times, which is impossible.

(i)$\Leftarrow$(ii):
Suppose that the cycle~$C$ satisfies~(ii).
Let $X$ and~$Y$ denote the sets of all $C$-chords inside or outside of~$C$, respectively.
For $Z\in\{X,Y\}$ there are two disjoint subpaths~$A_Z$ and~$B_Z$ of~$C$ such that each $C$-chord from~$Z$ has one endpoint on~$A_Z$ and the other on~$B_Z$, and the pair of first vertices~$(f(A_Z),f(B_Z))$ and the pair of last vertices~$(\ell(A_Z),\ell(B_Z))$ of~$A_Z$ and~$B_Z$ are each connected by such a $C$-chord from~$Z$.
We distinguish two cases:

(1)~One of the two paths~$A_X,B_X$ has a nonempty intersection with one of the paths~$A_Y,B_Y$.
Condition~(ii) rules out that one of these four paths has a nonempty intersection with both paths from the other pair.
Consequently, the desired partition of~$C$ is given by extending the union of the two paths that have a nonempty intersection with the union of the remaining two paths into maximal disjoint paths~$A$ and~$B$.

(2)~All four paths~$A_X,B_X,A_Y,B_Y$ are disjoint.
In this case the desired partition of~$C$ is given by extending $A_X\cup A_Y$ and~$B_X\cup B_Y$ into maximal disjoint paths~$A$ and~$B$.

This completes the proof of the theorem.
\end{proof}

By swapping the roles of the paths~$A$ and~$B$ in condition~(i) of Theorem~\ref{thm:3d-RS}, we see that any such Hamiltonian cycle~$C$ in~$H$ actually gives rise to two different rhombic strips for~$G(L(H))$; see Figure~\ref{fig:swap}.

\begin{figure}[h!]
\makebox[0cm]{ 
\includegraphics[page=1]{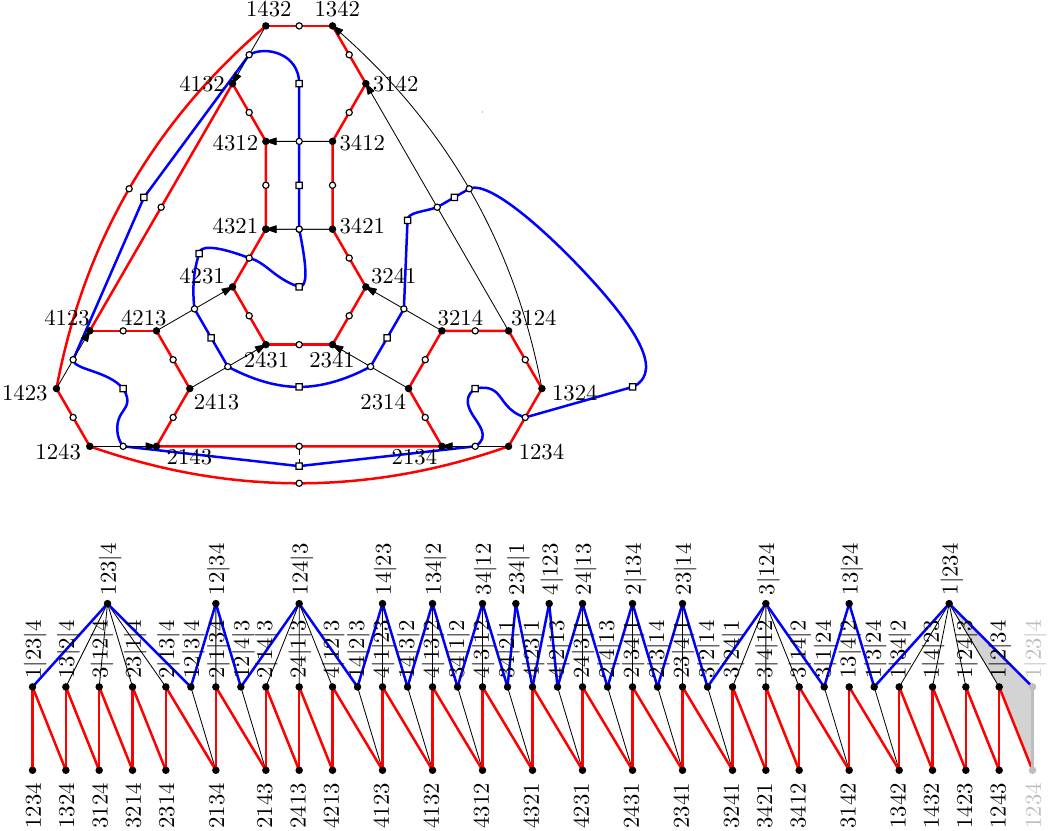}
}
\caption{Rhombic strip for the 3-dimensional permutahedron.
For clarity, the trivial cells at the bottom and top of the face lattice are not shown.
Also, the bars between any two entries of the number partitions at the bottom are omitted.}
\label{fig:RS-Perm3}
\end{figure}

\begin{figure}[p]
\makebox[0cm]{ 
\includegraphics[page=2]{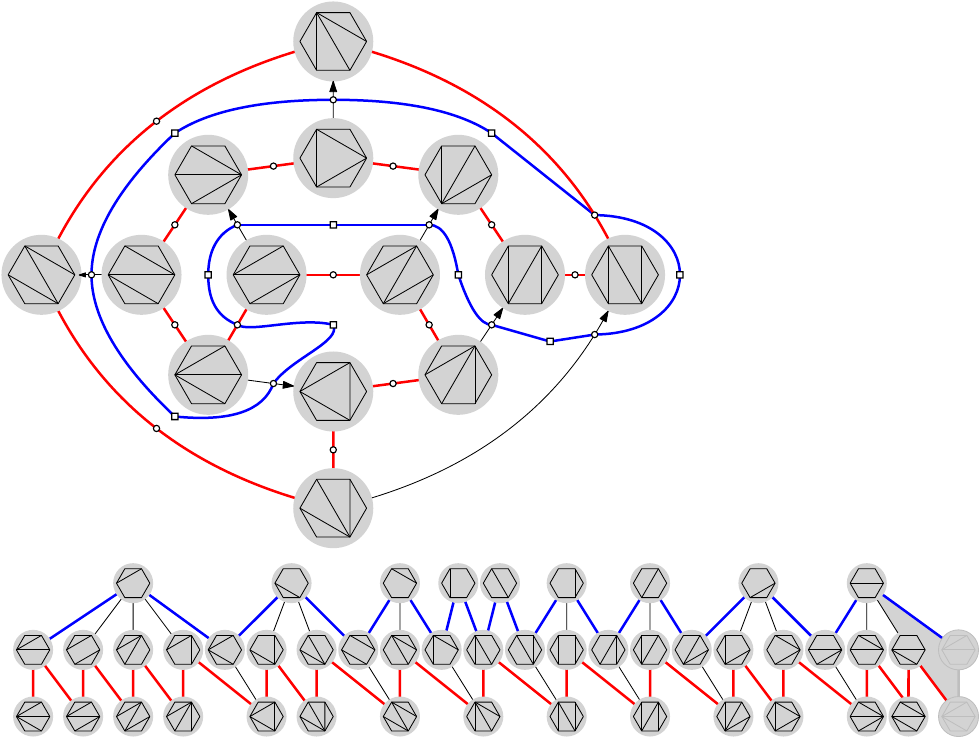}
}
\caption{Rhombic strip for the 3-dimensional associahedron.}
\label{fig:RS-Asso3}
\end{figure}

The 3-dimensional permutahedron and associahedron admit Hamiltonian cycles satisfying the conditions of Theorem~\ref{thm:3d-RS}, so we immediately obtain a rhombic strip for each of them; see Figures~\ref{fig:RS-Perm3} and~\ref{fig:RS-Asso3}, respectively.

\begin{figure}[p]
\includegraphics[page=1]{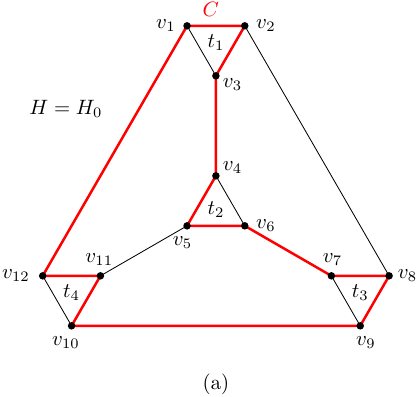}
\includegraphics[page=2]{trunc-tetra}
\caption{(a) Skeleton~$H_0$ of the truncated tetrahedron with one of its Hamiltonian cycles, and (b) the graph~$H_s$ of iterated truncations at the vertex~$v_5$, with a Hamiltonian cycle in its dual graph.
}
\label{fig:trunc-tetra1}
\end{figure}

Using Theorem~\ref{thm:3d-RS}~(ii), we now construct an infinite family of planar graphs~$H$, such that $H$ and its dual graph both admit a Hamiltonian cycle, but $G(L(H))$ does not admit a rhombic strip.
The \defi{truncated tetrahedron} is the polytope obtained from the tetrahedron by truncating every vertex to a triangle.
Its skeleton is the graph~$H$ on 12 vertices depicted in Figure~\ref{fig:trunc-tetra1}~(a).

We turn~$H$ into an infinite family of plane graphs~$H_s$, $s\geq 0$, by repeatedly truncating a vertex.
Specifically, for $s\geq 0$, the graph $H_s$ is obtained from~$H$ by subdividing the edge $(v_4,v_5)$ by $s$ additional vertices $a_1,a_2,\ldots,a_s$ and the edge $(v_6,v_5)$ by $s$ vertices $b_1,b_2,\ldots,b_s$, and by adding the edges $(a_i,b_i)$ for $i=1,\ldots,s$; see Figure~\ref{fig:trunc-tetra1}~(b).
Equivalently, $H_s$ is obtained from~$H_{s-1}$ by truncating the vertex~$v_5$.
Clearly, $H_s$ is a 2-connected cubic plane graph.

\begin{theorem}
\label{thm:3d-RS-counter}
For every $s\geq 0$, the 2-connected plane graph~$H_s$ and its dual graph both have a Hamiltonian cycle, but $G(L(H_s))$ does not admit a rhombic strip.
\end{theorem}

\begin{proof}
One can check directly that~$H_s$ and its dual graph both have a Hamiltonian cycle; see Figure~\ref{fig:trunc-tetra2} and~\ref{fig:trunc-tetra1}~(b), respectively.
In fact, up to symmetry, $H_0$ has a unique Hamiltonian cycle~$C$, depicted in Figure~\ref{fig:trunc-tetra1}~(a).
Specifically, such a cycle contains exactly two consecutive edges from each of the triangles~$t_1,t_2,t_3,t_4$.
It follows that, up to symmetry, $H_s$ has the two distinct Hamiltonian cycles~$C_1$ and~$C_2$ shown in Figure~\ref{fig:trunc-tetra2}.

\begin{figure}
\includegraphics[page=3]{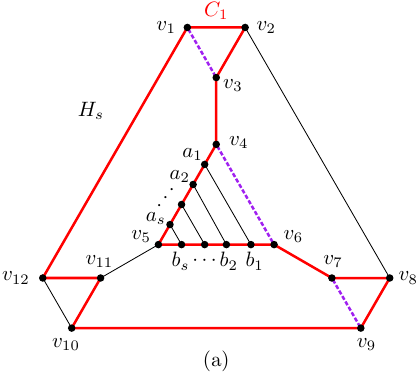}
\includegraphics[page=4]{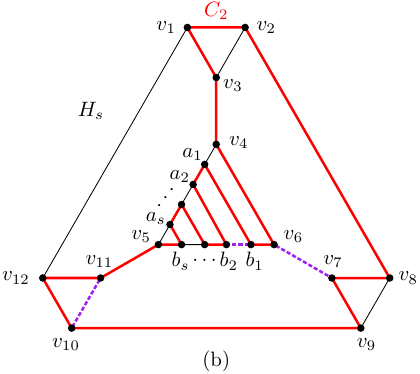}
\caption{The two Hamiltonian cycles~$C_1$ and~$C_2$ in the graph~$H_s$.
The $C_i$-chords violating condition~(ii) in Theorem~\ref{thm:3d-RS} are dashed.}
\label{fig:trunc-tetra2}
\end{figure}

We proceed to show that~$C_1$ and~$C_2$ both violate condition~(ii) of Theorem~\ref{thm:3d-RS}, and consequently $G(L(H_s))$ does not admit a rhombic strip.
For the cycle~$C_1$ shown in Figure~\ref{fig:trunc-tetra2}~(a), the $C_1$-chords~$(v_1,v_3)$, $(v_4,v_6)$ and~$(v_7,v_9)$ violate this condition.
For the cycle~$C_2$ shown in Figure~\ref{fig:trunc-tetra2}~(b), we can instead take the $C_2$-chords~$(v_{11},v_{10})$, $(v_7,v_6)$ and~$(b_1,x)$ with $x\ass v_5$ if $s=1$ and $x\ass b_2$ if $s\geq 2$.
\end{proof}

\begin{figure}[htb]
\makebox[0cm]{ 
\begin{tabular}{ccc}
\includegraphics[page=1]{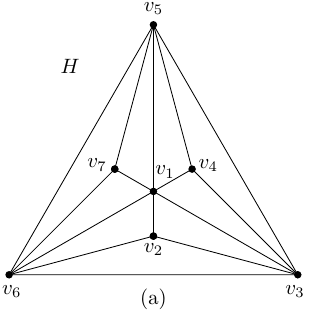} &
\includegraphics[page=2]{fano} &
\includegraphics[page=3]{fano}
\end{tabular}
}
\caption{The smallest plane graph~$H$ such that both~$H$ and its dual graph are Hamiltonian, but $G(L(H))$ does not admit a rhombic strip.
The $C_i$-chords violating condition~(ii) in Theorem~\ref{thm:3d-RS} are dashed.
}
\label{fig:fano}
\end{figure}

With the help of a computer, we determined the smallest 3-connected plane graph~$H$ such that both $H$ and its dual graph have a Hamiltonian cycle, but $G(L(H))$ does not admit a rhombic strip; see Figure~\ref{fig:fano}~(a).
The graph~$H$ has 7 vertices and is the incidence graph of the Fano plane.
Also, there are no 2-connected graphs with 6~vertices with those properties.

\begin{proposition}
\label{prop:fano}
The graph~$H$ and its dual graph both have a Hamiltonian cycle, but $G(L(H))$ does not admit a rhombic strip.
\end{proposition}

\begin{proof}
Up to symmetry, $H$ has two distinct Hamiltonian cycles, $C_1=(v_1, v_2,v_3,v_4,v_5,v_6,v_7)$ and $C_2=(v_1,v_6,v_2,v_3,v_4,v_5,v_7)$; see Figure~\ref{fig:fano}~(b)+(c).
A Hamiltonian cycle in the dual graph of~$H$ is shown in Figure~\ref{fig:fano}~(b).
For the cycle~$C_1$, the $C_1$-chords $(v_1,v_6),(v_6,v_3),(v_3,v_1)$ violate condition~(ii) in Theorem~\ref{thm:3d-RS}.
Similarly, for the cycle~$C_2$, the $C_2$-chords $(v_3,v_5),(v_5,v_6),(v_6,v_3)$ violate this condition.
Consequently, $G(L(H))$ does not admit a rhombic strip.
\end{proof}

\section{Graph associahedra of chordal graphs}
\label{sec:Gasso}

\subsection{The polytope and its face lattice}

Let $H=(V,E)$ be a graph.
A \defi{tube $T$} of~$H$ is a nonempty subset of vertices~$T\subset V$, such that the induced subgraph~$H[T]$ is connected.
Two tubes~$T$ and~$T'$ are \defi{compatible} if they are either nested, i.e., $T\subset T'$ or~$T\supset T'$, or non-adjacent, i.e., $H[T\cup T']$ is not connected.
A \defi{tubing} of~$H$ is a family of pairwise compatible tubes that includes the vertex sets of each connected component of~$H$.
The \defi{graph associahedron $A(H)$} of~$H$ is the polytope whose face lattice is isomorphic to the reverse inclusion order of all tubings of~$H$.
In particular, the vertices of $A(H)$ are the inclusion-maximal tubings of~$H$ (where the number of tubes equals the number of vertices of~$H$), and they are in one-to-one correspondence with elimination forests of~$H$.

\subsection{Hamiltonian cycles in the face lattice}

The following is the main result of this section.
A graph is \defi{chordal} if it has no induced cycles of length at least~4.

\begin{theorem}
\label{thm:LGAsso-HC}
For any chordal graph~$H$ with at least one edge, the graph $G(L(A(H)))$ has a Hamiltonian cycle.
\end{theorem}

If $H$ is a perfect matching with $n$ edges, then $A(H)$ is the hypercube~$Q_n$.
If $H$ is a complete graph on $n$ vertices, then $A(H)$ is the permutahedron~$\Pi_n$.
If $H$ is a path on $n$ vertices, then $A(H)$ is the (standard) associahedron~$A_{n+2}$.
All these graphs are chordal, and therefore Theorem~\ref{thm:LGAsso-HC} generalizes Theorems~\ref{thm:LQ-HC}, \ref{thm:LPi-HC} and~\ref{thm:LAsso-HC} presented before.

For a vertex~$v$ of a graph~$H$ we write $H-v$ for the graph obtained by deleting~$v$ (and all incident edges) from~$H$.
A \defi{clique} in~$H$ is an induced subgraph of~$H$ that is complete.
A well-known characterization of chordal graphs that we will be using is that they admit a perfect elimination order (PEO).
This is a total ordering of the vertices of~$H$ such that every vertex induces a clique with the neighboring vertices that come before it in the ordering.
Formally, a \defi{PEO graph}~$H=([n],E)$ is one that satisfies one of the following two recursive conditions:
$n=1$, i.e., $H$ is a single vertex graph; or $n>1$, the graph~$H-n$ is a PEO graph and the induced subgraph of the vertex~$n$ and its neighbors is a clique in~$H$.

The following lemma is an immediate consequence of the definition of compatible tubes.

\begin{lemma}
\label{lem:tubing-clique}
Let $H$ be a graph, let $C$ be a clique in~$H$, let $\cT$ be a tubing of~$H$, and let $T,T'\in \cT$ be two distinct tubes with $C\cap T\neq \emptyset$ and $C\cap T'\neq \emptyset$.
Then $T$ and~$T'$ are nested, i.e., we have $T\subset T'$ or $T\supset T'$.
\end{lemma}

\begin{proof}[Proof of Theorem~\ref{thm:LGAsso-HC}]
We consider $H$ as a PEO graph~$H=([n],E)$, and we assume w.l.o.g.\ that if it has exactly one edge and~$n>2$, then the vertex~$n$ is not incident to this edge.

The following definitions are illustrated in Figure~\ref{fig:tubing}.
Let $\cT$ be a tubing of~$H-n$, and let $T_1,T_2,\dots,T_k$ be the tubes of~$\cT$ containing at least one neighbor of~$n$ in $H$.
By Lemma~\ref{lem:tubing-clique} they are nested, so we can assume that $T_1\subset T_2\subset\cdots\subset T_k$.
Furthermore, let $\cR\ass\cT\setminus\{T_1,\dots,T_k\}$ be the remaining set of tubes that do not contain any neighbor of~$n$.

If $n$ is an isolated vertex of~$H$, then we have $\cR=\cT$ and $k=0$, and we define $\wc{c}_0(\cT)\ass \cR\cup\{\{n\}\}=\cT\cup\{\{n\}\}$.
Otherwise, if $n$ is not an isolated vertex, we define
\begin{align*}
\wc{c}_i(\cT)&\ass\cR\cup\bigcup_{j=1}^i \{T_j\} \cup \bigcup_{j=i}^k \big\{T_j\cup\{n\}\big\}, \quad i=0,\ldots,k, \\
\wh{c}_i(\cT)&\ass\cR\cup\bigcup_{j=1}^{i-1} \{T_j\} \cup \bigcup_{j=i}^k \big\{T_j\cup\{n\}\big\}, \quad i=1,\ldots,k,
\end{align*}
where $T_0\cup\{n\}$ in the equation for~$\wc{c}_0(\cT)$ should be interpreted as~$\{n\}$.

\begin{figure}[htb]
\includegraphics{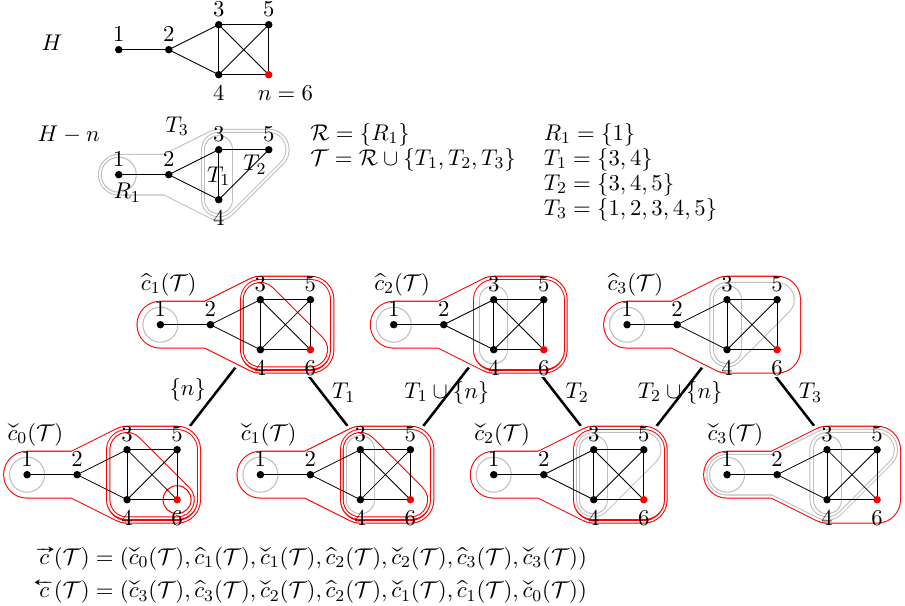}
\caption{Illustration of the proof of Theorem~\ref{thm:LGAsso-HC}.
The vertex~$n$ and all tubes containing it are highlighted.
The bold zigzag path shows edges in~$G(L(A(H)))$, labeled by the tube that is added/removed along this edge.
}
\label{fig:tubing}
\end{figure}

Note that both $\wc{c}_i(\cT)$ and $\wh{c}_i(\cT)$ are tubings of~$H$, where the latter has the same number of tubes as~$\cT$, and the former has one more tube than~$\cT$.
We define the sequence $\rvec{c}(\cT)$ of tubings of~$H$ by
\[ \rvec{c}(\cT)\ass\big(\wc{c}_0(\cT),\wh{c}_1(\cT),\wc{c}_1(\cT),\wh{c}_2(\cT),\wc{c}_2(\cT),\ldots,\wh{c}_{k-1}(\cT),\wc{c}_{k-1}(\cT),\wh{c}_k(\cT),\wc{c}_k(\cT)\big), \]
we note that it has length~$2k+1$ and describes a path in~$G(L(A(H)))$ that alternates between ranks~$r$ and~$r+1$, where $r$ is the rank of~$\cT$ in~$L(A(H-n))$; see Figure~\ref{fig:tubing}.
We write $\lvec{c}(X)\ass\rev(\lvec{c}(X))$ for the reverse sequence/path.

To prove the theorem, we construct a path~$P_n$ in $G(L(A(H)))$ that visits all faces of~$A(H)$ except~$\emptyset$ and that starts and ends at rank~0 faces (i.e., two inclusion-maximal tubings of~$H$ with $n$ tubes each), and therefore $(P_n,\emptyset)$ is the desired Hamiltonian cycle in~$G(L(A(H)))$.

The path~$P_n$ is constructed by induction on~$n$ as follows: For the base case of the construction where $H$ is a single edge connecting vertices~1 and~2 we take $P_2\ass\{\{1\},\{1,2\}\}, \{\{1,2\}\}, \{\{1,2\},\{2\}\}$.

For the induction step, let $P_n\assr(\cT_1,\ldots,\cT_N)$ be the path in~$G(L(A(H-n)))$.
Then we define
\[ P_{n+1}\ass(\lvec{c}(\cT_1),\rvec{c}(\cT_2),\lvec{c}(\cT_3),\rvec{c}(\cT_4),\ldots,\lvec{c}(\cT_N)). \]
It can be checked straightforwardly that the path~$P_{n+1}$ has the required properties.
The main observation is that for every tubing~$\cT$ of~$H$, there is a uniquely defined tubing~$\cT'$ of~$H-n$ and a unique integer~$i$ such that $\cT=\wc{c}_i(\cT')$ or~$\cT=\wh{c}_i(\cT')$.
\end{proof}

\section{Quotientopes}
\label{sec:quotient}

\subsection{The polytopes and Hamiltonian cycles in the face lattice}

We start by extending some of the terminology and notation introduced in Section~\ref{sec:prelim}.

Given a lattice~$(P,<)$, a \defi{lattice congruence} is an equivalence relation $\equiv$ on $P$ that preserves the join and meet operations, i.e., if $x\equiv x'$ and $y\equiv y'$, then $x\vee y\equiv x'\vee y'$ and $x\wedge y\equiv x'\wedge y'$.
For any $x\in P$ we write $[x]\ass\{y\in P\mid x\equiv y\}$ for the equivalence class of~$x$, and we define the set of all equivalence classes as $P/{\equiv}\ass\{[x]\mid x\in P\}$.
It is well-known and easy to show that every equivalence class is an interval.
This set inherits a lattice structure, where for any $X,Y\in P/{\equiv}$, we define $X<Y$ if there exist elements~$x\in X$ and~$y\in Y$ such that $x<y$.
We refer to $P/{\equiv}$ with this lattice structure as \defi{quotient lattice}.

\begin{figure}
\centering
\includegraphics[page=1]{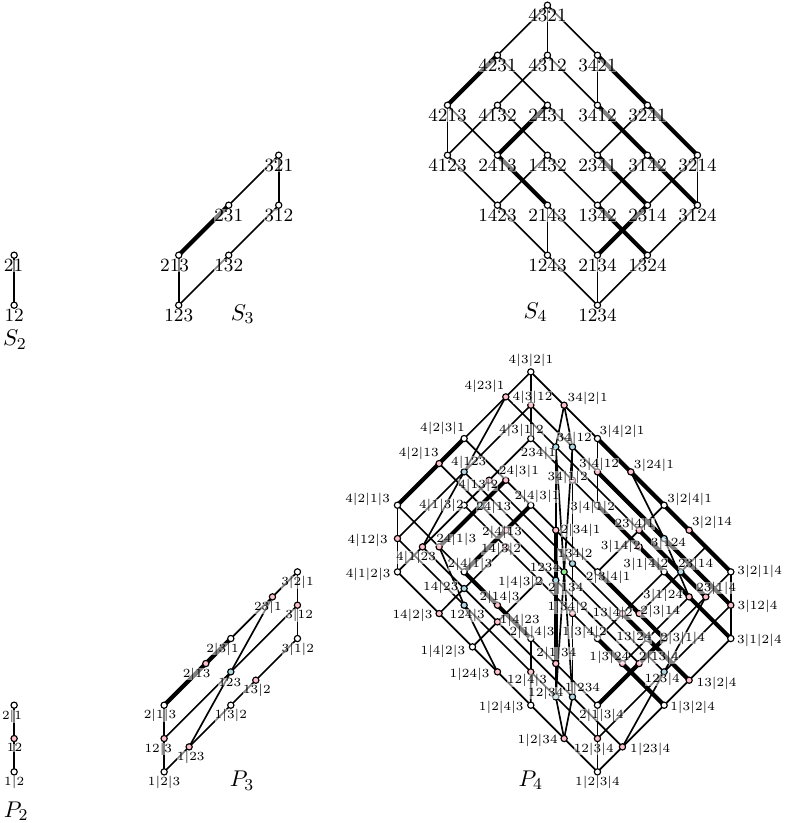}
\caption{Top: The weak order on~$S_n$ for $n=2,3,4$, together with a lattice congruence, namely the sylvester congruence.
The bold edges indicate pairs of permutations that are equivalent.
Bottom: The facial weak order on~$P_n$ for $n=2,3,4$, with the extensions of the lattice congruences from the top.
Set partitions corresponding to vertices, edges, 2- and 3-dimensional faces are indicated by white, red, blue and green bullets, respectively.
}
\label{fig:cong}
\end{figure}

In this paper we focus on lattice congruences of the so-called weak order on permutations.
Formally, the \defi{inversion set} of a permutation~$\pi=a_1\cdots a_n\in S_n$ is defined as
\[ \inv(\pi)\ass\{(a_i,a_j)\mid 1\leq i<j\leq n\mid a_i>a_j\}, \]
and the \defi{weak order} on~$S_n$ is defined by $\pi<\rho$ if and only if $\inv(\pi)\subset\inv(\rho)$; see the top part of Figure~\ref{fig:cong}.
The weak order is a lattice, namely, given two permutations~$\pi$ and~$\rho$, their meet~$\sigma\ass\pi\wedge \rho$ is given by $\inv(\sigma)=\inv(\pi)\cap \inv(\rho)$, and their join~$\sigma'\ass\pi\vee\rho$ has as its inversion set~$\inv(\sigma')$ the transitive closure of~$\inv(\pi)\cup\inv(\rho)$.
Cover relations in the weak order are adjacent transpositions, so the cover graph of the weak order on~$S_n$ is exactly the skeleton of the permutahedron~$\Pi_n$.
More generally, Pilaud and Santos~\cite{MR3964495} proved that for any lattice congruence~$\equiv$ of the weak order, the cover graph of $S_n/{\equiv}$ can be realized geometrically as the skeleton of a polytope, called \defi{quotientope} (see also \cite{MR4584712}), which we denote by~$\Pi_n/{\equiv}$.
Via suitable congruences, quotientopes generalize permutahedra, associahedra and hypercubes, and many other combinatorial polytopes, including permutreehedra~\cite{MR3856522} and rectangulotopes~\cite{MR4833064}.

The following is the main result of this section.

\begin{theorem}
\label{thm:LSn-HC}
For any $n\geq 1$ and any lattice congruence~$\equiv$ of the weak order on~$S_n$, the graph $G(L(\Pi_n/{\equiv}))$ has a Hamiltonian cycle, unless it is a single edge.
\end{theorem}

The case that $G(L(\Pi_n/{\equiv}))$ is a single edge occurs if and only if $\equiv$ has only one equivalence class (in particular, if $n=1$).

\begin{figure}[t!]
\makebox[0cm]{ 
\includegraphics[page=2,scale=0.9]{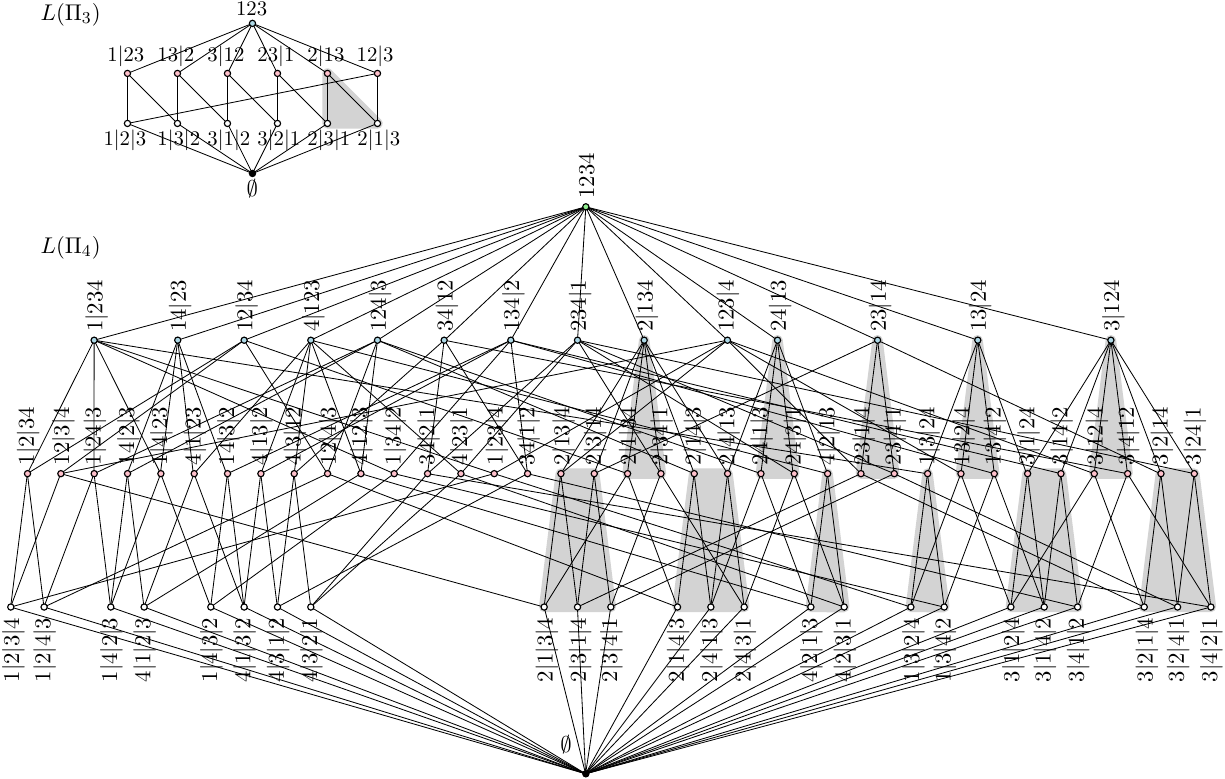}
}
\caption{The face lattices of the 2- and 3-dimensional permutahedron, namely the refinement orders on set partitions of~$[n]$, for $n=3$ and $n=4$, respectively, plus the trivial face~$\emptyset$ at the bottom.
The gray bubbles are the equivalence classes of the lattice congruences shown in Figure~\ref{fig:cong}.}
\label{fig:LPi-cong}
\end{figure}

\subsection{The face lattice of the permutahedron}

Recall from Section~\ref{sec:perm} that the face lattice~$L(\Pi_n)$ of the permutahedron~$\Pi_n$ can be modelled as the refinement order on set partitions of~$[n]$.
Specifically, we denote all set partitions of~$[n]$ by~$P_n$, and we write $\subset$ for the refinement order on~$P_n$, which corresponds to the inclusion order on the faces of~$\Pi_n$; see Figure~\ref{fig:LPi-cong}.
I.e., $x\cov y$ for $x,y\in P_n$ if and only if $x$ is obtained from~$y$ by splitting a block into two/adding one bar, or equivalently, $y$ is obtained from~$x$ by taking the unions of two blocks/removing one bar.
For example, we have $2|14|3|5\cov 124|3|5\cov 124|35$ and $2|14|3|5\not\subset 14|2|35$.
Set partitions with $n$ blocks can be identified naturally with permutations in~$S_n$.

\subsection{The facial weak order}

Finding an `explicit' combinatorial model of the faces of quotientopes is an open problem in general.
For the interesting special case of simple quotientopes, whose cover graphs are regular, the faces can be modeled by Schr\"oder separating trees \cite{MR5050522}.
\TM{We should spell out in terms of these trees (at least for once concrete example) what the zigzag algorithm is doing (in this special cases).}
In our paper we use a representation of the faces as equivalence classes of a lattice congruence of an order on the set~$P_n$ that is known as \defi{facial weak order}~\cite{krob_et_al:2011,MR3729508,MR4356247}, which is different from the refinement order discussed in the previous section.

To define the facial weak order, we need some definitions; see the bottom part of Figure~\ref{fig:cong}.
We say that two blocks $A,B\seq[n]$ of a set partition are \defi{separated} if $\max A<\min B$ or $\min A>\max B$.
In the first case we write $A<B$ and in the second case $A>B$.
If $A<B$, then we refer to $A$ as the \defi{smaller} block and to $B$ as the \defi{larger} block.
To define the facial weak order, we start with the weak order on~$S_n$, i.e., permutations ordered by containment of their inversion sets.
In terms of set partitions, the cover relations are $\pi\lessdot \rho$ for partitions $\pi=a_1|\cdots|a_i|a_{i+1}|\cdots|a_n\in S_n\seq P_n$ and $\rho=a_1|\cdots|a_{i-1}|a_{i+1}|a_i|a_{i+2}|\cdots|a_n\in S_n\seq P_n$ where $a_i<a_{i+1}$.
We extend this order to the entire set~$P_n$ by repeatedly adding, for every pair of set partitions $x=A_1|\cdots|A_i|A_{i+1}|\cdots|A_k\in P_n$ and $y=A_1|\cdots|A_{i+1}|A_i|\cdots|A_k\in P_n$ that have a separated pair of consecutive blocks $A_i$ and~$A_{i+1}$ such that $A_i<A_{i+1}$, the comparabilities $x<z<y$ for $z\ass A_1|\cdots|A_i\cup A_{i+1}|\cdots|A_k$.
In the first round, all set partitions with $n-1$ blocks are added (they encode edges of the permutahedron), and more generally, in the $r$th round all set partitions with $n-r$ blocks are added (they encode $r$-dimensional faces of the permutahedron).
Note that the cover relations in the facial weak order also correspond to adding/removing a bar in a set partition, but the order is different from the refinement order discussed in the previous section.
Furthermore, while in the face lattice of the permutahedron, every addition/removal of a bar is a cover relation, in the facial weak order, a bar can only be added to a block so as to create two separated smaller blocks, and we can only remove a bar between two separated blocks.
In the first case, we go down in the order when the smaller block is placed left of the added bar and the larger block is placed right of the bar, and we go up when placing them in the opposite order.
In the second case, we go up in the order when the block left of the removed bar is smaller than the block right of the bar, and we go down when they are in the opposite order.
For example, we have $12|34\lessdot 1234\lessdot 34|12$ in the facial weak order, and $12|34,34|12\cov 1234$ in the refinement order.
On the other hand, while $x\ass 13|24\cov 1234\assr y$ is a cover relation in the refinement order, $x$ and~$y$ are incomparable in the facial weak order.
In other words, the cover graph of the facial weak order is a (proper) spanning subgraph of the cover graph of the refinement order.
The facial weak order is known to be a lattice~\cite{krob_et_al:2011}.

The key observation is that every face of the permutahedron~$\Pi_n$ corresponds to an interval in the weak order on~$S_n$, more precisely, that the set of vertices (permutations) belonging to this face forms an interval in the weak order.
Namely, for a face of~$\Pi_n$ represented by the set partition~$x=A_1|\cdots|A_k\in P_n$, we let $x_0\in S_n$ and $x_1\in S_n$ be the permutations obtained from~$x$ by sorting the entries of each block~$A_i$, $i\in[k]$, increasingly or decreasingly, respectively.
Then the face $x$ of~$\Pi_n$ contains precisely the permutations of the interval~$[x_0,x_1]$ in the weak order on~$S_n$.

Using this observation, every lattice congruence of the weak order on~$S_n$ can be extended in a natural way to a lattice congruence of the facial weak order on~$P_n$.
Specifically, given an equivalence relation~$\equiv$ on the set~$S_n$, we define an equivalence relation~$\equiv^+$ on the set~$P_n$ as follows.
For two set partitions~$x,y\in P_n$ (representing faces of~$\Pi_n$), we consider the intervals~$[x_0,x_1]$ and~$[y_0,y_1]$ in the weak order on~$S_n$, and we define $x\equiv^+ y$ if and only if $x_0\equiv y_0$ and $x_1\equiv y_1$.

The following lemma provides an alternative, but equivalent, definition of this notion of extension.

\begin{lemma}
For any two set partitions $x,y\in P_n$ we have $x\equiv^+ y$ if and only if for every $x'\in [x_0,x_1]$ there exists a $y'\in [y_0,y_1]$ such that $x'\equiv y'$, and for every $y'\in[y_0,y_1]$ there exists an $x'\in[x_0,x_1]$ such that $x'\equiv y'$.
\end{lemma}

\begin{proof}
Let $x,y\in P_n$ be such that $x\equiv^+ y$, i.e., we have $x_0\equiv y_0$ and $x_1\equiv y_1$.
Furthermore, consider an arbitrary element $x'\in[x_0,x_1]$.
Since $x_0\le x'\leq x_1$, we clearly have $x'=x_0\vee x'$ and $x'=x'\wedge x_1$.
By the definition of lattice congruence, $\equiv$ preserves join and meet operations, and thus we have
\[
x' = x_0 \vee x' = x_0 \vee (x'\wedge x_1) \equiv y_0 \vee (x'\wedge y_1)\eqqcolon y'.
\]
Note that $y_0\le y_1$ and $x'\wedge y_1\le y_1$, implying that $y'\le y_1$.
Furthermore, we have $y_0\le y'$ since $y'$ is the join of $y_0$ and some other element.
Combining these observations proves that $y'\in [y_0,y_1]$, as desired.
The same argument also works with reversing the roles for~$x'$ and~$y'$, so one direction of the equivalence is proved.

To prove the other direction, we assume that for every $x'\in [x_0,x_1]$ there exists a $y'\in [y_0,y_1]$ such that $x'\equiv y'$, and vice versa.
Let $x'\in[x_0,x_1]$ and $y'\in[y_0,y_1]$ be such that $x'\equiv y_0$ $y'\equiv x_0$.
Using that $\equiv$ preserves the meet operation, we get $x_0=x_0\wedge x'\equiv y'\wedge y_0=y_0$, yielding $x_0\equiv y_0$.
An analogous argument with meets replaced by joins shows that $x_1\equiv y_1$, proving the other direction of the equivalence.
\end{proof}

\begin{theorem}[{\cite[Cor.~4.13+Thm.~4.22]{MR3729508}}]
\label{thm:cong-extends}
For any $n\geq 1$ and any lattice congruence~$\equiv$ of the weak order on~$S_n$, its extension~$\equiv^+$ is a lattice congruence of the facial weak order on~$P_n$.
The set $P_n/{\equiv^+}$ is in bijection with the set of faces of the quotientope~$\Pi_n/{\equiv}$ (excluding the trivial face~$\emptyset$).
\end{theorem}

Note, however, that not every lattice congruence of the facial weak order is an extension of a lattice congruence of the weak order.

In view of Theorem~\ref{thm:cong-extends}, our goal is to list equivalence classes of lattice congruences of the facial weak order on~$P_n$.
In fact, the following Hamiltonicity result is valid for \emph{any} lattice congruence~$\equiv$ of the facial weak order on~$P_n$ (in particular for those obtained by extending a congruence of the weak order on~$S_n$).
We refer to a path in the cover graph of~$P_n/{\equiv}$ as \defi{grounded}, if it starts and ends with an equivalence class that contains a permutation.

\begin{theorem}
\label{thm:fweak-cong}
For any $n\geq 1$ and any lattice congruence~$\equiv$ of the facial weak order on~$P_n$, the cover graph of the lattice quotient $P_n/{\equiv}$ has a Hamiltonian path.
If the congruence extends a lattice congruence of the weak order on~$S_n$, then the path is grounded.
\end{theorem}

\subsection{Proofs of Theorems~\ref{thm:fweak-cong} and~\ref{thm:LSn-HC}}

\begin{figure}[b!]
\makebox[0cm]{ 
\includegraphics{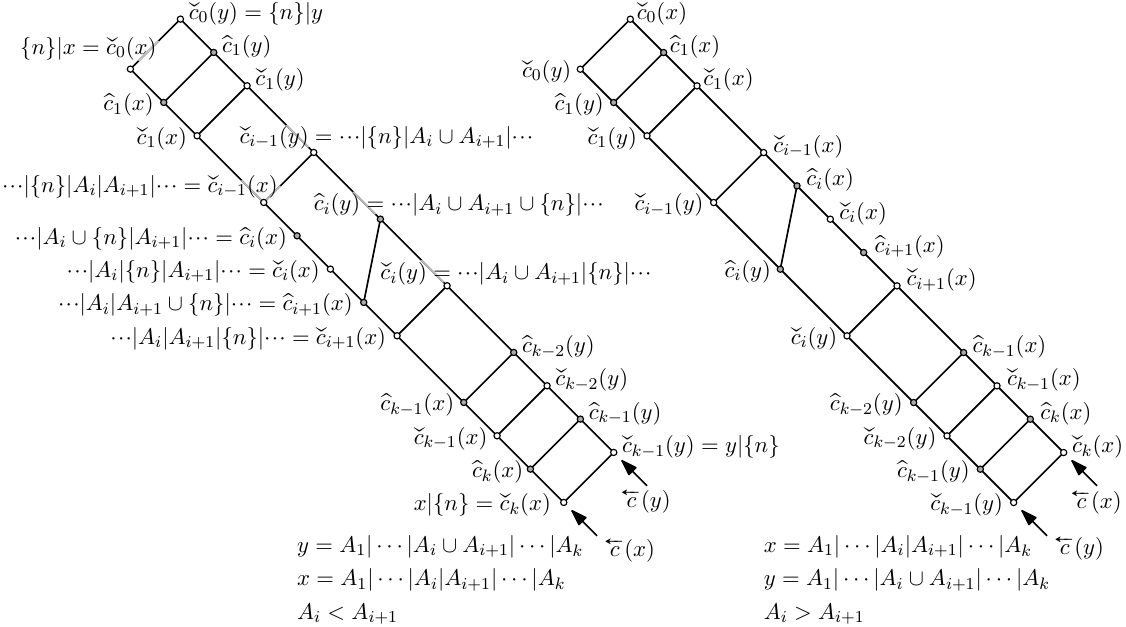}
}
\caption{Recursive structure of the facial weak order.}
\label{fig:forcing}
\end{figure}

The facial weak order on~$P_n$ has a simple recursive structure.
Specifically, for a set partition~$x=A_1|\cdots|A_k\in P_{n-1}$, we define
\begin{align*}
\wc{c}_i(x)&\ass A_1|\cdots |A_i |\{n\}| A_{i+1}|\cdots| A_k  &\text{for } i=0,\dots,k, \\
\wh{c}_i(x)&\ass A_1|\cdots |A_{i-1}|A_i\cup\{n\}|A_{i+1}|\cdots|A_k &\text{for } i=1,\dots,k,
\end{align*}
and
\[
\rvec{c}(x)\ass\big(\wc{c}_0(x),\wh{c}_1(x),\wc{c}_1(x),\wh{c}_2(x),\wc{c}_2(x),\ldots,\wh{c}_{k-1}(x),\wc{c}_{k-1}(x),\wh{c}_k(x),\wc{c}_k(x)\big)
\]
and $\lvec{c}(x)\ass\rev(\rvec{c}(x))$.
The sequence $\lvec{c}(x)$ is a chain in the facial weak order on~$P_n$.
Furthermore, consider set partitions $x,y\in P_{n-1}$ with $x=A_1|\cdots|A_k$ and $y=A_1|\cdots|A_i\cup A_{i+1}|\cdots|A_k$.
Then the two chains $\lvec{c}(x)$ and~$\lvec{c}(y)$ have exactly the following `ladder-like' cover relations; see Figure~\ref{fig:forcing}.
If $x<y$, i.e., if $A_i<A_{i+1}$, then we have
\begin{subequations}
\label{eq:ladder}
\begin{equation}
\begin{alignedat}{2}
\wc{c}_j(x)&\lessdot \wc{c}_j(y) \text{ for } j=0,\ldots,i-1, & & \text{   and   }
\wc{c}_j(x)\lessdot \wc{c}_{j-1}(y) \text{ for } j=i+1,\ldots,k, \\
\wh{c}_j(x)&\lessdot \wh{c}_j(y) \text{ for } j=1,\ldots,i-1, & & \text{   and   }
\wh{c}_j(x)\lessdot \wh{c}_{j-1}(y) \text{ for } j=i+1,\ldots,k.
\end{alignedat}
\end{equation}
On the other hand, if $y<x$, i.e., if $A_i>A_{i+1}$, then we have
\begin{equation}
\begin{alignedat}{2}
\wc{c}_j(y)&\lessdot \wc{c}_j(x) \text{ for } j=0,\ldots,i-1, & & \text{   and   }
\wc{c}_{j-1}(y)\lessdot \wc{c}_j(x) \text{ for } j=i+1,\ldots,k, \\
\wh{c}_j(y)&\lessdot \wh{c}_j(x) \text{ for } j=1,\ldots,i, & & \text{   and   }
\wh{c}_{j-1}(y)\lessdot \wh{c}_j(x) \text{ for } j=i+2,\ldots,k.
\end{alignedat}
\end{equation}
\end{subequations}
For any set partition~$x\in P_n$, we write $p(x)$ for the partition obtained from~$x$ by removing the element~$n$ from the block that contains it.
If $x$ has $\{n\}$ as a singleton block, then the block is removed entirely.
Clearly, we have $p(\wc{c}_i(x))=x$ and $p(\wh{c}_i(x))=x$.

\begin{proposition}
\label{prop:quot}
For any $n\geq 1$ and any lattice congruence~$\equiv$ of the facial weak order on~$P_n$, we have the following:
\begin{enumerate}[label=(\roman*),leftmargin=8mm,topsep=1mm]
\item The restrictions ${\equiv^*}\ass\{(x,y)\mid x,y\in P_{n-1}\text{ and }x|\{n\}\equiv y|\{n\}\}$ and ${\equiv'}\ass\{(\pi,\rho)\mid \pi,\rho\in S_{n-1}\text{ and }\pi|\{n\}\equiv \rho|\{n\}\}$ are lattice congruences of the facial weak order on~$P_{n-1}$ and the weak order on~$S_{n-1}$, respectively, and the former extends the latter.
\item For any equivalence class~$X\in P_n/{\equiv}$, the projection $p(X)\ass\{p(x)\mid x\in X\}$ is an equivalence class of~$\equiv^*$.
\item For any $x=A_1|\cdots|A_k\in P_{n-1}$ we define $\lvec{c}(x)\assr(c_1(x),\ldots,c_{2k+1}(x))$ and $\lambda(x)\ass|\{[c_i(x)]\mid i=1,\ldots,2k+1\}|$, which is the number of equivalence classes intersected by the chain~$\lvec{c}(x)$.
The quantity $\lambda(x)$ is either equal to~$1$ for all $x\in P_{n-1}$ or strictly larger than~1 for all $x\in P_{n-1}$.
\item If the cover graph of $P_{n-1}/{\equiv^*}$ has a grounded Hamiltonian path~$H^*$, then the cover graph of $P_n/{\equiv}$ has grounded a Hamiltonian path~$H$.
\end{enumerate}
\end{proposition}

\begin{proof}
The proofs of~(i)--(iii) are straightforward calculations using the basic definitions of lattice congruences, exploiting the ladder structure of the facial weak order captured by~\eqref{eq:ladder} and shown in Figure~\ref{fig:forcing}; for details, see \cite[Prop.~16~(i)-(iii)]{MR5042274}.

It remains to prove~(iv).
If $\lambda(x)=1$ for all $x\in P_{n-1}$, then $P_n/{\equiv}$ and $P_{n-1}/{\equiv^*}$ are isomorphic, so the claim is trivial.
Otherwise, by (iii) we have $\lambda(x)\geq 2$ for all $x\in P_{n-1}$, i.e., each of the chains $\lvec{c}(x)$ intersects at least two equivalence classes in~$P_n/{\equiv}$.
Given a Hamiltonian path~$H^*=([x_1],\ldots,[x_N])$ of the cover graph of~$P_{n-1}/{\equiv^*}$, where $x_i\in P_{n-1}$ and $N=|P_{n-1}/{\equiv^*}|$, we construct a Hamiltonian path~$H$ of the cover graph of~$P_n/{\equiv}$ as follows; see Figure~\ref{fig:fweak-zigzag}:
For each $x_i$, $i\in[N]$, the chain $\lvec{c}(x_i)$ intersects~$\lambda(x_i)\geq 2$ many equivalence classes of the congruence~$\equiv$, and by~(ii) those equivalence classes are independent of the representative~$x_i$ of the class~$[x_i]$ of the congruence~$\equiv^*$.
Thus we let $\lvec{c}'(x_i)$ be a subsequence of~$\lvec{c}(x_i)$, chosen so that it contains exactly one representative from each of the intersected equivalence classes of~$\equiv$.
Since $\lambda(x_i)\geq 2$ we may choose the subsequence $\lvec{c}'(x_i)$ so that it contains the first and last entry of~$\lvec{c}(x_i)$, namely the set partitions $x_i|\{n\}$ and $\{n\}|x_i$.
We also define $\rvec{c}'(x_i)\ass\rev(\lvec{c}'(x_i))$.
We then define
\[ H'\ass (\lvec{c}'(x_1),\rvec{c}'(x_2),\lvec{c}'(x_3),\rvec{c}'(x_4),\ldots) \]
as a sequence of representatives of all equivalence classes in~$P_n/{\equiv}$, and we let $H$ be the corresponding sequence of equivalence classes (apply $[\;]$ to each entry of~$H'$).
By construction, $H$ is a Hamiltonian path in the cover graph of~$P_n/{\equiv}$.

It remains to show that $H$ is grounded, using the assumption that $H^*$ is grounded.
We can assume that $x_1$ and~$x_N$ are permutations, i.e., $x_1,x_N\in S_{n-1}$.
The first entry of~$H'$ is the first entry of $\lvec{c}'(x_1)$, which is $x_1|\{n\}\in S_n$, so the first entry of~$H$ indeed contains a permutation.
Furthermore, the last entry of $H'$ is the first or last entry of $\lvec{c}'(x_N)$, which is $x_N|\{n\}\in S_n$ or $\{n\}|x_N\in S_n$, respectively (the first case occurs if $N$ is even and the second case occurs if $N$ is odd), so the last entry of~$H$ indeed contains a permutation.

\begin{figure}
\centering
\includegraphics[page=2]{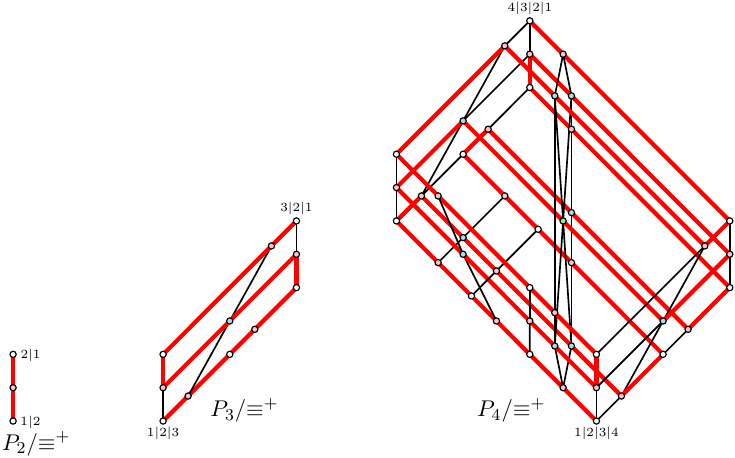}
\caption{(Grounded) Hamiltonian paths in the quotient lattices of the facial weak order for the lattice congruences from Figure~\ref{fig:cong}, constructed as described in the proof of Proposition~\ref{prop:quot}~(iv).
}
\label{fig:fweak-zigzag}
\end{figure}

This completes the proof.
\end{proof}

\begin{proof}[Proof of Theorem~\ref{thm:fweak-cong}]
The statement follows by induction on $n\geq 1$, using Proposition~\ref{prop:quot}~(i)+(iv).
\end{proof}

\begin{proof}[Proof of Theorem~\ref{thm:LSn-HC}]
Let $\equiv$ be any lattice congruence of the weak order on~$S_n$, and let $\equiv^+$ be its extension to a lattice congruence of the facial weak order on~$P_n$.
By Theorem~\ref{thm:cong-extends}, the equivalence classes of~$P_n/{\equiv^+}$ are in bijection with the faces of the quotientopes~$\Pi_n/{\equiv}$, so the Hamiltonian path $H$ in the cover graph of~$P_n/{\equiv^+}$ obtained from Theorem~\ref{thm:fweak-cong} visits every face of the quotientope~$\Pi_n/{\equiv}$ exactly once.
Furthermore, since cover relations in the facial weak order on~$P_n$ are also cover relations in the refinement order on~$P_n$, which is the face lattice of~$\Pi_n$ minus the trivial face~$\emptyset$, we obtain a listing of the faces of~$\Pi_n/{\equiv}$ along cover relations of $G(L(\Pi_n/{\equiv}))$.
Furthermore, as the path~$H$ is grounded, the listing starts and ends at a rank~0 face, i.e., a permutation.
Adding the trivial face~$\emptyset$ thus yields the desired Hamiltonian cycle in~$G(L(\Pi_n/{\equiv}))$.
\end{proof}

\section{Outlook and open questions}
\label{sec:open}

\begin{itemize}[leftmargin=4mm]
\item
Our proof of Conjecture~\ref{conj:LP-HC} for 3-dimensional polytopes (Theorem~\ref{thm:3d-HC}) essentially uses \defi{shellings} (see Lemma~\ref{lem:plane-path}), i.e., a particular order to build the polytope by adding the facets one after the other, while maintaining some invariant on a Hamiltonian cycle in the face lattice of the polytopal complex constructed so far.
This suggests a general strategy for tackling Conjecture~\ref{conj:LP-HC}, namely to combine a particularly nice shelling order with a particular invariant for the Hamiltonian cycle.
Unfortunately, so far we were not able to push this idea beyond the 3-dimensional case, mainly due to our limited knowledge about shellings.
In view of our results regarding special classes of polytopes, worthwhile next targets might be 4-dimensional polytopes, 0/1-polytopes (for example, the uniform matroid polytope), or graph associahedra of non-chordal graphs.

\item
Using the techniques developed in~\cite{MR5042274}, the methods presented in Section~\ref{sec:quotient} can be applied more generally to quotientopes arising from any supersolvable hyperplane arrangement (such as the braid arrangement that gives rise to the permutahedron).
This proves in particular that the face lattice of any of the type~$B$ quotientopes presented in~\cite{MR4584712} admits a Hamiltonian cycle.
This includes type~$B$ associahedra, which are also known as cyclohedra, i.e., graph associahedra of cycles.

\item
A necessary condition for the existence of a Hamiltonian cycle in a bipartite graph is the existence of a perfect matching.
Is it true that the cover graph of the face lattice of every polytope admits a perfect matching?
If yes, this can be seen as further evidence for Conjecture~\ref{conj:LP-HC}.
If no, then this would be a counterexample to the conjecture.

\item
The face lattice of the 3-dimensional permutahedron~$\Pi_4$ has a rhombic strip; see Figure~\ref{fig:RS-Perm3}.
Does the cover graph of the face lattice of the permutahedron~$G(L(\Pi_n))$ admit a rhombic strip for~$n\geq 5$?

\item
The face lattice of the 3-dimensional associahedron~$A_6$ has a rhombic strip; see Figure~\ref{fig:RS-Asso3}.
Does the cover graph of the face lattice of the associahedron~$G(L(A_n))$ admit a rhombic strip for~$n\geq 7$?

Huemer, Hurtado, Noy and Oma{\~{n}}a-Pulido~\cite{MR2510231} constructed a Hamiltonian cycle in a supergraph~$G^+_n$ of~$G(L(A_n))$, obtained by adding edges between pairs of dissections at rank~$k$ that differ in removing an inner edge and replacing it by another edge inside the subpolygon that is the union of the two polygons on both sides of the removed edge.
We note that in~$G(L(A_n))$ such pairs of dissections at rank~$k$ always have a common neighbor at rank~$k+1$, but not necessarily a common neighbor at rank~$k-1$.
Furthermore, not all pairs of dissections at rank~$k$ in~$G(L(A_n))$ that have a common neighbor on rank~$k+1$ are connected by an edge in~$G^+_n$.
The authors also considered the subgraphs of~$G^+_n$ obtained by restricting to fixed rank~$k$, and provided Hamiltonian cycles for them, but this does not solve our problem.

\item
As a partial step towards the previous two questions about the permutahedron and associahedron: Are there listings of the faces of fixed rank~$k$, such that any two consecutive faces have a common superface at rank~$k+1$ and a common subface at rank~$k-1$?
\end{itemize}

\section*{Acknowledgements}

We thank the reviewers of this paper for numerous helpful comments and suggestions that greatly improved the presentation.
In particular, we thank Vincent Pilaud for suggesting the facial weak order as a convenient tool for studying face lattices of quotientopes, and for his enthusiasm to share his knowledge and ideas on how to attack Conjecture~\ref{conj:LP-HC} with us.
Lastly, we thank Christoph Hertrich for suggesting cyclic polytopes as an interesting family of polytopes to investigate.

\bibliographystyle{alpha}
\bibliography{refs}

\end{document}